\documentclass[12pt]{amsart}

\newtheorem{theorem}{Theorem}[section]
\newtheorem{lemma}[theorem]{Lemma}

\newtheorem{proposition}[theorem]{Proposition}
\newtheorem{definition}[theorem]{Definition}

\usepackage[bookmarks]{hyperref}
\usepackage{cite}
\usepackage[makeroom]{cancel}
\usepackage{graphicx,amscd,amsthm,amsfonts,amsopn,amssymb,verbatim,enumerate,xcolor}

\usepackage[letterpaper,left=1in,top=1in,right=1in,bottom=1in]{geometry}
\usepackage{xcolor}

\def\nint{\mathop{\diagup\kern-13.0pt\int}}

\def\supp{{\operatorname{supp}}}

\def\bas{\begin{align*}}
\def\eas{\end{align*}}
\def\bi{\begin{itemize}}
\def\ei{\end{itemize}}

\def\emph#1{{\it #1}}

\DeclareMathOperator{\Tr}{Tr}
\DeclareMathOperator{\curl}{curl}
\theoremstyle{definition}
\newtheorem{remark}[theorem]{Remark}
\newtheorem{thm}[theorem]{Theorem}

\numberwithin{equation}{section}

\title[Strichartz and gallery waves for vacuum irrotational compressible Euler]{Counterexamples to Strichartz estimates for general irrotational compressible Euler backgrounds and positive estimates for a spherical isotropic background solution}
\author{Ovidiu-Neculai Avadanei}
\address{Department of Mathematics, California Institute of Technology, 1200 E.\ California Boulevard, Mail Code 253-37, Pasadena, CA 91125, USA}
\email{avadanei@caltech.edu}
\begin{document}
\begin{abstract}
We consider the free boundary problem for the irrotational compressible Euler equation in a physical vacuum setting. By using the irrotationality condition in the Eulerian formulation of Ifrim and Tataru, we previously showed that the problem can be rephrased in terms of the velocity potential function as an acoustic wave equation that turns out to be widely used in solar seismology. This paper is another step towards understanding what Strichartz estimates are achievable for the aforementioned equation. Our object of study is the corresponding linearized problem. In a previous work, we constructed a class of gallery modes in a model stationary case, which turned out to be counterexamples to Strichartz estimates, forcing a loss of derivatives. By a careful regularization of the domain defining function, we transplant these counterexamples to general Euler backgrounds, and obtain general counterexamples in the same range. On the other hand, we also prove that for the background corresponding to a concrete affine isotropic spherical solution for the irrotational compressible Euler equations, one has Strichartz estimates up to the aforementioned threshold, which shows that they are optimal for at least some explicit Euler backgrounds. To the best of our knowledge, these are the first results of this kind for explicit solutions of the irrotational compressible Euler equations in a physical vacuum.
\end{abstract}
\maketitle
\tableofcontents

\section{Introduction}

In this paper we are going to consider the free boundary problem for the compressible Euler equations describing a gas surrounded by a vacuum. The domain occupied by the gas at time $t$ is denoted by $\Omega_t$, with boundary $\Gamma_t$. These equations can be stated as follows:
\begin{equation}\label{rawEuler}
\begin{cases}\begin{aligned}
  &\rho_t+\nabla\cdot(\rho v)=0\\
  &\rho(v_t+(v\cdot\nabla)v)+\nabla p=0
  \end{aligned}
  \end{cases}
\end{equation}
As we have seen in \cite{Avadanei}, in the irrotational setting, the velocity field can be expressed in terms of a potential which evolves according to an acoustic wave equation, and the same holds true for the linearized counterpart of the compressible Euler system. An interesting question is to determine to which extent one can have Strichartz estimates for the linearized equations, as well as for the fully nonlinear problem. In \cite{Avadanei} we investigated this question in a model case, and it turned out that one has both Strichartz estimates and counterexamples within the respective ranges predicated by scaling considerations. One of the key ideas was to notice that the acoustic metric determines a geometry in which the geodesics reflect along the boundary multiple times, which allowed us to construct solutions that forced derivative losses in the Strichartz estimates.

In what follows, we shall assume that the pressure is given by constitutive laws of the form 
\begin{equation}p(\rho)=\rho^{\kappa+1}\text{, with }\kappa>0.
\end{equation}

In general, we can regard the Euler equations as a system consisting of a wave equation which governs the evolution of the variables $(\rho,\nabla\cdot v)$,  coupled with a transport equations that models the evolution of the vorticity $\omega=\curl v$. Since our gas is assumed to be irrotational, we have $\omega=0$ and $v=\nabla\phi$ for some potential $\phi$, which, as we know from \cite{Avadanei}, solves a wave equation, while its linearized counterpart solves an acoustic wave equation. In both of these equations, an important quantity is the speed of sound $c_s$, which describes the propagation for the wave components, and is given by the relation
\begin{equation}
    c_s^2=p'(\rho).
\end{equation}
In our context, we allow for the formation of vacuum states, which corresponds to the vanishing of the density $\rho$, and our domain $\Omega_t$ will therefore be given by $\displaystyle\Omega_t=\{x|\rho(t,x)>0\}$ with a moving boundary $\Gamma_t$, described as $\displaystyle\Gamma_t=\{x|\rho(t,x)=0\}$.

Physical considerations suggest that there exists only one stable physical regime, known as physical vacuum, and it is characterized by $\displaystyle d(x,\Gamma_t)\simeq c_s^2(x,t)$. This regime allows the free boundary to move with a bounded acceleration and velocity while also interacting with the interior, which allows internal waves traveling with velocity $c_s$ to reach the boundary in finite time. On the other hand, if the sound speed has a faster decay rate, then the particles on the boundary are expected to move linearly and independently, leading to a regime that can only last for a short time. When the decay rate is slower, this will instead correspond to an infinite initial acceleration of the boundary, which is another regime that is expected to break down in a short time.

Throughout history, two main approaches have been adopted in fluid dynamics: an Eulerian one, in which one fixes a reference frame and allows for the particles to move, and a Lagrangian one, in which the particles are instead fixed, and the reference frame is allowed to move. Both of these frameworks have been extensively studied in the context of the compressible Euler equations in the full space $\mathbb{R}^d$, where the matter of well-posedness is now well understood.

To be more specific, when one considers the compressible Euler equations in the full space $\mathbb{R}^d$, they can be regarded as a symmetric hyperbolic system. In particular,  one can then apply the results of Hughes-Kato-Marsden from \cite{Hughes-Kato-Marsden}  (see also Majda's work in \cite{Majda}), and deduce that the system is locally well-posed in $H^s$ in the Hadamard sense, where $s>\frac{d}{2}+1$, with the continuation  criterion
\begin{align*}
    \int_0^T\|\nabla(\rho,v)(t)\|_{L^\infty}\,dt<\infty
\end{align*}
If we restrict to the irrotational case, we can fully describe the problem purely in terms of the wave components $(\rho,\nabla\cdot v)$ and apply the general Strichartz estimates of Smith-Tataru from \cite{Smith-Tataru} in order to improve the local well-posedness regularity threshold to $\displaystyle s>\frac{d+1}{2}$, in dimensions $\displaystyle d=3,4,5$. On the other hand, it is not yet clear what condition on the vorticity would guarantee a similar result in the fully rotational setting; see the works of Disconzi-Luo-Mazzone-Speck \cite{Disconzi-Luo-Mazzone-Speck}, Wang \cite{Wang}, Andersson-Zhang \cite{Andersson-Zhang},\cite{Andersson-Zhang2}, and Zhang \cite{Zhang}.

By contrast, the free physical vacuum boundary problem hadn't been extensively studied in an Eulerian setting until recently, and all of the known results had been proved in Lagrangian coordinates, in a high regularity setting, or in indirectly defined spaces.

A fully Eulerian low regularity theory for the local well-posedness question for the compressible Euler was first developed by Ifrim and Tataru in \cite{Ifrim-Tataru}. In the aforementioned paper, they:
\begin{itemize}
    \item showed that the solutions are unique, under minimal regularity assumptions: $(\rho,v)\in\text{Lip}$. Moreover, they proved that the distance between different solution can be propagated in time, which in particular shows that the solutions are stable;
    \item constructed a scale of weighted Sobolev function spaces tailored to this problem;
    \item rigorously derived sharp, scale-invariant energy estimates within the aforementioned weighted Sobolev spaces, and also provided a minimal continuation criterion for the regularity of solutions, in particular requiring that $v\in\text{Lip}$;
    \item constructed smooth (regular) solutions in this newly developed Eulerian framework;
    \item by developing a nonlinear Littlewood-Paley theory, they were also able to construct rough solutions as limits of smooth solutions, and proved their continuous dependence on the initial data. The low regularity threshold that they obtained is also one degree of regularity above scaling, which morally corresponds to the $\frac{d}{2}+1$ threshold obtained by Hugh-Kato-Marsden in the whole space.
\end{itemize}
Later on, together with Disconzi, they obtained counterparts for these results for the relativistic Euler equations in a vacuum setting in \cite{Disconzi-Ifrim-Tataru}.

As for the fluid incompressible case, the problem poses different difficulties. For instance, as one approaches the boundary, the density does not tend to zero. Moreover, one must necessarily assume the Taylor sign condition, as proved by Ebin in \cite{Ebin}, and also include this parameter in the estimates. In this free boundary setting, analogous results to the aforementioned ones of Ifrim-Tataru and Disconzi-Ifrim-Tataru were obtained by Ifrim, Pineau, Tataru, and Taylor in \cite{Ifrim-Pineau-Tataru-Taylor}. We must also mention that as far fluid equations in a free boundary are concerned, Shatah and Zheng were the first to pursue an Eulerian approach in \cite{Shatah-Zeng,Shatah-Zeng2, Shatah-Zeng3}. They mainly focused on the free boundary problem with surface tension and showed the existence of a solution to the pure gravity problem in the zero surface tension limit, but their argument seems to require the curvature to be bounded, which would in turn require a higher degree of regularity than in \cite{Ifrim-Pineau-Tataru-Taylor}. The comparison between the latter  and the memoir of Wang, Zhang, Zhao, and Zeng \cite{Wang-Zhang-Zhao-Zeng} is easier to make, since one of the results in the latter establishes the existence and uniqueness of solutions at a level of regularity which is one derivative above scaling. However, they only considered graph domains with unbounded curvature. On the other hand, Ifrim, Pineau, Tataru, and Taylor significantly expanded the scope of the results to also allow for more complicated geometries, and were also the first to prove continuity of solutions with respect to the initial data, along with an enhanced uniqueness result, a construction of rough solutions as unique limits of smooth ones, refined low regularity energy estimates with pointwise geometric control parameters that only require very limited regularity, a new proof of the existence of smooth solutions, and an essentially scale invariant continuation criterion, akin to the one obtained by Beale-Kato-Majda for the incompressible Euler equations on the whole space. In particular, their approach differs significantly from the one in \cite{Wang-Zhang-Zhao-Zeng}, which in turn relies on the one adopted by Alazard-Burq-Zuily in \cite{Alazard-Burq-Zuily, Alazard-Burq-Zuily2}. See also the work of de Poyferr\'e \cite{dePoyferre}.

For the rest of this paper, we shall be concerned with the irrotational compressible Euler equations in a vacuum setting.

\subsection{Notations and the conserved energy}

The equations \eqref{rawEuler} can be rewritten in terms of the material derivative $\displaystyle D_t:=\partial_t+v\cdot\nabla$ as
\begin{equation}
\begin{cases}\begin{aligned}
  &D_t\rho+\rho\nabla\cdot v=0\\
  &\rho D_tv+\nabla p=0
  \end{aligned}
  \end{cases}
\end{equation}
Upon differentiating the equation describing the evolution of $\rho$ once again, we end up with the wave equation
\begin{align*}
    D_t^2\rho-\rho\nabla\cdot(\rho^{-1}p'(\rho)\nabla\rho)=\rho[(\nabla\cdot v)^2-\Tr((\nabla v)^2)]
\end{align*}
with propagation speed $c_s$. A similar computation shows that $\nabla\cdot v$ solves a similar equation.

It turns out that the vorticity $\omega=\curl v$ evolves according to the transport equation
\begin{align*}
    D_t\omega=-\omega(\nabla v)-(\nabla v)^T\omega.
\end{align*}
Therefore, the compressible Euler equations can be viewed as a coupled system consisting of a wave equation describing the evolutions of the variables $\rho$ and $\nabla\cdot v$, and a transport equation solved by the vorticity $\omega=\curl v$.

We can also see that the full Euler system admits the conserved energy given by
\begin{align*}
    E&=\int_{\Omega_t}\frac{1}{2}\rho |v|^2+\rho h(\rho)\,dx,
\end{align*}
where $h$ is the specific enthalpy, and is defined by
\begin{align*}
    h(\rho)=\int_0^\rho\frac{p(\lambda)}{\lambda^2}\,d\lambda.
\end{align*}
It also turns out that this energy can be regarded as a Hamiltonian in a certain formulation, see Chorin-Marsden \cite{Chorin-Marsden}, Ebin-Marsden \cite{Ebin-Marsden}, and Marsden-Ratiu-Weinstein \cite{Marsden-Ratiu-Weinstein}.

\subsection{The good variables}
In \cite{Ifrim-Tataru}, Ifrim and Tataru showed that the equation can be rewritten in a simpler way once we choose a pair of convenient variables. We shall briefly motivate their choice. 

When $\kappa=1$, we can immediately check that the original variables $(\rho,v)$ are convenient to use. However, this is no longer the case when $\kappa\neq 1$, but it turns out that we can also find a proper choice of variables in that setting. We first note that the constitutive law implies that the speed of sound satisfies the relation $c_s^2=(\kappa+1)\rho^\kappa$. Since we would like our solutions to exhibit decay and have linear behavior, a quick computation shows that for  $r'=\rho^{-1}p'(\rho)$.\, the equations are significantly simpler. 

We can see that this relation is guaranteed to hold once we set $\displaystyle r=\frac{\kappa+1}{\kappa}\rho^\kappa$, and this will enable us to rewrite the equations as
\begin{equation}\label{Euler}
\begin{cases}\begin{aligned}
  &D_tr+\kappa r(\nabla\cdot v)=0\\
  &D_tv+\nabla r=0
  \end{aligned}
  \end{cases}
\end{equation}
Throughout this paper we shall work with these equations, coupled with the condition $\curl v=0$. Our system then takes the form
\begin{equation}\label{irrotational Euler}
\begin{cases}\begin{aligned}
  &D_tr+\kappa r(\nabla\cdot v)=0\\
  &D_tv+\nabla r=0\\
  &\curl v=0.
  \end{aligned}
  \end{cases}
\end{equation}
Their linearized counterparts take the form
\begin{equation}\label{linearized Euler}
\begin{cases}\begin{aligned}
  &D_ts+w\cdot\nabla r+\kappa s(\nabla\cdot v)+\kappa r(\nabla\cdot w)=0\\
  &D_tw+(w\cdot\nabla)v+\nabla s=0,
  \end{aligned}
  \end{cases}
\end{equation}
respectively 
\begin{equation}\label{linearized irrotational Euler}
\begin{cases}\begin{aligned}
  &D_ts+w\cdot\nabla r+\kappa s(\nabla\cdot v)+\kappa r(\nabla\cdot w)=0\\
  &D_tw+(w\cdot\nabla)v+\nabla s=0\\
  &\curl w=0.
  \end{aligned}
  \end{cases}
\end{equation}

As we have seen in \cite{Avadanei}, these equations admit potential function formulations: there exist potential functions $\phi$ and $\psi$ such that 
  \begin{equation}\label{Stream function formulation}
  \begin{aligned}
  v&=\nabla\phi\\
      D_t\phi&=\frac{|v|^2}{2}-r,
      \end{aligned}
      \end{equation}
      where the velocity potential $\phi$ also solves
      \begin{equation}\label{Stream function formulation wave equation}
  \begin{aligned}
      D_t^2\phi&+\nabla r\cdot\nabla\phi-\kappa r\Delta\phi=0.
      \end{aligned}
  \end{equation}
 
  Since the linearized velocity is also irrotational, it also admits a potential $\psi$, which satisfies the following equations:
  \begin{equation}\label{Linearized stream function formulation}
  \begin{aligned}
  w&=\nabla\psi\\
      D_t\psi&=-s.
      \end{aligned}
      \end{equation}
     We also showed that the linearized velocity potential $\psi$ evolves according to the following acoustic wave equation:
      \begin{equation}\label{Linearized stream function formulation wave equation}
  \begin{aligned}
      D_t^2\psi&-\nabla r\cdot\nabla\psi-\kappa r\Delta\psi=\kappa s(\nabla\cdot v).
      \end{aligned}
  \end{equation}
This acoustic wave equation appears in helioseismology and models sollar oscillations. It arises as a simplified version of a scalar equation which was first introduced  by Gizon-Barucq-Duruflé-Hanson-Leguèbe-Birch-Chabassier-Fournier-Hohage-Papini in \cite{GBDHLBCFHP}; see also \cite{AHN1,AHN2,MHFG }. For uniqueness of the associated inverse problems under various hypotheses, see Muller's PhD thesis. Our equation is a simplified version of the one of stellar oscillations (see \cite{LO} and \cite{C}) and of Galbrun's equation, which was introduced in \cite{Galbrun}, but in spite of this, it still describes solar dynamics quite accurately; see also \cite{BMMPP,HB,BJM}, as well as \cite{HH,HLS,H}. In the field of helioseismology, two main approaches have been adopted: a global one, which relies on "hearing" the shape of the Sun by solving inverse spectral problems  (one of the results achieved through this methodology involves the radially symmetric model, see \cite{CGT,Cetal}), and a local one, in which one constructs two or three dimensional images of the interior of the Sun (see \cite{DJHP,Getal}) by solving inverse boundary problems and by taking into account data on travel times and correlations of solar oscillations. The field of helioseismology is crucial in understanding the solar cycle and space weather prediction. 

\subsection{Function spaces and energies}
With respect to the new variables $(r,v)$, the conserved energy can be written as
\begin{align*}
    E&=\int_{\Omega_t}r^{\frac{1-\kappa}{\kappa}}\left(r^2+\frac{\kappa+1}{2}rv^2\right)\,dx
\end{align*}
Based on this conserved quantity, Ifrim and Tataru introduced the base energy space $\mathcal{H}$ with norm
\begin{align*}
    \|(s,w)\|^2_{\mathcal{H}}=\int_{\Omega_t}r^{\frac{1-\kappa}{\kappa}}\left(s^2+\kappa rw^2\right)\,dx,
\end{align*}
in the context of the linearized equations \eqref{linearized Euler}
for functions $(s,w)$, which are defined almost everywhere in the gas domain $\Omega_t$.

In order to define higher regularity spaces, they considered the leading order operator $D_t^2-\kappa r\Delta$ from the associated wave equation, which is naturally associated to the acoustic metric
\begin{align*}
    g=r^{-1}dx^2.
\end{align*}
This was their motivation in defining the higher order Sobolev spaces $\mathcal{H}^{2k}$ for elements $(s,w)\in\mathcal{D}'(\Omega_t)\times\mathcal{D}'(\Omega_t)$, with norm
\begin{align*}
    \|(s,w)\|^2_{\mathcal{H}^{2k}}&=\sum_{\substack{|\beta|\leq 2k}}^{|\beta|-\alpha\leq k}\|r^\alpha\partial^\beta(s,w)\|^2_{\mathcal{H}},
\end{align*}
where $0\leq\alpha\leq k$. For fractional $k$, one can define the spaces $\mathcal{H}^{2k}$ by complex interpolation. See \cite{Ifrim-Tataru} for the precise definition and a detailed discussion.

 We also note that the compressible Euler system admits the following scaling symmetry:
    \begin{align*}
        r_\lambda(t,x)&=\lambda^{-2}r(\lambda t,\lambda^2x)\\
        v_\lambda(t,x)&=\lambda^{-1}v(\lambda t,\lambda^2x),
    \end{align*}
    in the sense that if $(r,v)$ is a solution, then so is $(r_\lambda,v_\lambda)$.

    We can also see that the critical exponent $k$, which is defined as the one for which the homogeneous counterpart of $\mathcal{H}^{2k}$ is invariant to scaling is given by
    \begin{align*}
        2k_0=d+1+\frac{1}{\kappa}
    \end{align*}
 In \cite{Ifrim-Tataru}, in order to obtain their energy estimates, Ifrim and Tataru also introduced the control parameters 
    \begin{align*}
        A=\|\nabla r-N\|_{L^\infty}+\|v\|_{\dot{C}^{\frac{1}{2}}},
    \end{align*}
    which is associated with the critical Sobolev exponent $2k_0$, and
    \begin{align*}
        B=\|\nabla r\|_{\tilde{C}^{0,\frac{1}{2}}}+\|\nabla v\|_{L^\infty},
    \end{align*}
    where
    \begin{align*}
    \|f\|_{\tilde{C}^{0,\frac{1}{2}}}=\sup_{\substack{x,y\in\Omega_t}}\frac{|f(x)-f(y)|}{r(x)^{\frac{1}{2}}+r(y)^{\frac{1}{2}}+|x-y|^{\frac{1}{2}}}
    \end{align*}
    They also defined the phase space
    \begin{align*}
        \mathbf{H}^{2k}=\{(r,v)|(r,v)\in\mathcal{H}^{2k}\}.
    \end{align*}
    Here, the defining function $r$ of the domain $\Omega_t$ directly determines the $\mathcal{H}^{2k}$ norm, which means that the phase space $H^{2k}$  is morally an infinite dimensional manifold. Their main local well-posedness result is as follows:
    
\begin{thm}\label{Initial Theorem}[Ifrim, Tataru,\cite{Ifrim-Tataru}]
    The system \eqref{Euler} is locally well-posed in the space $\mathbf{H}^{2k}$ for $k\in\mathbb{R}$ satisfying
    \begin{align*}
    2k>2k_0+1.
    \end{align*}
\end{thm}
One of the main ideas in their proof was to obtain energy estimates of the form
\begin{align*}
\|(r,v)(t)\|_{\mathcal{H}^{2s}}&\lesssim e^{\int_0^tC(A)B(u)\,du}\|(r,v)(0)\|_{\mathcal{H}^{2s}}.
\end{align*}

We also note that for the previously defined control parameters $A$ and $B$ we have the following Sobolev estimates
\begin{align*}
    A&\lesssim \|(r,v)\|_{\mathbf{H}^{2k}},k>k_0\\
    B&\lesssim \|(r,v)\|_{\mathbf{H}^{2k}},k>k_0+\frac{1}{2}.
\end{align*}
In light of these estimates, the $2k>2k_0+1$ threshold in Theorem \ref{Initial Theorem} arises naturally when one tries to control $B(t)$ pointwise in time. On the other hand, this control parameter appears in the energy estimates as being integrated in time, hence a natural question is whether this could be used in order to lower the regularity threshold. This in turn has served our motivation to look into whether our problem admits Strichartz estimates.

\subsection{Main results}

\begin{definition}
    Let $d,q,r\geq 2$, $\gamma\in\mathbb{R}$.
    \begin{enumerate}
        \item We say that $\displaystyle (q,r)$ is \textit{wave-admissible} if
    \begin{align*}
        \frac{1}{q}+\frac{d-1}{2r}\leq \frac{d-1}{4}.
    \end{align*}
    If we have equality in the previous inequality, we say that $(q,r)$ is \textit{sharp wave-admissible}.
    \item We say that $(q,r,\gamma)$ is a \text{wave Strichartz triple} if $(q,r)$ is wave-admissible and
    \begin{align*}
        \frac{1}{q}+\frac{d}{r}=\frac{d}{2}-\gamma.
    \end{align*}
    \end{enumerate}
\end{definition}
For reference, we now recall the Strichartz estimates for the wave equations in $\mathbb{R}^d$:
\begin{theorem}
    Let $d\geq 2$, $(q,r,\gamma)$ be wave-admissible, with $(d,q,r)\neq (3,2,\infty)$, $(q,r)\neq(\infty,\infty)$, and let $u$ be a solution to the initial value problem
    \begin{align*}
        (\partial_t^2-\Delta)u&=0\\
        (u,u_t)(0)&=(u_0,u_1).
    \end{align*}
    Then,
    \begin{align*}
        \|u\|_{L_t^qL_x^r(\mathbb{R}\times\mathbb{R}^d)}&\lesssim \|u_0\|_{\dot{H}_x^{\gamma}(\mathbb{R}^d)}+\|u_1\|_{\dot{H}_x^{\gamma-1}(\mathbb{R}^d)}.
    \end{align*}
\end{theorem}

It turns out that these solutions do indeed provide counterexamples to Strichartz estimates in the context of the wave equation. This is what we prove in our first main result:

\begin{theorem}\label{Counterexample Strichartz}
Let $(r,v)\in C([0,1];\mathbf{H}^{2k})$ with $k>k_0$ be a fixed irrotational solution to \eqref{rawEuler} on the time interval $[0,1]$ with $\displaystyle \int_0^1B(t)\,dt<\infty$, and let $(q,r,\gamma)$ be wave-admissible in dimension $d\geq 2$. Then, for every $\displaystyle s>k_0$ with $\displaystyle s<\frac{1}{2q}+\gamma+\frac{1}{2\kappa}+1$, there exists a sequence of solutions $\displaystyle(\psi^j)_{j\geq 0}$ to the initial value problem
\begin{align*}
    (D_t^2-\kappa r\Delta-\nabla r\cdot\nabla+\kappa(\nabla\cdot v)D_t)\psi^j&=0\\
    (\psi^j,D_t\psi^j)(0)=(\psi^j_0,\psi^j_1),
\end{align*}
 such that 
\begin{align*}
    \sup_{\substack{j\geq 0}}(\|(\psi^j_1,\nabla_x\psi^j_0)\|_{\mathcal{H}^{2s}})\leq 1,
\end{align*}
    and whenever $\displaystyle \alpha\in \left(0,\frac{1}{2q}+\gamma+\frac{1}{2\kappa}+1-s\right)$,
we have
\begin{align*}
    \lim_{\substack{j\rightarrow\infty}}2^{-2j\alpha}\|\nabla^2\psi^j(t,\cdot)\|_{L_t^qL_x^r([0,1]\times\Omega)}=\infty.
\end{align*}

\end{theorem}
We have seen that our problem admits a scaling symmetry. At the level of the linearized potential $\psi$, this takes the form

\begin{align*}
    \psi_\lambda(t,x)=\lambda^{-3}\psi(\lambda t,\lambda^2x)
\end{align*}
We shall next investigate Strichartz estimates for our problem in this intrinsic scaling. We can see that these estimates would need to have the form
\begin{align*}
    \||\nabla_x|^{\gamma}\psi\|_{L_t^qL_x^r([0,1]\times\Omega)}\lesssim \|(D_t\psi(0),\nabla_x\psi(0))\|_{\mathcal{H}}.
\end{align*}
Once we impose the scaling invariance, we naturally end up with the following
\begin{definition}
    Let $d,q,r\geq 2$, $\kappa>0$, and $\gamma\in\mathbb{R}$. We say that $\displaystyle (q,r,\gamma)$ is an \textit{Euler Strcihartz} triple if $(q,r)$ is wave-admissible and
    \begin{align*}
        \frac{1}{2q}+\frac{d}{r}=\frac{d}{2}+\frac{1}{2\kappa}+\gamma-1.
    \end{align*}
\end{definition}
Our second main result is the following
\begin{theorem}\label{Counterexample Strichartz 2}
Let $(r,v)\in C([0,1];\mathbf{H}^{2k})$ with $k>k_0$ be a fixed irrotational solution to \eqref{rawEuler} on the time interval $[0,1]$ with $\displaystyle \int_0^1B(t)\,dt<\infty$, and let $(q,r,\gamma)$ be an Euler Strichartz triple in dimension $d\geq 2$. Then, for every $\displaystyle s\geq 0$ with $\displaystyle 2k_0<2s<4-2\gamma$, there exists a sequence of solutions $\displaystyle(\psi^j)_{j\geq 0}$ to the initial value problem
\begin{align*}
    (D_t^2-\kappa r\Delta-\nabla r\cdot\nabla+\kappa(\nabla\cdot v)D_t)\psi^j&=0\\
    (\psi^j,D_t\psi^j)(0)=(\psi^j_0,\psi^j_1),
\end{align*}
 such that 
\begin{align*}
    \sup_{\substack{j\geq 0}}(\|(\psi^j_1,\nabla_x\psi^j_0)\|_{\mathcal{H}^{2s}})\leq 1,
\end{align*}
    and whenever $\displaystyle \alpha\in \left(0,2-\gamma-s\right)$,
we have
\begin{align*}
    \lim_{\substack{j\rightarrow\infty}}2^{-2j\alpha}\|\nabla^2\psi^j(t,\cdot)\|_{L_t^qL_x^r([0,1]\times\Omega)}=\infty.
\end{align*}
\end{theorem}
\begin{remark}
    In particular, Theorems \ref{Counterexample Strichartz} and \ref{Counterexample Strichartz 2} guarantee that that whenever $\displaystyle (q,r)=(2,\infty)$, $d>3$, and $\displaystyle 2k_0<2s<2k_0+\frac{1}{2}$, we can construct a sequence of solutions $\displaystyle(\psi^j)_{j\geq 0}$ to the initial value problem
\begin{align*}
    (D_t^2-\kappa r\Delta-\nabla r\cdot\nabla+\kappa(\nabla\cdot v)D_t)\psi^j&=0\\
    (\psi^j,D_t\psi^j)(0)=(\psi^j_0,\psi^j_1),
\end{align*}
 such that 
\begin{align*}
    \sup_{\substack{j\geq 0}}(\|(\psi^j_1,\nabla_x\psi^j_0)\|_{\mathcal{H}^{2s}})\leq 1,
\end{align*}
    and whenever $\displaystyle \alpha\in \left(0,k_0+\frac{1}{4}-s\right)$,
we have
\begin{align*}
    \lim_{\substack{j\rightarrow\infty}}2^{-2j\alpha}\|\nabla^2\psi^j(t,\cdot)\|_{L_t^2L_x^{\infty}([0,1]\times\Omega)}=\infty.
\end{align*}
In particular, this shows that the obstruction caused by gallery waves lies at least at $\frac{1}{2}$ degrees of regularity above scaling, and as we shall see in Theorem \ref{Strichartz Estimates}, this is actually the true threshold for at least one specific isotropic spherical Euler background $(r,v)$.
\end{remark}
The key idea behind the proofs of Theorems \ref{Counterexample Strichartz} and \ref{Counterexample Strichartz 2} is to transplant a construction that we carried out in \cite{Avadanei} for the model background $r=x_d$ and $v=0$.  In that previous situation, it turned out that for every given frequency, there exist geodesics with multiple periodic reflections next to the boundary, and exact solutions to the acoustic wave equations that are highly localized in tangential frequency and propagate along them, at least up until a certain time of coherence. In analogy with Ivanovici \cite{Ivanovici}, we shall call them \textit{gallery waves} or \textit{gallery modes}. For a full heuristic discussion of gallery waves in the case of the model background $r=x_d$ and $v=0$, see \cite{Avadanei}. In order to successfully transplant our solution, we construct a careful regularization of the variable $r$, which we reconcile with the forward flow of $v$.

We shall actually see that for a specific irrotational isotropic Euler background solution, the solutions to our linearized acoustic wave equation also satisfy Strichartz estimates, up to the threshold described in our counterexamples in Theorems \ref{Counterexample Strichartz} and \ref{Counterexample Strichartz 2}. The solution in question was first considered by Sideris in \cite{Sideris2} as part of an affine family and is given by
\begin{align*}
    \Omega_t&=a(t)\mathbb{B}\\
    v(t,x)&=\frac{\dot{a}(t)}{a(t)}x\\
    r(t,x)&=ca(t)^{-d\kappa}\left(1-\frac{|x|^2}{a(t)^2}\right)\\
    \ddot{a}(t)&=2ca^{-d\kappa-1}.
\end{align*}
Here, $\mathbb{B}$ is the unit ball, $\displaystyle c>0$, $\displaystyle a(0)=1$, and $\displaystyle a(t)>0$ for every $\displaystyle t\in[0,T]$. These are irrotational solutions, and in three dimensions, they were studied by Sideris in \cite{Sideris} as part of an affine family of solutions.

In particular, this shows that our thresholds are optimal for this specific background. We have the following
\begin{theorem}\label{Strichartz Estimates}
   Let $d\geq 2$, $\kappa>0$, and let $\displaystyle (r,v)$ be the aforementioned spherical background solution. Let $T_0>0$ be such that $a(t)>0$ on $[0,T_0]$, and $(q,r,\gamma)$ be a wave-admissible Strichartz triple so with $(d,q,r)\neq(3,2,\infty)$. Let $\displaystyle s>\frac{1}{2q}+\gamma+\frac{1}{2\kappa}+1$. Then, the solution $u$ to the initial value problem
   \begin{align*}
       (D^2-\kappa r\Delta-\nabla r\cdot \nabla+\kappa(\nabla\cdot v)D_t)u&=0\\
       (u,D_tu)(0)&=(u_0,u_1),
   \end{align*}
   satisfies
   \begin{align*}
       \|\nabla_x^2u\|_{L_t^qL_x^r([0,T_0]\times\Omega)}&\lesssim_{\kappa,q,r} \|(u_1,\nabla_xu_0)\|_{\mathcal{H}^{2s}(\Omega)}.
   \end{align*}
   When $(d,q,r)=(2,4,\infty)$, we instead have the estimate
   \begin{align*}
       \|\nabla_x^2u\|_{L_t^qL_x^r([0,T_0]\times\Omega)}&\lesssim_{T_0,\kappa,q,r} \|(u_1,\nabla_xu_0)\|_{\mathcal{H}^{2s}(\Omega)}.
   \end{align*}
\end{theorem}
\begin{remark}
    In particular, Theorem \ref{Strichartz Estimates} implies that when $d\geq 4$ and $\displaystyle (q,r)=(2,\infty)$, for every $\displaystyle s\in\mathbb{R}$ with $\displaystyle 2s>2k_0+\frac{1}{2}$, the solution to the initial value problem
\begin{align*}
       (D_t^2-\kappa r\Delta-\nabla r\cdot\nabla+\kappa(\nabla\cdot v)D_t)u&=0\\
       (u,D_tu)(0)&=(u_0,u_1),
   \end{align*}
   satisfies
   \begin{align*}
       \|\nabla_x^2u\|_{L_t^2L_x^\infty([0,T_0]\times\Omega)}&\lesssim_{\kappa,q,r} \|(u_1,\nabla_xu_0)\|_{\mathcal{H}^{2s}(\Omega)},
   \end{align*}
   whenever $(2,\infty)$ is admissible (we note that this requires $d\geq 4$, which excludes the previous forbidden endpoint). Together with the counterexample theorem, this implies that the threshold $\displaystyle 2k_0+\frac{1}{2}$ is optimal. However, we do not discuss the behavior at the aforementioned threshold. 
\end{remark}
\subsection{Further historical comments}
The compressible Euler equations have been extensively studied over time, and they have also received a lot of interest from physicists. When one allows for vacuum states, which correspond to the vanishing of the density, the speed of sound vanishes at the boundary, which forces the acoustic geometry to degenerate, and significantly complicates the behavior of the gas. However, this is also the natural boundary conditions for compressible gases. In this context, the space is regarded as consisting of a particle region $\Omega_t$ and a vacuum one, both of which are separated by a moving free boundary $\Gamma_t=\partial\Omega_t$. There two main scenarios which can occur here, depending on the behavior of the sound speed $c_s$ at the boundary, which is in turn determined by the density:
\begin{itemize}
    \item A fluid case, in which we assume that the density and the sound speed have a nonzero positive limit at the free boundary.
    \item A gas scenario, in which on the other hand the density is assumed to tend to zero near the free boundary. Throughout this paper, we shall focus on this case.
\end{itemize}
General fluid flows were analyzed by Christodoulou-Miao \cite{Christodoulou-Miao} and Lindblad \cite{Lindblad}, whereas incompressible ones were investigated by Lindblad-Luo \cite{Lindblad-Luo}. In the incompressible case, we once again mention the more recent results of Ifrim-Pineau-Tataru-Taylor \cite{Ifrim-Pineau-Tataru-Taylor} that we have already discussed in this introduction.

Depending on the behavior of the speed of sound relatively to the distance to the boundary, we distinguish three cases:
\begin{enumerate}
    \item[a)]  $c_s\lesssim d_{\Gamma_t}$, which is known as the fast decay scenario. Physical considerations suggest that in this case, the vacuum boundary will evolve linearly, which will in turn prevent internal waves form  reaching it arbitrarily fast. In particular, this scenario will survive at least on short times, and one can apply the known well-posedness results on symmetric hyperbolic systems in this context, see for instance DiPerna \cite{DiPerna}, Chen \cite{Chen}, and Lions \cite{Lions}, as well as Kawashima-Makino-Ukai \cite{Kawashima-Makino-Ukai}, Liu-Yang \cite{Liu-Yang}, and Chemin \cite{Chemin}. In particular, it follows that this case does not correspond to a genuine free boundary problem, and furthermore, the geometry breaks down in finite time, as Chemin \cite{Chemin} proved in the one dimensional case.
    \item[b)] $c_s\gg d_{\Gamma_t}$, which is known as the slow decay scenario. Under this condition, it turns out that waves can reach the boundary arbitrarily fast, a phenomenon which causes the internal flow and the one of the free boundary to be strongly coupled, which in particular means that this scenario yields a genuine free boundary problem. We distinguish a natural family of decay rates given by $c_s\approx d_{\Gamma_t}^\beta$, with $\beta\in(0,1)$. Physical and mathematical considerations suggest that among all of these decay rates, there exists a single stable one, which occurs when $\beta=\frac{1}{2}$. The other situations have been conjectured to be unstable and to instantaneously switch into the stable regime. However, we must mention that to the best of our knowledge, no rigorous results have been proved in this direction, and it is seems very likely that there are significant differences between the regimes $\beta<\frac{1}{2}$ and $\beta>\frac{1}{2}$.
    \end{enumerate}

 The first case that has been considered in the physical vacuum scenario was in one dimension, and it was investigated in Coutand-Shkoller \cite{Coutand-Shkoller} and Jang-Masmoudi \cite{Jang-Masmoudi}. While Coutand and Shkoller prove some energy estimates and provide a procedure to construct solutions, the function spaces are not completely defined, and the initial data is not directly described. This matter is addressed by Jang and Masmoudi, who introduced the Lagrangian counterparts of Ifrim and Tataru's weighted Sobolev spaces, and prove the existence and uniqueness of solutions in sufficiently regular spaces. 

 Another case that received significant interest was the three dimensional one. In this scenario, Coutand-Lindblad-Shkoller \cite{Coutand-Lindblad-Shkoller} formally derived energy estimates in the case $\kappa=1$. They later on also proved existence by using a parabolic regularization procedure in \cite{Coutand-Shkoller2}, but also in a functional setting that is not completely described; moreover, the difference bound that they proved required more regularity than their existence result.  Jang and Masmoudi \cite{Jang-Masmoudi2} also independently used a parabolic regularization in order to show that solutions and exist and are unique for $\kappa>0$, with a completely different proof of the energy estimates. However, they only proved their results on the torus in a Lagrangian setting, and briefly outlined the general case.

 All of the aforementioned results for the vacuum case were proved in a Lagrangian setting,  and as we have already mentioned in this introduction, the first fully Eulerian results were obtained by Ifrim and Tataru  \cite{Ifrim-Tataru}. Together with Disconzi, they also showed that the counterparts of these results hold for the relativistic equations \cite{Disconzi-Ifrim-Tataru}. Later on, together with Pineau and Taylor they also proved the counterparts of these results in the incompressible fluid case \cite{Ifrim-Pineau-Tataru-Taylor}.

We must also note that the question of global existence and nonlinear stability of expanding affine solutions had also been studied before. The first results of these kind were proved by Hadzi\' c and Jang in \cite{Hadzic-Jang} in the case $0<\kappa\leq\frac{2}{3}$, while Shkoller and Sideris later on extended these results to the full range of adiabatic exponents in \cite{Shkoller-Sideris}.

 Strichartz estimates for the wave and Schr\" odinger equation have a very long history, and the first setting in which they were obtained was the case $q=r$ for the wave and classical Schr\" odinger equations by Strichartz in \cite{Strichartz}. Later on, these results were extended to mixed $L_t^qL_x^r$ norms by Ginibre and Velo \cite{Ginibre-Velo} in the case of the the Schr\" odinger equation, and then independently by Ginibre-Velo \cite{Ginibre-Velo2} and Lindblad-Sogge \cite{Lindblad-Sogge} for the wave equation, following earlier work by Kapitanskii \cite{Kapitanski}. Keel-Tao \cite{Keel-Tao} settled the remaining endpoints. 
 
 The first such estimates for smooth variable coefficients were proved independently by Mockenhaupt-Seeger-Sogge \cite{Mockenhaupt-Seeger-Sogge} and Kapitanskii \cite{Kapitanskii2}. In the case of rough coefficients, the first estimates of this kind were proved by Smith \cite{Smith} for metrics of class $C^2$, who used wave packet techniques. Later on, Smith and Sogge proved that the result may fail for metrics of H\"older regularity $C^\alpha$, with $\alpha<2$.
 
 Later on, Bahouri-Chemin \cite{Bahouri-Chemin3} and Tataru \cite{Tataru-I} independently obtained improved Stricartz estimates with a $\frac{1}{4}$-derivative loss.  Tataru then extended and improved these results in all dimensions in \cite {Tataru-II} and \cite{Tataru-III}, where he proved Strichartz estimates with a $\frac{1}{6}$-derivative loss.  Bahouri and Chemin also improved their own results in the meantime in \cite{Bahouri-Chemin2} up to a loss of derivatives slightly better than $\frac{1}{5}$. On the other hand, Smith and Tataru proved afterwards in \cite{Smith-Tataru-counter} that a loss of $\frac{1}{6}$ derivatives is optimal for general metrics of class $C^1$. Smith and Tataru \cite{Smith-Tataru} also obtained lossless Strichartz estimates for quasilinear wave equations in dimensions $d=2$ and $d=3$.

In the case of manifolds with smooth and strictly geodesically concave boundary, Smith and Sogge \cite{Smith-Sogge} proved Strichartz estimates in the non-endpoint cases for the corresponding wave equation by using the Melrose and Taylor parametrix. It turns out that the concavity condition is crucial, since one could otherwise have multiply reflecting geodesics, and gliding rays that arise as their limits, which would prevent the existence of the aforementioned parametrix.

"Log-loss" estimates for the spectral clusters on compact manifolds without boundary were proved by Koch, Smith, and Tataru \cite{Koch-Smith-Tataru}. For results in the context of the fractional Schr\" odinger and wave equations on compact manifolds  without boundary, see Dinh's work \cite{Dinh}. In the case of manifolds with boundaries,  Burq, Lebeau, and Planchon \cite{Burq-Lebeau-Planchon} used the $L^r(\Omega)$ estimates that had been obtained by Smith and Sogge \cite{Smith-Sogge2} for a class of spectral operators in order to prove Strichartz estimates in this setting. On the other hand, we must note that the range of indices $(q,r)$ is determined by the admissible range for $r$ in the square function estimate for the wave equation, which in turn controls $u$ in the space $\displaystyle L^r(\Omega,L^2((-T,T))$. For instance, this requires $\displaystyle q,r\geq 5$ in three dimensions. The result from \cite{Burq-Lebeau-Planchon} was later extended by Blair, Smith and Sogge \cite{Blair-Smith-Sogge}, who showed that if $\Omega$ is a compact manifold with boundary, and $(q,r,\beta)$ a triple with
\begin{align*}
    \frac{1}{q}+\frac{d}{r}=\frac{d}{2}-\beta,
\end{align*}
subject to the admissibility conditions
\begin{align*}
    \frac{3}{q}+\frac{d-1}{r}\leq\frac{d-1}{2},\text{ when }d\leq 4,\,\text{ and }\frac{1}{q}+\frac{1}{r}\leq\frac{1}{2},\text{ when }d\geq 4,\text{ and }q\in(2,\infty],\quad r\in[2,\infty),
\end{align*}
then one has the expected Strichartz estimates for the wave equations with homogeneous Dirichlet or Neumann boundary conditions. Here, the implicit constants only depend on $\Omega$ and $T$. 

For strictly convex manifolds with boundary, Ivanovici constructed counterexamples to Strichartz estimates in \cite{Ivanovici} in the case of the upper half-space, and she also investigated the extent to which they hold. Later on, Ivanovici, Lebeau, and Planchon studied dispersive estimates for the general wave equation in the Friedlander model case. In the case of the Klein-Gordon and wave equations, Ivanovici obtained large time dispersive estimates in \cite{Ivanovici2}. In the general case, we refer the reader to the works of Ivanovici, Lascar, Lebeau, and Planchon \cite{Ivanovici-Lebeau-Planchon2}, and of Ivanovici, Lebeau, and Planchon \cite{Ivanovici-Lebeau-Planchon3}. The same authors also constructed counterexamples to Strichartz estimates in two dimensional convex domains in \cite{Ivanovici-Lebeau-Planchon4}. In the case of the wave equations with Neumann boundary conditions inside general strictly convex domains, see the results of Ivanovici and Sanchez-Sanz \cite{Ivanovici-Sanz}.  Ivanovici also proved Strichartz and dispersive estimates for the quantum ball model in \cite{Ivanovici2025}. 

\subsection{An outline of the paper}
Here, we briefly describe the main structure of our paper, and some of the main ideas in each section. 
 \subsubsection{Preliminaries} In this section we recall the construction of Ifrim and Tataru's weighted Sobolev spaces $\mathcal{H}^{2s}$, along with their Sobolev embeddings and interpolation properties. In the end, we also recall the potential function formulations for the irrotational compressible Euler system and its linearized counterpart that we obtained in \cite{Avadanei}. 
 \subsubsection{A regularized domain defining function} We also introduce a regularization $\tilde{r}$ of the domain defining function $r$, which is adapted in frequency to every single boundary layer of the domain $\Omega_t$. This construction will in particular enable us to transplant our construction of the counterexample from \cite{Avadanei} for general backgrounds.
 
 \subsubsection{Gallery modes in modified coordinates} In this section we shall transplant our construction for the counterexamples from \cite{Avadanei}. We do so by essentially carrying out two steps at the same time: we use the regularized function $\tilde{r}$ as an approximate normal coordinate which captures the geometrical condition $A\ll 1$, as well as the flow of $v$ in order to account for the Galilean transport of our solution.
 \subsubsection{A physical spherical solution} In this short section we define an explicit affine spherical solution to compressible irrotational Euler equations which was first considered by Sideris in \cite{Sideris}. This is the solution that will provide the background for which we shall prove our desired Strichartz estimates in the desired open range and obtain results analogous to the ones in \cite{Avadanei}.
 \subsubsection{A spectral analysis of the acoustic operators} In this subsection we consider the acoustic wave operator associated to our spherical background, and study its spectral properties. In particular, we describe its eigenvalues and eigenfunctions in terms of radial functions and spherical harmonics. This step will be essential in the proof of our main results, which relies on a decomposition of the solution into a sum of dyadic cluster components.
 \subsubsection{The proof of the Strichartz estimates} In this section we provide the proofs of our main results. We decompose the solutions to our acoustic wave equations into dyadical spherical clusters, and prove a dispersive estimate by using the $TT^{\ast}$ argument. The key ingredients here consist of estimates for the corresponding kernel of the $TT^{\ast}$ operator.

 \subsection{Acknowledgements} The author would like to thank Daniel Tataru and Mihaela Ifrim for very many helpful discussions.

 The author also used ChatGPT 5.6 Sol and OpenAI Codex GPT-6 Astra for exploratory discussion, assistance in checking calculations and arguments, and editorial feedback during the preparation of this manuscript. The mathematical statements, proofs, and final formulations are humanly written and are the author’s responsibility.

 Part of this work was done while the author was participating in the "Mathematical Problems in Fluid Dynamics Reunion" at the Simons Laufer Mathematical Sciences Institute (July 2023) and the "Nonlinear Waves and Relativity" thematic program at the Erwin Schr\"odinger International Institute for Mathematics and Physics (May 2024).

 The author was also supported by the NSF
grant DMS-2054975, the Simons Foundation, as well as the NSF grant DMS-1928930 while in residence at the Simons Laufer Mathematical Sciences Institute (formerly MSRI) in Berkeley, California. The author would also like to thank the Yale University Mathematics Department for its hospitality during the $2025-2026$ academic year.
\section{Preliminaries}
In this section we are going to introduce the main function spaces in which we are going to work. We are going to work in a class of Sobolev spaces that were first introduced by Ifrim-Tataru. Throughout this section, we shall follow the exposition from Section 2 in \cite{Ifrim-Tataru}.
\subsection{Weighted Sobolev spaces}
\begin{definition}
    Let $\sigma>-\frac{1}{2}$ and $j\geq 0$. We define $\displaystyle H^{j,\sigma}(\Omega)$ as the space of all distributions $f$ in $\Omega$ with the property that the norm given by
    \begin{align*}
        \|f\|^2_{H^{j,\sigma}}&=\sum_{\substack{|\alpha|\leq j}}\|r^{\sigma}\partial^\alpha f\|^2_{L^2}
    \end{align*}
\end{definition}
Of course, we can also define the homogeneous counterparts for these spaces by replacing the norm with
\begin{align*}
     \|f\|^2_{\dot{H}^{j,\sigma}}&=\sum_{\substack{|\alpha|= j}}\|r^{\sigma}\partial^\alpha f\|^2_{L^2}
\end{align*}
By complex interpolation, we can extend this definition to the fractional case, defining $H^{s,\sigma}(\Omega)$ and $\dot{H}^{s,\sigma}(\Omega)$, for $s\geq 0$. 

The reason why we require $\sigma>-\frac{1}{2}$ is that we don't make any vanishing assumptions on the free boundary for the elements in the function spaces.

We now enumerate a few embedding properties that these spaces satisfy:
\begin{lemma}
    We assume that $s_1>s_2>0$, and that $\sigma_1>\sigma_2>-\frac{1}{2}$, with the property that $s_1-s_2=\sigma_1-\sigma_2$. Then,
    \begin{align*}
        H^{s_1,\sigma_1}\subset H^{s_2,\sigma_2}.
    \end{align*}
    In addition, if $\sigma>0$ and $\sigma\leq j$, then
    \begin{align*}
        H^{j,\sigma}\subset H^{j-\sigma}.
    \end{align*}
\end{lemma}
\begin{proof}
    See the proof of Lemmas 2.2 and 2.3 in \cite{Ifrim-Tataru}.
\end{proof}
We also record the following Morrey type embeddings, which follow from usual Sobolev embeddings:
\begin{lemma}
    For $\displaystyle\sigma\geq 0$ and $\displaystyle 0\leq s\leq j-\sigma-\frac{d}{2}$, we have the embedding
    \begin{align*}
        H_r^{j,\sigma}\subset C^s,
    \end{align*}
    where the equality $\displaystyle s=j-\sigma-\frac{d}{2}$ can occur only if $s$ is not an integer.
\end{lemma}
\subsection{Weighted Sobolev spaces adapted to the compressible Euler equation}
The first motivation for the choice of spaces is the following quantity, which is conserved as long as $(r,v)$ is a solution to the compressible Euler equation:
\begin{align*}
    E=\int_{\Omega_t}r^{\frac{1-\kappa}{\kappa}}\left(r^2+\frac{\kappa+1}{2}rv^2\right)\,dx
\end{align*}

This motivated the choice of the following energy functional, which was introduced and used by Ifrim-Tataru in order to study the linearized equation:
\begin{align*}
    E_{\text{lin}}(s,w)&=\int_{\Omega_t}r^{\frac{1-\kappa}{\kappa}}(s^2+\kappa rw^2)\,dx
\end{align*}

This motivated them to define the baseline space $\mathcal{H}$, with norm given by
\begin{align*}
    \|(s,w)\|^2_{\mathcal{H}}=E_{\text{lin}}(s,w).
\end{align*}
Equivalently,
\begin{align*}
    \mathcal{H}:=H^{0,\frac{1-\kappa}{2\kappa}}\times H^{0,\frac{1}{2\kappa}}
\end{align*}
In order to define higher order Sobolev spaces, one needs to determine the correct scale that gives the balance between time and spatial derivatives. This can be done by writing the wave component of our equation as a second order equation, with leading order operator $\displaystyle D_t^2-r\Delta$, which motivates the following definition:
\begin{align*}
    \mathcal{H}^{2k}:=H^{2k,k+\frac{1-\kappa}{2\kappa}}\times H^{2k,k+\frac{1}{2\kappa}}
\end{align*}
with $k\geq 0$. When $k$ is an integer, they can be defined by the explicit norm:
\begin{align*}
    \|(s,w)\|^2_{\mathcal{H}^{2k}}:=\sum_{\substack{|\beta|\leq 2k\\
    0\leq \alpha\leq k}}^{|\beta|-\alpha\leq k}\|r^\alpha\partial^\beta(s,w)\|^2_{\mathcal{H}},
\end{align*}
where $\alpha\geq 0$. We can analogously define their homogeneous counterparts $\dot{\mathcal{H}^{2k}}$.

These are all Hilbert spaces, and due to the fact that in this case all interpolation methods give the same result , one can define $\mathcal{H}^{2k}$ by using an interpolation method that is akin to the $J$-method. This is the object of the following
\begin{lemma}
    Let $0<k<N$. We define $\mathcal{H}^{2k}$ as the space of distributions $(s,w)$ that admit a representation
    \begin{align*}
(s,w)=\sum_{\substack{l=0}}^\infty(s_l,w_l),
    \end{align*}
    where $\{(s_l,w_l)\}$ has the property that
    \begin{align*}
        |||\{(s_l,w_l)\}|||^2:=\sum_{\substack{l=0}}^\infty 2^{2kl}\|(s_l,w_l)\|^2_{\mathcal{H}}+2^{2l(k-N)}\|(s_l,w_l)\|^2_{\mathcal{H}^{2N}}<\infty
        \end{align*}
        The norm of $(s,w)$ is defined by
        \begin{align*}
            \|(s,w)\|^2_{\mathcal{H}^{2k}}:=\inf |||\{(s_l,w_l)\}|||^2,
        \end{align*}
        where the infimum on the right hand side is taken over all representations of $(s,w)$.
\end{lemma}
\subsection{The state space $\mathbf{H}^{2k}$} $\mathbf{H}^{2k}$ is defined for $2k>2k_0$ (above scaling) as the set of pairs $(r,v)$ of functions on $\Omega_t$ that satisfy the following properties:
\begin{itemize}
    \item Boundary regularity: $\Gamma_t$ is a Lipschitz surface.
    \item Nondegeneracy: $r$ is a Lipschitz function on $\overline{\Omega}_t$, positive in the interior of $\Omega_t$, which vanishes simply on the boundary $\Gamma_t$
    \item Regularity: $(r,v)\in\mathcal{H}^{2k}$.
\end{itemize}
As $\Omega_t$ depends on $r$, $\mathbf{H}^{2k}$ cannot be regarded as a linear space, but as an infinite dimensional manifold. Due to the fact that our particle flows are going to be merely continuous on $\mathbf{H}^{2k}$, it suffices to only define a topology on these manifolds:
\begin{definition}\label{Topology}
    A sequence $((r_n,v_n))_{n\geq 0}$ converges to $(r,v)$ in $\mathbf{H}^{2k}$ if the following hold:
    \begin{enumerate}
        \item[a)] Uniform nondegeneracy: $|\nabla r_n|\geq c>0$ at the boundary for all $n\geq 0$ for some positive constant $c$;
        \item[b)] Domain convergence: $\displaystyle\lim_{\substack{n\rightarrow\infty}}\|r_n-r\|_{\text{Lip}}=0$;
        \item[c)] Norm convergence: For every $\varepsilon>0$, there exist smooth functions $(\tilde{r}_n,\tilde{v}_n)$ in $\Omega_n$, and $(\tilde{r},\tilde{v})$ in $\Omega$, such that
        \begin{align*}
          (\tilde{r}_n,\tilde{v}_n)\rightarrow (\tilde{r},\tilde{v})\text{ as }n\rightarrow\infty\text{ in }C^\infty  
        \end{align*}
        and
        \begin{align*}
            \|(\tilde{r}_n,\tilde{v}_n)-(r_n,v_n)\|_{\mathcal{H}^{2k}(\Omega)}&\leq\varepsilon\\
            \|(\tilde{r},\tilde{v})-(r,v)\|_{\mathcal{H}^{2k}(\Omega)}&\leq\varepsilon.
        \end{align*}
    \end{enumerate}
    
\end{definition}
This definition will allow us to prove two key properties of the solutions to our equation:
    \begin{enumerate}
        \item The solutions $(r,v)$ are continuous in time with values in $\mathbf{H}^{2k}$.
        \item The solutions $(r,v)$ in $\mathbf{H}^{2k}$ depend continuously on the initial data in $\mathbf{H}^{2k}$.
    \end{enumerate}

\subsection{Regularizing kernels}
A very important tool in our analysis consists of a family of regularization operators that are akin to the Littlewood-Paley projectors, adapted to our free boundary domains, at each dyadic scale $2^h$, where $h\geq 0$. They will need to have the following properties:
\begin{itemize}
    \item Fixed domain regularization. Given $(s,w)\in\mathcal{H}^{2k}(\Omega)$, the operators need to produce regularizations $\Psi^h(s,w)$ within the same $\mathcal{H}^{2j}(\Omega)$ scale of spaces. 
    \item State and domain regularization. Given $(r,v)\in\mathbf{H}^{2k}(\Omega)$, where $r$ defines a domain $\Omega$, the operators need to produce regularizations $\Psi^h(r,v)$ within the same $\mathbf{H}^{2j}(\Omega)$ scale of spaces, where the regularized domains are defined by the function $r^h$, as $\Omega_h:=\{x\in\mathbb{R}^d|r^h(x)>0\}$.
\end{itemize}

In order to better describe our analysis, to each fixed dyadic frequency $2^h$, we shall associate the following:
\begin{enumerate}
    \item[1.] Dyadic boundary layer: $\displaystyle\Omega^{[h]}=\{x\in\Omega,r(x)\approx 2^{-2h}\}$
    \item[2.] Full boundary strip: $\displaystyle\Omega^{[>h]}=\{x\in\Omega,r(x)\lesssim 2^{-2h}\}$
    \item[3.] Interior region: $\displaystyle\Omega^{[<h]}=\{x\in\Omega,r(x)\gtrsim 2^{-2h}\}$
\end{enumerate}
We also introduce dyadic enlargements of $\Omega$, defined by
\begin{align*}
    \tilde{\Omega}^{[h]}=\{x\in\mathbf{R}^d,d(x,\Omega)\leq c2^{-2h}\},
\end{align*}
and
\begin{align*}
    \tilde{\Omega}^{[>h]}=\{x\in\mathbf{R}^d,d(x,\Gamma)\leq c2^{-2h}\},
\end{align*}
where $c$ is a small positive universal constant.

Given a domain $\Omega$ defined by a nondegenerate Lipschitz function $r$, and the functions $(s,w)$ defined on $\Omega$, we define our regularizations by the relation 
\begin{align*}
(s^h,w^h)&=\Psi^h(s,w)=\int K^h(x,y)(s,w)(y)\,dy,
\end{align*}
where the kernels $K^h$ form a family with the following properties:
\begin{enumerate}
    \item[1.]Domain and localization: \begin{align*}
      K^h:\tilde{\Omega}^{[h]}\times\Omega\rightarrow\mathbf{R},  
    \end{align*} with support properties
    \begin{align*}
      \supp K^h\subset\{(x,y)\in\tilde{\Omega}^{[h]}\times\Omega^{<h},|x-y|\lesssim\delta y^h:=2^{-2h}+2^{-h}r(y)^{\frac{1}{2}}\}  
    \end{align*}
    \item[2.]Size and regularity:
    \begin{align*}
        |\partial_x^\alpha\partial_y^\beta K^h(x,y)|\lesssim (2^{-2h}+2^{-h}r(y)^{\frac{1}{2}})^{-d-|\alpha|-|\beta|},|\alpha|+|\beta|\leq 4N,
    \end{align*}
    where $N$ is large enough.
    \item[3.] Approximate identity:
    \begin{align*}
        \int K^h(x,y)\,dy&=1\\
        \int K^h(x,y)(x-y)^\alpha\,dy&=0,1\leq|\alpha|\leq 2N.
    \end{align*}
\end{enumerate}
The first property ensures that the regularizations $(s^h,w^h)$ in the extended domain $\tilde{\Omega}^{[h]}$ have the following mapping dyadic properties: 
\begin{itemize}
    \item For $j<h$, the regularizations $(s^h,w^h)$ in $\Omega^{[j]}$ are determined by $(s,w)$ in $\Omega^{[j]}$.
    \item For the $h$ layers, the regularizations $(s^h,w^h)$ in $\tilde{\Omega}^{[>h]}$ are determined by $(s,w)$ only in $\Omega^{[h
    ]}$.
\end{itemize}
The existence of these regularization kernels is guaranteed by Lemma 2.9 in \cite{Ifrim-Tataru}.

Given a family of solutions $(r^h,v^h)$ to \eqref{Euler}, for every $j$, we also consider the the family of sets $\displaystyle \Omega^{h,j}=\{x\in\Omega^h|r^h(x)\approx 2^{-2j}\}$, where $\Omega^h=\{x\in\mathbb{R}^d|r^h(x)>0\}$, along with the clear analogues $\Omega^{h,>j}$ and $\Omega^{h,<j}$, $\tilde{\Omega}^{h,j}$, and $\tilde{\Omega}^{h,>j}$.

The regularization operators $\Psi^h$ satisfy some properties that are akin to Bernstein's inequalities for Littlewood-Paley projectors. They are outlined in the following
\begin{proposition}
    Let $r_1$ be a nondegenerate defining function such that $|r-r_1|\ll 2^{-2h}$. Then, the following hold:
    \begin{enumerate}
    \item[a)]Regularization bound:
        \begin{align*}
            \|\Psi^h(s,w)\|_{\mathcal{H}_{r_1}^{2k+2j}}\lesssim 2^{2jh}\|(s,w)\|_{\mathcal{H}_r^{2k}},j\geq 0
        \end{align*}
        \item[b)] Difference bound:
        \begin{align*}
            \|(\Psi^{h+1}-\Psi^{h})(s,w)\|_{\mathcal{H}_{r_1}^{2k+2j}}\lesssim 2^{2jh}\|(s,w)\|_{\mathcal{H}_r^{2k}},-k\leq j\leq 0,
        \end{align*}
        \item[c)] Error bound:
        \begin{align*}
            \|(I-\Psi^h)(s,w)\|_{\mathcal{H}_r^{2k+2j}}\lesssim 2^{2jh}\|(s,w)\|_{\mathcal{H}_r^{2k}},-k\leq j\leq 0.
        \end{align*}
    \end{enumerate}
\end{proposition}
\begin{proof}
    These follow by taking $(s^h,w^h):=\Psi^h(s,w)$. See \cite{Ifrim-Tataru} for further details.
\end{proof}
\subsection{Interpolation inequalities} The main result is the following
\begin{proposition}
    Let $\sigma_0,\sigma_m\in\mathbb{R}$, and $1\leq p_0,p_m\leq\infty$. Let $\displaystyle \theta_j=\frac{j}{m}$, $\displaystyle \frac{1}{p_j}=\frac{1-\theta_j}{p_0}+\frac{\theta_j}{p_m}$, $\displaystyle \sigma_j=\sigma_0(1-\theta_j)+\sigma_m\theta_j$.

    We also assume that
    \begin{align*}
m-\sigma_m-d\left(\frac{1}{p_m}-\frac{1}{p_0}\right)>-\sigma_0,\sigma_j>-\frac{1}{p_j}.
    \end{align*}

    Then, for $0<j<m$, we have
    \begin{align*}
        \|r^{\sigma_j}\partial^jf\|_{L^{p_j}}&\lesssim \|r^{\sigma_0}f\|^{1-\theta_j}_{L^{p_0}}\|r^{\sigma_m}\partial^mf\|^{\theta_j}_{L^{p_m}}
    \end{align*}
\end{proposition}
\begin{proof}
    See the proof of Proposition 2.12 in \cite{Ifrim-Tataru}.
\end{proof}
We make the convention that the following three inequalities will only be used in local coordinate patches. Whenever global estimates are necessary, we also need to include lower order terms arising from localizations. However, they will be controlled by the inhomogeneous counterparts of these spaces, only requiring lower regularity control.
\begin{proposition}\label{Interpolation 1}
Let $\sigma_m>-\frac{1}{2}$ and $\displaystyle m-\sigma_m-\frac{d}{2}>0$. We define $\displaystyle \theta_j=\frac{j}{m}$, $\displaystyle \frac{1}{p_j}=\frac{\theta_j}{2}$, $\displaystyle \sigma_j=\sigma_m\theta_j$. Then, for $0<j<m$, we have
\begin{align*}
    \|r^{\sigma_j}\partial^jf\|_{L^{p_j}}&\lesssim \|f\|^{1-\theta_j}_{L^{\infty}}\|r^{\sigma_m}\partial^mf\|^{\theta_j}_{L^2}.
\end{align*}
\end{proposition}
\begin{proof}
    See the proof of Proposition 2.14 in \cite{Ifrim-Tataru}.
\end{proof}
\begin{proposition}
Let $\sigma_m>-\frac{1}{2}$ and $\displaystyle m-\frac{1}{2}-\sigma_m-\frac{d}{2}>0$. We define $\displaystyle \theta_j=\frac{2j-1}{2m-1}$, $\displaystyle \frac{1}{p_j}=\frac{\theta_j}{2}$, $\displaystyle \sigma_j=\sigma_m\theta_j$. Then, for $0<j<m$, we have
\begin{align*}
    \|r^{\sigma_j}\partial^jf\|_{L^{p_j}}&\lesssim \|f\|^{1-\theta_j}_{\dot{C}^{\frac{1}{2}}}\|r^{\sigma_m}\partial^mf\|^{\theta_j}_{L^2}.
\end{align*}
\end{proposition}
\begin{proof}
    See the proof of Proposition 2.15 in \cite{Ifrim-Tataru}.
\end{proof}
\begin{proposition}\label{Interpolation 2}
Let $\sigma_m>\frac{m-2}{2}$ and $\displaystyle m-\frac{1}{2}-\sigma_m-\frac{d}{2}>0$. We define $\displaystyle \theta_j=\frac{j}{m}$, $\displaystyle \frac{1}{p_j}=\frac{\theta_j}{2}$, $\displaystyle \sigma_j=\sigma_m\theta_j-\frac{1-\theta_j}{2}$. Then, for $0<j<m$, we have
\begin{align*}
    \|r^{\sigma_j}\partial^jf\|_{L^{p_j}}&\lesssim \|f\|^{1-\theta_j}_{\tilde{C}^{0,\frac{1}{2}}}\|r^{\sigma_m}\partial^mf\|^{\theta_j}_{L^2}.
\end{align*}
\end{proposition}
\begin{proof}
    See the proof of Proposition 2.16 in \cite{Ifrim-Tataru}.
\end{proof}

\subsection{Energy and distance functional estimates}
We are going to use the following crucial results:
\begin{thm}
    For each integer $k\geq 0$, there exists an energy functional $E^{2k}$ with the followin properties
    \begin{enumerate}
        \item[a)] Coercivity: as long as $A\ll 1$, we have
        \begin{align*}
            E^{2k}(r,v)\approx \|(r,v)\|^2_{\mathcal{H}^{2k}}
        \end{align*}
        \item[b)] Energy estimates for solutions:
        \begin{align*}
             \frac{d}{dt}E^{2k}(r,v)&\lesssim_A B\|(r,v)\|^2_{\mathcal{H}^{2k}}.
        \end{align*}
        In particular, Gr\" onwall's inequality implies that
        \begin{align*}
          \|(r,v)(t)\|^2_{\mathcal{H}^{2k}}&\lesssim e^{\int_0^TC(A)B(s)\,ds} \|(r,v)(0)\|^2_{\mathcal{H}^{2k}} 
        \end{align*}
    \end{enumerate}
\end{thm}
\begin{remark}
    The analysis carried out in Section 7 of \cite{Ifrim-Tataru}, and Proposition 3.1 from the same paper show that
    \begin{equation}\label{Usual energy estimates}
    \begin{aligned}
          \|(r,v)(t)\|^2_{\mathcal{H}^{2k}}&\lesssim e^{\int_0^TC(A)B(s)\,ds} \|(r,v)(0)\|^2_{\mathcal{H}^{2k}}\\
          \|(s,w)\|^2_{\mathcal{H}}&\lesssim e^{\int_0^TC(A)B(s)\,ds} \|(s,w)(0)\|^2_{\mathcal{H}}
          \end{aligned}
        \end{equation}
        for every $k\geq 0$, not necessarily an integer. These are the energy estimates for the fully nonlinear equation and for its linearized counterpart that we shall use throughout this paper.
\end{remark}
We are also going to compare different solutions.

Let
\begin{align*}
    D_{\mathcal{H}}((r_1,v_1),(r_2,v_2))&=\int_{\Omega_t}(r_1+r_2)^{\sigma-1}((r_1-r_2)^2+\kappa(r_1+r_2)(v_1-v_2)^2)\,dx,
\end{align*}
where $\displaystyle \sigma=\frac{1}{\kappa}$ and $\displaystyle \Omega_t=\Omega_{1,t}\cap\Omega_{2,t}$.

As a consequence of Lemma 4.2 and Proposition 4.3 in \cite{Ifrim-Tataru}, we have
\begin{equation}\label{Metric estimate}
  D_{\mathcal{H}}((r_1,v_1),(r_2,v_2))(t)\lesssim e^{\int_0^TC(A_1+A_2)(B_1+B_2)(s)\,ds} D_{\mathcal{H}}((r_1,v_1),(r_2,v_2))(0) 
\end{equation}
\subsection{Potential function formulations}

We shall now recall the potential function formulations of \eqref{irrotational Euler} and \eqref{linearized irrotational Euler}.

We recall our equations:
\begin{equation*}
\begin{cases}\begin{aligned}
  &D_tr+\kappa r(\nabla\cdot v)=0\\
  &D_tv+\nabla r=0\\
  &\curl v=0,
  \end{aligned}
  \end{cases}
\end{equation*}
as well as the linearized counterparts:
\begin{equation*}
\begin{cases}\begin{aligned}
  &D_ts+w\cdot\nabla r+\kappa s(\nabla\cdot v)+\kappa r(\nabla\cdot w)=0\\
  &D_tw+(w\cdot\nabla)v+\nabla s=0\\
  &\curl w=0.
  \end{aligned}
  \end{cases}
\end{equation*}

We recall the following reformulation from \cite{Avadanei}:
\begin{proposition}
  There exist velocity functions $\phi$ and $\psi$, such that 
  \begin{align*}
      D_t\phi&=\frac{|v|^2}{2}-r\\
      D_t^2\phi&+\nabla r\cdot\nabla\phi-\kappa r\Delta\phi=0,
  \end{align*}
  and
  \begin{align*}
      D_t\psi&=-s\\
      D_t^2\psi&-\nabla r\cdot\nabla\psi-\kappa r\Delta\psi=\kappa s(\nabla\cdot v),
  \end{align*}
\end{proposition}

\section{A regularized domain defining function}
In this section we shall introduce a regularization $\tilde{r}$ of the domain defining function $r$, which is adapted in frequency to every single boundary layer of the domain $\Omega_t$. This construction will in particular enable us to transplant our construction of the counterexample from \cite{Avadanei} for general backgrounds with $A\ll 1$ and integrable $B$.

We start with a partition of unity of the form $\chi(2^{2j}r)$, where $\chi$ is supported in the region $\{\lambda\approx 1\}$. For every $j$, we have
\begin{align*}
|2^{2j}r-2^{2j}\Psi^j(r)|&=2^{2j}|(I-\Psi^j)(r)|\lesssim 2^{2j}\|(I-\Psi^j)(r-x_d)\|_{L_x^\infty}\\
&\lesssim 2^{2j}2^{-2j}\|\nabla r-e_d\|_{L_x^\infty}\lesssim A\ll 1
\end{align*}
This shows that only finitely many terms in the sum $\displaystyle\sum_{\substack{j}}\chi(2^{2j}\Psi^j(r))$ are nonzero, and that its value is bounded from below by $1-O(A)$, thus allowing us to consider a partition of unity of the form $\displaystyle\chi_j(x)=\frac{\chi(2^{2j}\Psi^j(r))}{\sum_{\substack{k}}\chi(2^{2k}\Psi^k(r))}$, with the properties 
\begin{align*}
   1&=\sum_{\substack{j}}\chi_j(x)\\
   0&=\sum_{\substack{j}}\partial_{x_d}\chi_j(x)\\
   0&=\sum_{\substack{j}}\partial^2_{x_d}\chi_j(x). 
\end{align*}

Let
\begin{align*}
\tilde{r}:=\sum_{\substack{j}}\chi_j(x)\Psi^j(r).
\end{align*}
We have the following 
\begin{lemma}\label{Regularized coordinates}
Let $\tilde{r}$ be defined as above. Then, the change of coordinates $(x',x_d)\rightarrow (x',\tilde{r}(x',x_d))$ is invertible, and $\tilde{r}$ satisfies the following estimates:
\begin{enumerate}
    \item[a)] $\displaystyle \|\nabla\tilde{r}-e_d\|_{L_x^\infty}\lesssim A\ll 1$
    \item[b)] $\displaystyle \|\tilde{r}\nabla^2\tilde{r}\|_{L_x^\infty}+\|r\nabla^2\tilde{r}\|_{L_x^\infty}\lesssim A\ll 1$
    \item[c)] $\displaystyle \|r^{-\frac{1}{2}}\nabla(r-\tilde{r})\|_{L_x^\infty}\lesssim \|\nabla r\|_{\dot{C}^{\frac{1}{2}}}$.
\end{enumerate}

\end{lemma}
\begin{proof}
We have
\begin{align*}
 \partial_{x_d}\tilde{r}&=\sum_{\substack{j}}\chi_j(x)(\Psi^j(\partial_{x_d}r)+D\Psi^j(r))+\partial_{x_d}\chi_j(x)\Psi^j(r)\\
 &=\sum_{\substack{j}}\chi_j(x)(\Psi^j(\partial_{x_d}r)+(D\Psi^j)(r))-\sum_{\substack{j}}\partial_{x_d}\chi_j(x)(I-\Psi^j)(r-x_d).   
\end{align*}

We consider the map $\displaystyle (x',x_d)\rightarrow (x',\tilde{r})$. We can see that the Jacobian of this transformation is $|\partial_{x_d}\tilde{r}|$. 
In the definition of $\tilde{r}$, it is clear that only finitely many terms in our expression for $\partial_{x_d}\tilde{r}$ are nonzero. We also have
\begin{align*}
 \partial_{x_d}\tilde{r}-1&=\sum_{\substack{j}}\chi_j(x)\Psi^j(\partial_{x_d}r-1)+\chi_j(x)(D\Psi^j)(r)-\sum_{\substack{j}}\partial_{x_d}\chi_j(x)(I-\Psi^j)(r-x_d).   
\end{align*}
There are only finitely many indices in the previous expression which are nonzero. Let $j$ be such an index.

We have
\begin{align*}
   \partial_{x_d}\chi_j(x)=\frac{\chi'(2^{2j}\Psi^j(r))2^{2j}\Psi^j(\partial_{x_d}r)}{\sum_{\substack{k}}\chi(2^{2k}\Psi^k(r))}-\frac{\chi(2^{2j}\Psi^j(r))\sum_{\substack{k}}\chi'(2^{2k}\Psi^k(r))2^{2k}\Psi^k(\partial_{x_d}r)}{\left(\sum_{\substack{k}}\chi(2^{2k}\Psi^k(r))\right)^2}.
\end{align*}
In the numerator from the second term, all indices $k$ for which the terms don't vanish satisfy $2^{2k}\approx 2^{2j}$, hence
\begin{align*}
\left|\partial_{x_d}\chi_j(x)(I-\Psi^j)(r-x_d)\right|&\lesssim (1+A)2^{2j}\|(I-\Psi^j)(r-x_d)\|_{L_x^\infty}\\
&\lesssim (1+A)2^{2j}2^{-2j}\|(I-\Psi^j)(
\nabla r-e_d)\|_{L_x^\infty}\lesssim A,\\
\left|\chi_j(x)\Psi^j(\partial_{x_d}r-1)\right|&\lesssim \|\Psi^j(
\nabla r-e_d)\|_{L_x^\infty}\lesssim A,
\end{align*}
The terms involving $D\Psi^j$ can be estimated similarly. This means that $|\partial_{x_d}\tilde{r}-1|\lesssim A\ll 1$, hence this change of coordinates is invertible. Moreover, when $\displaystyle r(x)\approx 2^{-2j}$, the only indices $j'$ for which the terms in the expression of $\tilde{r}$ are nonzero satisfy $2^{-2j'}\approx 2^{-2j}$. We write
\begin{align*}
r-\tilde{r}=\sum_{\substack{j}}\chi_j(x)(I-\Psi^j)(r)
\end{align*}
For such an index $j'$, we have
\begin{align*}
    |\chi_{j'}(x)(I-\Psi^{j'})(r)|&\lesssim|(I-\Psi^{j'})(r)|=|(I-\Psi^{j'})(r-x_d)|\lesssim 2^{-2j'}\|(I-\Psi^{j'})(\nabla r-e_d)\|_{L_x^\infty}\\
    &\lesssim 2^{-2j'}A\lesssim 2^{-2j}A
\end{align*}
As $A\ll 1$, this shows that $r\approx 2^{-2j}$ if and only if $\tilde{r}\approx 2^{-2j}$. We may therefore use $r$ and $\tilde{r}$ interchangeably.

For $\displaystyle \partial^2_{x_d}\tilde{r}$, we write
\begin{align*}
    \partial_{x_d}^2\tilde{r}&=\sum_{\substack{j}}\chi_j(x)\Psi^j(\partial^2_{x_d}r)+2\sum_{\substack{j}}\partial_{x_d}\chi_j(x)\Psi^j(\partial_{x_d}r)+\sum_{\substack{j}}\partial^2_{x_d}\chi_j(x)\Psi^j(r)\\
    &=\sum_{\substack{j}}\chi_j(x)\Psi^j(\partial^2_{x_d}r)-2\sum_{\substack{j}}\partial_{x_d}\chi_j(x)(I-\Psi^j)(\partial_{x_d}r)-\sum_{\substack{j}}\partial^2_{x_d}\chi_j(x)(I-\Psi^j)(r)\\
    &=\sum_{\substack{j}}\chi_j(x)\Psi^j(\partial^2_{x_d}r)-2\sum_{\substack{j}}\partial_{x_d}\chi_j(x)(I-\Psi^j)(\partial_{x_d}r-1)-\sum_{\substack{j}}\partial^2_{x_d}\chi_j(x)(I-\Psi^j)(r-x_d)
\end{align*}

When $\tilde{r}
\approx 2^{-2j'}$, all indices $j$ for which the terms in the previous sum are nonzero satisfy $2^{-2j}\approx 2^{-2j'}$.

Let $j$ be such an index. We note that
\begin{align*}
 \partial^2_{x_d}\chi_j(x)&=\frac{\chi''(2^{2j}\Psi^j(r))(2^{2j}\Psi^j(\partial_{x_d}r))^2+\chi'(2^{2j}\Psi^j(r))2^{2j}\Psi^j(\partial^2_{x_d}r)}{\sum_{\substack{k}}\chi(2^{2k}\Psi^k(r))}\\
 &-\frac{\chi'(2^{2j}\Psi^j(r))2^{2j}\Psi^j(\partial_{x_d}r)\sum_{\substack{k}}\chi'(2^{2k}\Psi^k(r))2^{2k}\Psi^k(\partial_{x_d}r)}{\left(\sum_{\substack{k}}\chi(2^{2k}\Psi^k(r))\right)^2}\\
 &-\frac{\chi'(2^{2j}\Psi^j(r))\Psi^j(\partial_{x_d}r)\sum_{\substack{k}}\chi'(2^{2k}\Psi^k(r))2^{2k}\Psi^k(\partial_{x_d}r)}{\left(\sum_{\substack{k}}\chi(2^{2k}\Psi^k(r))\right)^2}\\
 &-\frac{\chi(2^{2j}\Psi^j(r))\sum_{\substack{k}}\chi''(2^{2k}\Psi^k(r))(2^{2k}\Psi^k(r))^2}{\left(\sum_{\substack{k}}\chi(2^{2k}\Psi^k(r))\right)^2}\\
 &-\frac{\chi'(2^{2k}\Psi^k(r))2^{2k}\Psi^k(\partial^2_{x_d}r)}{\left(\sum_{\substack{k}}\chi(2^{2k}\Psi^k(r))\right)^2}+\frac{2\chi(2^{2j}\Psi^j(r))\left(\sum_{\substack{k}}\chi'(2^{2k}\Psi^k(r))2^{2k}\Psi^k(\partial_{x_d}r)\right)^2}{\left(\sum_{\substack{k}}\chi(2^{2k}\Psi^k(r))\right)^3},
\end{align*}
and that all indices corresponding to non-vanishing terms in the numerators satisfy $2^{2k}\approx 2^{2j}\approx 2^{2j'}$. Going back into the expansion of $\partial_{x_d}\tilde{r}$, it immediately follows that all terms essentially fall into four categories, which can be estimated as follows (here we also use the fact that the denominator is positive, and uniformly bounded away from $0$): 
\begin{align*}
  |\chi''(2^{2j}\Psi^j(r))&2^{2j}\Psi^j(\partial_{x_d}r)2^{2j}\Psi^j(\partial_{x_d}r)(I-\Psi^j)(r-x_d)|\lesssim 2^{2j}(1+O(A))^2A\approx 2^{2j'}A\\
  |\chi'(2^{2j}\Psi^j(r))&2^{2j}\Psi^j(\partial^2_{x_d}r)(I-\Psi^j)(r-x_d)|\lesssim 2^{2j}(1+O(A))A\approx 2^{2j'}A \\
  |\chi'(2^{2j}\Psi^j(r))&2^{2j}\Psi^j(\partial_{x_d}r)(I-\Psi^j)(\partial_{x_d}r-1)|\lesssim 2^{2j}(1+O(A))A\approx 2^{2j'}A\\
  |\chi(2^{2j}\Psi^j(r))&\Psi^j(\partial^2_{x_d}r)|=|\chi(2^{2j}\Psi^j(r))\Psi^j(\partial_{x_d}(\partial_{x_d}r-1))|\lesssim 2^{2j}A\approx 2^{2j'}A.
\end{align*}
This implies that
\begin{align*}
 \|\tilde{r}\partial^2_{x_d}\tilde{r}\|_{L_x^\infty((r\approx 2^{-2j'}))}&\lesssim A.   
\end{align*}
The estimates for the other derivatives follow similarly.

This finishes the proof.

\end{proof}

\section{Gallery modes in modified coordinates}
In this section we shall transplant our construction for the counterexamples from \cite{Avadanei}. We do so by essentially carrying out two steps at the same time: we use the regularized function $\tilde{r}$ as an approximate normal coordinate which captures the geometrical condition $A\ll 1$, as well as the forward flow of $v$ in order to account for the Galilean transport of our solution.

We recall the following construction from \cite{Avadanei}:
Let $\displaystyle a\in C_0^\infty(\mathbb{R}^d)$ be an even function with 
\begin{align*}\supp\, a \in\left\{\tilde{\xi'}|1-\varepsilon(\kappa)\leq|\tilde{\xi'}|\leq 1+\varepsilon(\kappa)\right\}.
\end{align*}

For $j\geq 0$, we define
\begin{align*}
        W^j(t,y)&=\frac{1}{(2\pi)^{d-1}}\int e^{it|\xi'|^{\frac{1}{2}}+iy'\cdot\xi'}B\left(0,|\xi'|y_d\right)a(2^{-2j}\xi')\,d\xi'
    \end{align*}
We note that
\begin{align*}
    W^j(t,y)&=\frac{1}{(2\pi)^{d-1}}\int e^{it|\xi'|^{\frac{1}{2}}+iy'\cdot\xi'-|\xi'|y_d}a(2^{-2j}\xi')\,d\xi'\\
        &=\frac{2^{2j(d-1)}}{(2\pi)^{d-1}}\int e^{it2^j|\eta|^{\frac{1}{2}}+i2^{2j}y'\cdot\eta-2^{2j}|\eta|y_d}a(\eta)\,d\eta\\
        &=2^{2j(d-1)}W^{0}(2^jt,2^{2j}y)
\end{align*}

It is clear that $W^j$ is an exact solution to our wave equation

\begin{align*}
    (\partial_t^2-\kappa y_d\Delta_{y'}-\kappa y_d\partial_d^2-\partial_d)W^j&=0.
\end{align*}

Let $\displaystyle \tilde{W}^j(t,x)=W^j(t,\Phi(x))$, where $\displaystyle \Phi(x)=(x',\tilde{r}(x))$. We have
\begin{align*}
    (\kappa \tilde{r}\Delta+\nabla\tilde{r}\cdot\nabla)\tilde{W}^j&=\kappa\tilde{r}\Delta_{x',\tilde{r}}W^j+\partial_{\tilde{r}}W^j\\
    &+2\kappa\tilde{r}\nabla_{x'}\partial_{\tilde{r}}W^j\cdot\nabla_{x'}\tilde{r}+\kappa\tilde{r}\partial^2_{\tilde{r}}W^j|\nabla_{x'}\tilde{r}|^2+\kappa\tilde{r}\partial^2_{\tilde{r}}W^j((\partial_{x_d}\tilde{r})^2-1)+\kappa\tilde{r}\partial^2_{x_d}\tilde{r}\partial_{\tilde{r}}W^j\\
    &+\nabla_{x'}W^j\cdot \nabla_{x'}\tilde{r}+\partial_{\tilde{r}}W^j|\nabla_{x'}\tilde{r}|^2+\partial_{\tilde{r}}W^j((\partial_{x_d}\tilde{r})^2-1)+\kappa\tilde{r}\Delta_{x'}\tilde{r}\partial_{\tilde{r}}W^j,
\end{align*}
which is in turn equivalent to
\begin{align*}
    (\kappa r\Delta+\nabla r\cdot\nabla)\tilde{W}^j&=\partial_t^2W^j+\kappa (r-\tilde{r})\Delta W^j+\nabla(r-\tilde{r})\cdot \nabla W^j\\
    &+2\kappa\tilde{r}\nabla_{x'}\partial_{\tilde{r}}W^j\cdot\nabla_{x'}\tilde{r}+\kappa\tilde{r}\partial^2_{\tilde{r}}W^j|\nabla_{x'}\tilde{r}|^2+\kappa\tilde{r}\partial^2_{\tilde{r}}W^j((\partial_{x_d}\tilde{r})^2-1)+\kappa\tilde{r}\partial^2_{x_d}\tilde{r}\partial_{\tilde{r}}W^j\\
    &+\nabla_{x'}W^j\cdot \nabla_{x'}\tilde{r}+\partial_{\tilde{r}}W^j|\nabla_{x'}\tilde{r}|^2+\partial_{\tilde{r}}W^j((\partial_{x_d}\tilde{r})^2-1)+\kappa\tilde{r}\Delta_{x'}\tilde{r}\partial_{\tilde{r}}W^j\\
    &:=\partial_t^2W^j+H^j.
\end{align*}
Lemma \ref{Regularized coordinates} implies that that
\begin{align*}
    \|y_d^{\frac{1}{2\kappa}-\frac{1}{2}}H^j\|_{L^2}&\lesssim A(\|y_d^{\frac{1}{2\kappa}+\frac{1}{2}}\nabla^2_yW^j\|_{L^2}+\|y_d^{\frac{1}{2\kappa}-\frac{1}{2}}\nabla_yW^j\|_{L^2}).
\end{align*}
By setting $\displaystyle Y=2^{2j}y$ and $\displaystyle \tau=2^jt$, we have $\displaystyle \|y_d^{\frac{1}{2\kappa}+\frac{1}{2}}\nabla^2_yW^j\|_{L^2}=2^{j\left(d+1-\frac{1}{\kappa}\right)}\|Y_d^{\frac{1}{2\kappa}+\frac{1}{2}}\nabla^2_YW^0\|_{L^2}$ and $\displaystyle \|y_d^{\frac{1}{2\kappa}-\frac{1}{2}}\nabla_yW^j\|_{L^2}=2^{j\left(d+1-\frac{1}{\kappa}\right)}\|Y_d^{\frac{1}{2\kappa}-\frac{1}{2}}\nabla_YW^0\|_{L^2}$. For $\displaystyle c\in(0,1)$, we have
\begin{align*}
    \sup_{0\leq\tau\leq c}\left(\|Y_d^{\frac{1}{2\kappa}+\frac{1}{2}}\nabla^2_YW^0\|_{L^2}+\|Y_d^{\frac{1}{2\kappa}-\frac{1}{2}}\nabla_YW^0\|_{L^2}\right)<\infty
\end{align*}
This implies that
\begin{align*}
    \sup_{0\leq t\leq c2^{-j}}\left(\|Y_d^{\frac{1}{2\kappa}+\frac{1}{2}}\nabla^2_YW^j\|_{L^2}+\|Y_d^{\frac{1}{2\kappa}-\frac{1}{2}}\nabla_YW^j\|_{L^2}\right)\lesssim 2^{j\left(d+1-\frac{1}{\kappa}\right)}, 
\end{align*}
hence
\begin{align*}
    \sup_{0\leq t\leq c2^{-j}}\|y_d^{\frac{1}{2\kappa}-\frac{1}{2}}H^j\|_{L^2}\lesssim A2^{j\left(d+1-\frac{1}{\kappa}\right)}.
\end{align*}
It is also clear that
\begin{align*}
    (\partial_t^2-\kappa r\Delta-\nabla r\cdot\nabla)\tilde{W}^j=-H^j.
\end{align*}

Let $X_t$ be the forward flow associated to $v$, in the sense that $\displaystyle \frac{d}{dt}X_t(a)=v(t,X_t(a))$, and $\displaystyle X_0(a)=a$. Let $\displaystyle \Phi_t(a):=\Phi(X_t^{-1}(a))$, so that in particular, $\Phi_t(X_t(a))=\Phi(a)$. We define $\displaystyle W_\ast^j(t,x)=\tilde{W}^j(t,X^{-1}_t(x))=W^j(t,\Phi_t(x))$. In particular, $\displaystyle \tilde{W}^j(t,a)=W_\ast^j(t,X_t(a))$.

Let $\displaystyle J_t:=\det DX_t$, $\displaystyle m(t):=J_t^\kappa$, and $\displaystyle G(t):=(DX_t)^{-1}(DX_t)^{-T}$.

We have 
\begin{align*}
\frac{d}{dt}r(t,X_t(a))=(D_tr)(t,X_t(a))=-\kappa r(\nabla\cdot v).
\end{align*}
We also know that $\displaystyle \frac{d}{dt}J_t(a)=\frac{d}{dt}DX_t(a)=(\nabla \cdot v)(t,X_t(a))J_t(a)$. It follows that
\begin{align*}
    \frac{d}{dt}(J_t^\kappa r(t,X_t(a)))&=\kappa J_t^{\kappa-1}J_t(\nabla_x\cdot v)(t,X_t(a))r(t,X_t(a))-\kappa J_t^\kappa r(\nabla\cdot v)=0.
\end{align*}
Thus, in particular, $\displaystyle r(t,X_t(a))=J_t^{-\kappa}r_0(a)$, hence $\displaystyle r^{\frac{1}{\kappa}}(t,X_t(a))=J_t^{-1}r^{\frac{1}{\kappa}}_0(a)$ and \newline $\displaystyle r^{\frac{1}{\kappa}-1}(t,X_t(a))=J_t^{\kappa-1}r^{\frac{1}{\kappa}-1}_0(a)$. 

We now have
\begin{align*}
    \partial_t(m(t)\partial_t\tilde{W}^j(t,a))&=\partial_t(m(t))\partial_t\tilde{W}^j(t,a)+m(t)\partial^2_t\tilde{W}^j(t,a)\\
    &=\kappa m(t)(\nabla\cdot v)(t,X_t(a))(D_t W^j_\ast)(t,X_t(a))+(D^2_tW^j_\ast)(t,X_t(a))m(t)\\
    &=m(t)(D^2_tW^j_\ast+\kappa(\nabla\cdot v)D_tW^j_\ast)(t,X_t(a))
\end{align*}

We also have
\begin{align*}
    r^{\frac{1}{\kappa}-1}(\kappa r\Delta_xW^j_\ast&+\nabla_x r\cdot \nabla_xW^j_\ast)(t,X_t(a))\\
    &=\kappa \text{div}_x(r^{\frac{1}{\kappa}}\nabla_xW^j_\ast)(t,X_t(a))\\
    &=J_t^{-1}\kappa \text{div}_a(J_t(a)r^{\frac{1}{\kappa}}(t,X_t(a))(DX_t)^{-1}(a)(DX_t)^{-T}(a)\nabla_a\tilde{W}^j(t,a))\\
    &=\kappa J_t^{-1}(a)\text{div}_a(r_0^{\frac{1}{\kappa}}(a)G(t,a)\nabla_a\tilde{W}^j(t,a))
\end{align*}
This implies that
\begin{align*}
    r_0^{\frac{1}{\kappa}-1}(a)m(t)J_t^{-1}(a)(\kappa r\Delta_xW^j_\ast&+\nabla_x r\cdot \nabla_xW^j_\ast)(t,X_t(a))\\
    &=\kappa J_t^{-1}(a)\text{div}_a(r_0^{\frac{1}{\kappa}}(a)G(t,a)\nabla_a\tilde{W}^j(t,a)),
\end{align*}
hence
\begin{align*}
    m(t)(\kappa r\Delta_xW_\ast^j&+\nabla_x r\cdot \nabla_xW_\ast^j)(t,X_t(a))\\
    &=\kappa r_0^{1-\frac{1}{\kappa}}(a)\text{div}_a(r_0^{\frac{1}{\kappa}}(a)G(t,a)\nabla_a\tilde{W}^j(t,a))\\
    &=(\kappa r_0\Delta_a\tilde{W}^j(a)+\nabla_a r_0\cdot \nabla_a\tilde{W}^j(a))\\
    &+\kappa r_0^{1-\frac{1}{\kappa}}(a)\text{div}_a(r_0^{\frac{1}{\kappa}}(a)(G(t,a)-I)\nabla_a\tilde{W}^j(t,a)).
\end{align*}
By combining the two equations, we get that
\begin{align*}
    m(t)(D^2_tW_\ast^j&+\kappa(\nabla\cdot v)D_tW_\ast^j-\kappa r\Delta_xW_\ast^j-\nabla_x r\cdot \nabla_xW_\ast^j)(t,X_t(a))\\
    &=\partial_t(m(t)\partial_t\tilde{W}^j(t,a))-(\kappa r_0\Delta_a\tilde{W}^j(a)+\nabla_a r_0\cdot \nabla_a\tilde{W}^j(a))\\
    &-\kappa r_0^{1-\frac{1}{\kappa}}(a)\text{div}_a(r_0^{\frac{1}{\kappa}}(a)(G(t,a)-I)\nabla_a\tilde{W}^j(t,a))\\
    &=\partial^2_t\tilde{W}^j(t,a))-(\kappa r_0\Delta_a\tilde{W}^j(a)+\nabla_a r_0\cdot \nabla_a\tilde{W}^j(a))\\
    &+\partial_t((m(t)-1)\partial_t\tilde{W}^j(t,a))-\kappa r_0^{1-\frac{1}{\kappa}}(a)\text{div}_a(r_0^{\frac{1}{\kappa}}(a)(G(t,a)-I)\nabla_a\tilde{W}^j(t,a))\\
     &=-H^j(t,a)+\partial_t((m(t)-1)\partial_t\tilde{W}^j(t,a))-\kappa r_0^{1-\frac{1}{\kappa}}(a)\text{div}_a(r_0^{\frac{1}{\kappa}}(a)(G(t,a)-I)\nabla_a\tilde{W}^j(t,a))
\end{align*}
We also define $\displaystyle r_j=\Psi^jr$, $\displaystyle \Phi^j(x)=(x',r^j(x))$, as well as an approximation of the initial data $\tilde{W}^j(0)$ tailored for the $2^{-2j}$ boundary layer and spatial frequency $2^{2j}$, given by 
\begin{align*}
   W^{j,\sharp}(t,x)=W^j(t,\Phi^j(x)). 
\end{align*} 
We then have the following
\begin{lemma}\label{Boundary layer approximation}
  Let $2k>2k_0$. Under the previous notations, we have
  \begin{align*}
      \|(\partial_tW^{j,\sharp}(0),\nabla W^{j,\sharp}(0))\|_{\mathcal{H}_r^{2s}}&\lesssim 2^{j\left(2s+d-\frac{1}{\kappa}\right)}\\
      \|(\partial_tW^{j,\sharp}(0)-\partial_t\tilde{W}^j(0),\nabla (W^{j,\sharp}(0)-\tilde{W}^j(0)))\|_{\mathcal{H}_r}&\lesssim A2^{j\left(d-\frac{1}{\kappa}\right)}
  \end{align*}
\end{lemma}
\begin{proof}
    From \cite{Ifrim-Tataru}, we know that
    \begin{align*}
        r_j(x)&=\int_{\Omega}K(x,y)r(y)\,dy
    \end{align*}
    The kernels satisfy the bounds
    \begin{align*}
        |\partial_x^{\beta}K^j(x,y)|&\lesssim (2^{-2j}+2^{-j}r(x)^{\frac{1}{2}})^{-d-|\beta|}.
    \end{align*}
    These in turn imply that
    \begin{align*}
       \int_\Omega|\partial_x^{\beta}K^j(x,y)|\,dy&\lesssim (2^{-2j}+2^{-j}r(x)^{\frac{1}{2}})^{-|\beta|}\\
       \int_\Omega|\partial_x^{\beta}K^j(x,y)||x-y|\,dy&\lesssim (2^{-2j}+2^{-j}r(x)^{\frac{1}{2}})^{1-|\beta|}.
    \end{align*}
    Now, since $\Psi^j$ reproduces affine functions, we have
    \begin{align*}
        r_j(x)-r(x)&=\int_\Omega K^j(x,y)((r(y)-y_d)-(r(x)-x_d))\,dy.
    \end{align*}
    It then follows that
    \begin{align*}
        |r_j(x)-r(x)|&\lesssim A(2^{-2j}+2^{-j}r(x)^{\frac{1}{2}})\\
        |\nabla r_j(x)-e_d|+|\nabla(r_j(x)-r(x))|&\lesssim A\\
        |D^mr_j(x)|&\lesssim A(2^{-2j}+2^{-j}r(x)^{\frac{1}{2}})^{1-m}\lesssim A2^{2j(m-1)},\text{ for }m\geq 2.
    \end{align*}

    We now recall that
    \begin{align*}
        |\tilde{r}(x)-r(x)|&\lesssim A\tilde{r}(x)\simeq Ar(x)\\
        |\nabla(\tilde{r}(x)-r(x))|&\lesssim A
    \end{align*}
    Thus,
    \begin{align*}
      |\tilde{r}(x)-r_j(x)| &\lesssim A(2^{-2j}+r(x))\\
      |\nabla(\tilde{r}(x)-r_j(x))| &\lesssim A.
    \end{align*}
In particular, we also have the estimate
\begin{align*}
    2^{2j}|\tilde{r}(x)-r_j(x)| &\lesssim A(1+2^{2j}\tilde{r}(x)) 
\end{align*}
In particular, if we a priori choose $A\ll 1$, we have
\begin{align*}
    2^{2j}r_j(x) &\geq \frac{2^{2j}\tilde{r}(x)}{2}-\frac{1}{4}. 
\end{align*}
    We recall that
    \begin{align*}
        W^j(t,y)&=\frac{1}{(2\pi)^{d-1}}\int e^{it|\xi'|^{\frac{1}{2}}+iy'\cdot\xi'-|\xi'|y_d}a(\xi')\,d\xi'
    \end{align*}
    It follows that we have the bounds
    \begin{align*}
        |\partial^l_\tau\partial^\beta_YW^0(0,2^{2j}x',2^{2j}r_j(x)))|&\lesssim_{l,\beta,m}\langle 2^{2j}x'\rangle^{-m}e^{-c2^{2j}\tilde{r}(x)},
    \end{align*}
    for some constant $c>0$. 
    
    We now prove the higher order derivatives. For a generic function $f$, Fa\`a di Bruno's generalized chain rule implies that
   \begin{align*}
       D^m(f\circ(2^{2j}\Phi^j))&=\sum_{\substack{j_i\geq 2,\sum(j_i-1)=m-q\\
       0\leq p\leq q,q\geq 1}}D^qf\circ(2^{2j}\Phi^j)2^{2j(q-p)}(D\Phi^j)^{q-p}\prod_{\substack{i=1}}^p2^{2j}D^{j_i}\Phi
   \end{align*}
  
  Together with our previous bounds, this implies that
   \begin{align*}
      |\partial_x^\beta\partial_t^lW^{j,\sharp}(0,x)|&\lesssim_{l,\beta,m} 2^{j(2d-2+l+2|\beta|)}\langle 2^{2j}x'\rangle^{-M}e^{-c2^{2j}\tilde{r}(x)}.
   \end{align*}
  In this case, it follows that for every nonnegative integer $\alpha$,
  \begin{align*}
    \|r^{\alpha+\frac{1}{2\kappa}-\frac{1}{2}}\partial_tW^{j,\sharp}(0,x)\|_{L^2(dx)}&\lesssim 2^{j(2d-1)}\|r^{\alpha+\frac{1}{2\kappa}-\frac{1}{2}}\langle 2^{2j}x'\rangle^{-M}e^{-c2^{2j}\tilde{r}(x)}\|_{L^2(dx)}\\
    &\lesssim 2^{j(2d-1)}\|\tilde{r}^{\alpha+\frac{1}{2\kappa}-\frac{1}{2}}\langle 2^{2j}x'\rangle^{-M}e^{-c2^{2j}\tilde{r}(x)}\|_{L^2(dx)}\\
    &\lesssim 2^{j\left(d-\frac{1}{\kappa}+2(|\beta|-\alpha)\right)}\|Y_d^{\alpha+\frac{1}{2\kappa}-\frac{1}{2}}\langle Y'\rangle^{-M}e^{-cY_d}\|_{L^2(dY)}\\
    &\lesssim 2^{j\left(d-\frac{1}{\kappa}+2(|\beta|-\alpha)\right)}
  \end{align*}
  At the last step we made the change of variables $\displaystyle (Y',Y_d)=(2^{2j}x',2^{2j}\tilde{r}(x))$.
  We similarly get that
  \begin{align*}
    \|r^{\alpha+\frac{1}{2\kappa}}\nabla_xW^{j,\sharp}(0,x)\|_{L^2(dx)}&\lesssim  2^{j\left(d-\frac{1}{\kappa}+2(|\beta|-\alpha)\right)}
  \end{align*}
   When $k$ is an integer, in an $\mathcal{H}^{2k}$ norm we have $|\beta|-\alpha\leq k$, it follows that 
   \begin{align*}
   \|(\partial_tW^{j,\sharp}(0),\nabla W^{j,\sharp}(0))\|_{\mathcal{H}_r^{2k}}&\lesssim 2^{j\left(2k+d-\frac{1}{\kappa}\right)}.
   \end{align*}
   For fractional $s$, the conclusion follows by interpolation.

   Finally, we prove the difference bounds. The fundamental theorem of calculus, together with the estimate $\displaystyle 2^{2j}|r_j(x)-\tilde{r}(x)|\lesssim A(1+2^{2j}\tilde{r}(x))$, give that
   \begin{align*}
        |\partial^l_\tau\partial^\beta_YW^0(0,2^{2j}x',2^{2j}r_j(x)))-\partial^l_\tau\partial^\beta_YW^0(0,2^{2j}x',2^{2j}\tilde{r}(x)))|&\lesssim_{l,\beta,m}A(1+2^{2j}\tilde{r}(x))\langle 2^{2j}x'\rangle^{-m}e^{-c2^{2j}\tilde{r}(x)}.
   \end{align*}
   We then have the bound
   \begin{align*}
    \|r^{\frac{1}{2\kappa}-\frac{1}{2}}(\partial_tW^{j,\sharp}(0,x)-\partial_t\tilde{W}^j(0,x))\|_{L^2(dx)}&\lesssim A2^{j(2d-1)}\|(1+2^{2j}\tilde{r}(x))r^{\frac{1}{2\kappa}-\frac{1}{2}}\langle 2^{2j}x'\rangle^{-M}e^{-c2^{2j}\tilde{r}(x)}\|_{L^2(dx)}\\
    &\lesssim A2^{j(2d-1)}\|(1+2^{2j}\tilde{r}(x))\tilde{r}^{\alpha+\frac{1}{2\kappa}-\frac{1}{2}}\langle 2^{2j}x'\rangle^{-M}e^{-c2^{2j}\tilde{r}(x)}\|_{L^2(dx)}\\
    &\lesssim 2^{j\left(d-\frac{1}{\kappa})\right)}\|Y_d^{\alpha+\frac{1}{2\kappa}-\frac{1}{2}}\langle Y'\rangle^{-M}e^{-cY_d}\|_{L^2(dY)}\\
    &\lesssim 2^{j\left(d-\frac{1}{\kappa})\right)}
  \end{align*}
  We similarly deduce that
   \begin{align*}
    \|r^{\frac{1}{2\kappa}}(\nabla W^{j,\sharp}(0,x)-\nabla\tilde{W}^j(0,x))\|_{L^2(dx)}&\lesssim 2^{j\left(d-\frac{1}{\kappa}\right)},
  \end{align*}
  hence
  \begin{align*}
      \|(\partial_tW^{j,\sharp}(0)-\partial_t\tilde{W}^j(0),\nabla (W^{j,\sharp}(0)-\tilde{W}^j(0)))\|_{\mathcal{H}_r}&\lesssim A2^{j\left(d-\frac{1}{\kappa}\right)},
  \end{align*}
  as desired.
\end{proof}
We now consider the solution $\theta^j$ to the problem

\begin{align*}
(D_t^2&+\kappa (\nabla\cdot v)D_t-\kappa r\Delta-\nabla r\cdot\nabla)\theta^j(t,X_t(a))\\
&=m(t)^{-1}(H^j(t,a)-\partial_t((m(t)-1)\partial_t\tilde{W}^j(t,a))+\kappa r_0^{1-\frac{1}{\kappa}}(a)\text{div}_a(r_0^{\frac{1}{\kappa}}(a)(G(t,a)-I)\nabla_a\tilde{W}^j(t,a)))\\
\theta^j_{t=0}&=W^{j,\sharp}(0)-W^j_\ast(0),D_t\theta^j_{t=0}=\partial_tW^{j,\sharp}(0)-D_tW^j_\ast(0).
\end{align*}
We have the following
\begin{lemma}
 If $\theta^j$ is a solution to
 \begin{align*}
(D_t^2&+\kappa (\nabla\cdot v)D_t-\kappa r\Delta-\nabla r\cdot\nabla)\theta^j(t,X_t(a))\\
&=m(t)^{-1}(H^j(t,a)-\partial_t((m(t)-1)\partial_t\tilde{W}^j(t,a))+\kappa r_0^{1-\frac{1}{\kappa}}(a)\text{div}_a(r_0^{\frac{1}{\kappa}}(a)(G(t,a)-I)\nabla_a\tilde{W}^j(t,a)))\\
\theta^j_{t=0}&=W^{j,\sharp}(0)-W^j_\ast(0),D_t\theta^j_{t=0}=\partial_tW^{j,\sharp}(0)-D_tW^j_\ast(0).
\end{align*}
then there exists a sequence $\displaystyle (\delta_j)_{j\geq 0}\subset(0,\infty)$ with $\displaystyle\lim_{\substack{j\rightarrow\infty}}\delta_j=0$, such that for large enough $j$, 
\begin{align*}
    \sup_{0\leq t\leq c2^{-j}}\|(D_t\theta^j,\nabla\theta^j)\|_{\mathcal{H}}&\lesssim (A+\delta_j)2^{j\left(d-\frac{1}{\kappa}\right)}.
\end{align*}
\end{lemma}
\begin{proof}
    We have the following:
    \begin{align*}
\frac{d}{dt}\int_{\Omega_t} r^{\frac{1}{\kappa}-1}(|D_t\theta^j|^2+\kappa r|\nabla\theta^j|^2)\,dx&=\int_{\Omega_t} D_t(r^{\frac{1}{\kappa}-1}(|D_t\theta^j|^2+\kappa r|\nabla\theta^j|^2))\,dx\\
&+\int_{\Omega_t} (\nabla\cdot v)r^{\frac{1}{\kappa}-1}(|D_t\theta^j|^2+\kappa r|\nabla\theta^j|^2)\,dx\\
&=2\int_{\Omega_t} r^{\frac{1}{\kappa}-1}(D_t\theta^jD_t^2\theta^j+\kappa r\nabla\theta^j\cdot\nabla D_t\theta^j)\,dx\\
&-\left(1-\kappa\right)\int_{\Omega_t} (\nabla\cdot v)r^{\frac{1}{\kappa}-1}|D_t\theta^j|^2\,dx\\
&-\kappa\int_{\Omega_t} (\nabla\cdot v)r^{\frac{1}{\kappa}}|\nabla\theta^j|^2\,dx-2\kappa\int_{\Omega_t} r^{\frac{1}{\kappa}}\nabla\theta^j\cdot\nabla v\nabla\theta^j\,dx\\
&+\int_{\Omega_t} (\nabla\cdot v)r^{\frac{1}{\kappa}-1}(|D_t\theta^j|^2+\kappa r|\nabla\theta^j|^2)\,dx\\
&=2\int_{\Omega_t} r^{\frac{1}{\kappa}-1}D_t\theta^j(D_t^2\theta^j+\kappa (\nabla\cdot v)D_t\theta^j-\kappa r\Delta\theta^j-\nabla r\cdot\nabla \theta^j)\,dx\\
&+II(t),
    \end{align*}
    where
    \begin{align*}
        II(t)&=-\kappa \int_{\Omega_t} r^{\frac{1}{\kappa}-1}(\nabla\cdot v)(D_t\theta^j)^2\,dx-2\kappa\int_{\Omega_t} r^{\frac{1}{\kappa}}\nabla\theta^j\cdot\nabla v\nabla\theta^j\,dx.
    \end{align*}
    In particular, we deduce the estimate
    \begin{align*}
        |II(t)|&\lesssim \|\nabla v(t)\|_{L_x^\infty}\|(D_t\theta^j,\nabla\theta^j)(t)\|^2_{\mathcal{H}}.
    \end{align*}
    We now define $\displaystyle \tilde{\theta}^j(t,a)=\theta^j(t,X_t(a))$. In particular, $\displaystyle \partial_t\tilde{\theta}^j(t,a)=(D_t\theta^j)(t,X_t(a))$. Let 
    \begin{align*}
        I(t)&=2\int_{\Omega_t} r^{\frac{1}{\kappa}-1}D_t\theta^j(D_t^2\theta^j+\kappa (\nabla\cdot v)D_t\theta^j-\kappa r\Delta\theta^j-\nabla r\cdot\nabla \theta^j)\,dx\\
        &=2\int_{\Omega_0} J_t(a)^{\kappa-1}m(t)^{-1}J_t(a)r_0^{\frac{1}{\kappa}-1}\partial_t\tilde{\theta}^jH^j\,da\\
        &-2\int_{\Omega_0}J_t(a)^{\kappa-1} m(t)^{-1}J_t(a)r_0^{\frac{1}{\kappa}-1}\partial_t\tilde{\theta}^j\partial_tm(t)\partial_t\tilde{W}^j(t,a)\,da\\
        &-2\int_{\Omega_0} J_t(a)^{\kappa-1}m(t)^{-1}J_t(a)r_0^{\frac{1}{\kappa}-1}\partial_t\tilde{\theta}^j(m(t)-1)\partial^2_t\tilde{W}^j(t,a)\,da\\
        &+2\kappa\int_{\Omega_0}J_t(a)^{\kappa-1} m(t)^{-1}J_t(a)\partial_t\tilde{\theta}^j\text{div}_a(r_0^{\frac{1}{\kappa}}(a)(G(t,a)-I)\nabla_a\tilde{W}^j(t,a)))\,da\\
        &=I_a(t)+I_b(t)+I_c(t)+I_d(t)
    \end{align*}

    Let $\displaystyle \delta_j=\int_0^{c2^-j}\|\nabla v(t)\|\,dt$. In particular, 
    \begin{align*}
        \left|\int_0^tII(s)\,ds\right|&\lesssim\delta_j \sup_{\substack{0\le t\leq c2^{-j}}}\|(\partial_t\tilde{\theta}^j,\nabla\tilde{\theta}^j)\|^2_{\mathcal{H}}
    \end{align*}

    We note that when $T\leq c2^{-j}$ and $j$ is large enough, we immediately get that
    \begin{align*}
       \|m(t)-1\|_{L_a^\infty} +\|G(t)-I\|_{L_a^\infty}&\lesssim \delta_j\\
       \|m_t(t)\|_{L_a^\infty}+\|\partial_tG(t)\|_{L_a^\infty}&\lesssim \|\nabla v(t)\|_{L_a^\infty}.
    \end{align*}
    Moreover, our construction of $W^j$ guarantees that
    \begin{align*}
        \|(\partial_tW^j_\ast,\nabla_aW^j_\ast)\|_{\mathcal{H}}&\lesssim 2^{j\left(d-\frac{1}{\kappa}\right)}\\
        \|(\partial^2_tW^j_\ast,\nabla_a\partial_tW^j_\ast)\|_{\mathcal{H}}&\lesssim 2^{j\left(d+1-\frac{1}{\kappa}\right)}.
    \end{align*}
    For every $t\in[0,T_j]$, we then have the bounds
    \begin{align*}
        \left|\int_0^tI_a(s)\,ds\right|&\lesssim \|(\partial_t\tilde{\theta}^j,\nabla\tilde{\theta}^j)\|_{\mathcal{H}}c2^{-j}A2^{j\left(d+1-\frac{1}{\kappa}\right)}\lesssim cA2^{j\left(d-\frac{1}{\kappa}\right)}\|(\partial_t\tilde{\theta}^j,\nabla\tilde{\theta}^j)\|_{\mathcal{H}}\\
        \left|\int_0^tI_b(s)\,ds\right|&\lesssim \|(\partial_t\tilde{\theta}^j,\nabla\tilde{\theta}^j)\|_{\mathcal{H}}\delta_j2^{j\left(d-\frac{1}{\kappa}\right)}\\
        \left|\int_0^tI_c(s)\,ds\right|&\lesssim 2^j\|(\partial_t\tilde{\theta}^j,\nabla\tilde{\theta}^j)\|_{\mathcal{H}}c2^{-j}\delta_j2^{j\left(d-\frac{1}{\kappa}\right)}\lesssim \|(\partial_t\tilde{\theta}^j,\nabla\tilde{\theta}^j)\|_{\mathcal{H}}\delta_j2^{j\left(d-\frac{1}{\kappa}\right)}.
    \end{align*}
    We are left to estimate $\displaystyle I_d$. We integrate by parts once in space and we obtain
    \begin{align*}
        I_d(t)&=2\kappa\int_{\Omega_0} \partial_t\tilde{\theta}^j\text{div}_a(r_0^{\frac{1}{\kappa}}(a)(G(t,a)-I)\nabla_a\tilde{W}^j(t,a)))\,da\\
        &=-2\kappa \int_{\Omega_0} \nabla_a\partial_t\tilde{\theta}^j\cdot r_0^{\frac{1}{\kappa}}(a)(G(t,a)-I)\nabla_a\tilde{W}^j(t,a))\,da.
    \end{align*}
    Thus, upon integrating by parts once in time we get that
    \begin{align*}
        \int_0^tI_d(s)\,ds&=-2\kappa \int_{\Omega_0} \nabla_a\tilde{\theta}^j\cdot r_0^{\frac{1}{\kappa}}(a)(G(t,a)-I)\nabla_a\tilde{W}^j(t,a))\,da\\
        &+2\kappa \int_0^t\int_{\Omega_0} \nabla_a\tilde{\theta}^j\cdot r_0^{\frac{1}{\kappa}}(a)\partial_t((G(t,a)-I)\nabla_a\tilde{W}^j(t,a)))\,da\,ds.
    \end{align*}
    It now immediately follows from our bounds on $G$ and on $\tilde{W}^j$ that
    \begin{align*}
        \left|\int_0^tI_d(s)\,ds\right|&\lesssim \|(\partial_t\tilde{\theta}^j,\nabla\tilde{\theta}^j)\|_{\mathcal{H}}\delta_j2^{j\left(d-\frac{1}{\kappa}\right)}
    \end{align*}
    We deduce that for every $T\geq 0$,
    \begin{align*}
        \|(\theta^j_t(T),\nabla \theta^j(T))\|^2_{\mathcal{H}}&\lesssim  \|(\theta^j_t(0),\nabla \theta^j(0))\|^2_{\mathcal{H}}+\left|\int_0^tI_a(s)\,ds\right|+\left|\int_0^tI_b(s)\,ds\right|+\left|\int_0^tI_c(s)\,ds\right|\\
        &+\left|\int_0^tI_d(s)\,ds\right|+\left|\int_0^tII(s)\,ds\right|\\
        &\lesssim CA^22^{2j\left(d-\frac{1}{\kappa}\right)}+C\delta_j\sup_{\substack{0\leq t\leq c2^{-j}}}\|(\partial_t\tilde{\theta}^j,\nabla\tilde{\theta}^j)\|^2_{\mathcal{H}}\\
        &+C(cA+\delta_j)2^{j\left(d-\frac{1}{\kappa}\right)}\sup_{\substack{0\leq t\leq c2^{-j}}}\|(\partial_t\tilde{\theta}^j,\nabla\tilde{\theta}^j)\|_{\mathcal{H}},
    \end{align*}
    where the bound for the initial data follows from Lemma \ref{Boundary layer approximation}.
   For $j$ large enough, we have $\displaystyle \int_0^{T_j}\|\nabla v(s)\|_{L_x^\infty}\,ds=\delta_j\ll 1$, and since $A\ll 1$, upon shrinking $c>0$ if necessary, we get that

    \begin{align*}
        \sup_{0\leq t\leq c2^{-j}}\|(\theta^j_t,\nabla\theta^j)\|_{\mathcal{H}}\lesssim (A+\delta_j)2^{j\left(d-\frac{1}{\kappa}\right)},
    \end{align*}
    as desired.
\end{proof}
We now return to the proof of our result. From our expression for $W^j$ we have that
\begin{align*}
    \partial^2_{Y_d}W^0(0,0,1)&=\frac{1}{(2\pi)^{d-1}}\int e^{-|\eta|}a(\eta)\,d\eta>0
\end{align*}
Upon choosing a smooth compactly supported function $\displaystyle\varphi\in C_c^\infty$ supported near $(0,1)$, and further shrinking $c>0$, we have
\begin{align*}
    \inf_{0\leq\tau\leq c}\int \partial^2_{Y_d}W^0(\tau,Y)\varphi(Y)\,dY>0
\end{align*}
We now set $\displaystyle\varphi^j(x)=\varphi(2^{2j}\Phi_t(x))$.
Let $\displaystyle Y=2^{2j}\Phi_t(x)$. We have
\begin{align*}
    \int_\Omega\partial^2_{x_d}&W_\ast^j(t,x)\varphi^j(x)\,dx=-\int_\Omega\partial_{x_d}W_\ast^j(t,x)\partial_{x_d}\varphi^j(x)\,dx\\
    &=-2^{2j}\int_\Omega (1+O(A+\delta_j))((e_d+O(A+\delta_j))\cdot\nabla_{Y}W^0(\tau,Y))((e_d+O(A+\delta_j))\cdot\nabla_{Y}\varphi(Y))\,dY\\
     &=-2^{2j}\int_\Omega \partial_{Y_d}W^0(\tau,Y)) \partial_{Y_d}\varphi(Y)\,dY+O(2^{2j}(A+\delta_j))\\
     &=2^{2j}\int_\Omega \partial^2_{Y_d}W^0(\tau,Y)) \varphi(Y)\,dY+O(2^{2j}(A+\delta_j)).
\end{align*}

Since $\displaystyle A\ll 1$, for $j$ sufficiently large it follows that
\begin{align*}
    \inf_{\substack{0\leq t\leq c2^{-j}}}\int_\Omega\partial^2_{x_d}\tilde{W}^j(t,x)\varphi^j(x)\,dx&=2^{2j(d-1)}\int\partial^2_{x_d}\tilde{W}^0(2^jt,2^{2j}x)\varphi^j(x)\,dx\gtrsim 2^{2j}.
\end{align*}
We also have
\begin{align*}
 \left|\int_\Omega\partial^2_{x_d}\theta^j(t,x)\varphi^j(x)\,dx\right|&=\left|\int_\Omega\partial_{x_d}\theta^j(t,x)\partial_{x_d}\varphi^j(x)\,dx\right| \\
 &\lesssim \|r^{\frac{1}{2\kappa}}\nabla\theta^j\|_{L^2}\|r^{-\frac{1}{2\kappa}}\partial_{x_d}\varphi^j(x)\|_{L^2}\\
 &\lesssim 2^{j\left(2+\frac{1}{\kappa}-d\right)}\|(\partial_t\theta^j,\nabla\theta^j)\|_{\mathcal{H}}\\
 &\lesssim (A+\delta_j)2^{j\left(d-\frac{1}{\kappa}\right)}2^{j\left(2+\frac{1}{\kappa}-d\right)}\lesssim (A+\delta_j)2^{2j}
\end{align*}
Let $\displaystyle V^j=W_\ast^j+\theta^j$. It is clear that $V^j$ solves
\begin{align*}
(D_t^2+\kappa (\nabla\cdot v)D_t&-\kappa r\Delta-\nabla r\cdot\nabla)V^j=0\\
V^j_{t=0}&=W^{j,\sharp}(0),D_tV^j_{t=0}=\partial_tW^{j,\sharp}(0),
\end{align*}
However, the previous estimates imply that
\begin{align*}
    \inf_{\substack{0\leq t\leq c2^{-j}}}\int_\Omega\partial^2_{x_d}V^j(t,x)\varphi^j(x)\,dx\gtrsim 2^{2j}. 
\end{align*}
Let $\displaystyle p\in[1,\infty]$. It is clear that $\displaystyle \|\varphi^j\|_{L^{p'}}\lesssim 2^{-\frac{2jd}{p'}}$. Thus,
\begin{align*}
     \inf_{\substack{0\leq t\leq c2^{-j}}}\|\nabla_x^2V^j\|_{L_x^{p}}&\gtrsim \|\partial_{x_d}^2V^j\|_{L_x^{p}}\gtrsim 2^{j\left(2d+2-\frac{2d}{p}\right)},
\end{align*}
which in turn implies that
\begin{align*}
     \|\nabla_x^2V^j\|_{L_t^qL_x^{p}([0,T]\times\Omega)}&\gtrsim \|\nabla_x^2V^j\|_{L_t^qL_x^{p}([0,c2^{-j}]\times\Omega)}\gtrsim (c2^{-j})^{\frac{1}{q}}\inf_{\substack{0\leq t\leq c2^{-j}}}\|\partial_{x_d}^2V^j\|_{L_x^{p}}\gtrsim 2^{j\left(2d+2-\frac{2d}{p}-\frac{1}{q}\right)},
\end{align*}

On the other hand, from Lemma \ref{Boundary layer approximation}, 
\begin{align*}
    \|(D_tV^j(0),\nabla V^j(0))\|_{\mathcal{H}^{2s}}&= \|(\partial_tW^{j,\sharp}(0),\nabla W^{j,\sharp}(0))\|_{\mathcal{H}^{2s}}\lesssim  2^{j\left(2s+d-\frac{1}{\kappa}\right)}.
\end{align*}

We now set
\begin{align*}
    U^j=\frac{V^j}{2^{j\left(2s+d-\frac{1}{\kappa}\right)}}.
\end{align*}
It then follows that
\begin{align*}
    \|(D_tU^j(0),\nabla U^j(0))\|_{\mathcal{H}^{2s}}\lesssim 1,
\end{align*}
while
   \begin{align*}
     \|\nabla_x^2U^j\|_{L_t^qL_x^{p}([0,T]\times\Omega)}&\gtrsim \|\nabla_x^2U^j\|_{L_t^qL_x^{p}([0,c2^{-j}]\times\Omega)}\gtrsim(c2^{-j})^{\frac{1}{q}}\inf_{\substack{0\leq t\leq c2^{-j}}}\|\partial_{x_d}^2U^j\|_{L_x^{p}}\\
     &\gtrsim 2^{j\left(d+2+\frac{1}{\kappa}-\frac{2d}{p}-\frac{1}{q}-2s\right)}\gg 1,
\end{align*} 
since $\displaystyle d+2+\frac{1}{\kappa}-\frac{2d}{p}-\frac{1}{q}-2s>0$. This finishes the proof.

\section{A physical solution}
In this short section we define an explicit affine spherical solution to compressible irrotational Euler equations which was first considered by Sideris in \cite{Sideris}. This is the solution that will provide the background for which we shall prove our desired Strichartz estimates in the desired open range and obtain results analogous to the ones in \cite{Avadanei}.

We consider the isotropic Euler background solution given by
\begin{align*}
    \Omega_t&=a(t)\mathbb{B}\\
    v(t,x)&=\frac{\dot{a}(t)}{a(t)}x\\
    r(t,x)&=ca(t)^{-d\kappa}\left(1-\frac{|x|^2}{a(t)^2}\right)\\
    \ddot{a}(t)&=2ca^{-d\kappa-1}.
\end{align*}
Here $\displaystyle c>0$, $\displaystyle a(0)=1$, and $\displaystyle a(t)>0$ for every $\displaystyle t\in[0,T]$, and $\mathbb{B}$ is the unit ball.

These describe irrotational solutions, and in three dimensions, they were studied by Sideris in \cite{Sideris} as part of an affine family of solutions. This construction is the isotropic specialization of the general affine family that he analyzed in \cite{Sideris2}.

We have 
\begin{align*}
 \frac{d}{dt}(a(t)y)= \frac{\dot{a}(t)}{a(t)}a(t)y=\frac{\dot{a}(t)}{a(t)}v(t,a(t)y) , 
\end{align*}
hence for fixed $y$, $t\rightarrow a(t)y$ describes the particle trajectory originating from $y$.

This motivates us to set $\displaystyle U(t,y):=\psi(t,a(t)y)$. 
Then, $\displaystyle D_t^2\psi(t,a(t)y)=\partial_t^2U(t,y)$, and $U$ solves the equation
\begin{align*}
    \partial_t^2 U- ca(t)^{-d\kappa-2}(\kappa R\Delta-2y\cdot\nabla)+d\kappa\frac{\dot{a}(t)}{a(t)}\partial_tU&=0,
\end{align*}
where $\displaystyle R(y)=1-|y|^2$. We introduce the change of variable $\displaystyle \tau=\sqrt{c}\int_0^ta(s)^{-\frac{d\kappa+2}{2}}\,ds$.

Now, if we set $\displaystyle V$ given by $\displaystyle U(t,y)=a(t)^{\frac{2-d\kappa}{4}}V(\tau(t),y)$, we obtain the equation
\begin{align*}
    \partial_\tau^2 V-\kappa R\Delta V+2y\cdot \nabla V+Q(\tau)V&=0,
\end{align*}
where
\begin{align*}
    Q(\tau)&=\frac{1}{(\tau'(t))^2}\left(\frac{2-d\kappa}{4}\frac{a''(t)a(t)-(a'(t))^2}{(a(t))^2}+\frac{(2-d\kappa)(2+3d\kappa)}{16}\frac{(a'(t))^2}{(a(t))^2}\right)\\
    &=\frac{(a(t))^{2+d\kappa}}{c}\left(\frac{2-d\kappa}{4}\frac{a''(t)}{a(t)}+\frac{(2-d\kappa)(-2+3d\kappa)}{16}\frac{(a'(t))^2}{(a(t))^2}\right)\\
    &=\frac{2-d\kappa}{2}+\frac{(2-d\kappa)(3d\kappa-2)}{16c}(a(t))^{d\kappa}(a'(t))^2.
\end{align*}
Here we have used the equation for $a$, which tells us that $\displaystyle a''(t)=2c(a(t))^{-d\kappa-1}$. In particular, it is clear that $\displaystyle Q(\tau)\in L_t^1([0,T])$.

We shall first prove the estimates in the case $Q(\tau)=0$. For example, this is the case if $d\kappa=2$.
The initial data becomes
\begin{align*}
    V(0)=\psi(0),V_\tau(0)=\frac{1}{\sqrt{c}}\left(\partial_t\psi(0)-\frac{2-d\kappa}{4}\dot{a}(0)\psi(0)\right)
\end{align*}
We also have $\displaystyle \nabla_x^2\psi(t,a(t)y)=a(t)^{-\frac{d\kappa+6}{4}}\nabla_y^2V(\tau(t),y)$. We have $\displaystyle dx=a(t)^d\,dy$ and $\displaystyle dt=(\lambda(t))^{-1}d\tau$, and since all of these quantities are bounded from above and below away from zero, proving the estimates for $\psi$ and $U$ reduces to proving that $\displaystyle \|\nabla_y^2V\|_{L_\tau^pL_y^q}\lesssim \|(V_\tau,\nabla_y V)\|_{\mathcal{H}_R^{2s}}$.
\section{A spectral analysis of the acoustic operators}

For $b>-1$, we define the operator $\displaystyle \mathcal{L}_b=-\kappa R\Delta+2\kappa(b+1)y\cdot\nabla$, so that our equation is
\begin{align*}
    \partial_\tau^2 V+\mathcal{L}_{\frac{1}{\kappa}-1}V&=0.
\end{align*}
Here we regard $\mathcal{L}_b$ as the nonnegative self-adjoint operator on $\displaystyle L^2(\mathbb{B},R^b\,dy)$ corresponding to the closure of the quadratic form
\begin{align*}
  q_b(f,g)&=\kappa\int_{\mathbb{B}}R^{b+1}\nabla f\cdot\overline{\nabla g}\,dy,\quad f,g\in C_c^\infty(\overline{\mathbb{B}}).  
\end{align*}
In our case, the closure is taken with respect tot the form $\displaystyle \|f\|_{L^2(R^b\,dy)}+\|\nabla f\|_{L^2(R^{b+1}\,dy)}$.

Our goal in this section is to study the spectral properties of these operators. In particular, we shall seek to diagonalize $\mathcal{L}_b$. We first note that the spectrum and the eigenfunctions of these operators were studied by Seis in \cite{Seis}. For the sake of completeness, we shall rederive them here.

Let $y=z\omega$, where $\displaystyle 0\leq z\leq 1$ and $\displaystyle \omega\in\mathbb{S}^{d-1}$. We consider the normalized spherical harmonic functions $Y_{\nu,l}$ satisfying the equations
\begin{align*}
    -\Delta_{\mathbb{S}^{d-1}}Y_{\nu,l}&=l(l+d-2)Y_{\nu,l}.
\end{align*}
We shall seek eigenfunctions for $\mathcal{L}_b$ of the form $\displaystyle F(y)=f(z)Y_{\nu,l}(\omega)$. We know that $\displaystyle \Delta_y=\partial_z^2+\frac{d-1}{z}\partial_z+\frac{1}{z^2}\Delta_{\mathbb{S}^{d-1}}$ and that $\displaystyle y\cdot\nabla_y=z\partial_z$.
We therefore have the relations
\begin{align*}
    \Delta_yF&=\left(f''(z)+\frac{d-1}{z}f'(z)-\frac{l(l+d-2)}{z^2}f(z)\right)Y_{l,\nu}(\omega)\\
    y\cdot \nabla_yF&=zf'(z)Y_{l,\nu}(\omega).
\end{align*}
These in turn imply that the eigenvalue equation $\displaystyle \mathcal{L}_bF=\lambda F$ is equivalent to
\begin{align*}
    -\kappa(1-z^2)\left(f''+\frac{d-1}{z}f'-\frac{l(l+d-2)}{z^2}f\right)+2\kappa(b+1)zf'=\lambda f
\end{align*}
We now set $\displaystyle f(z)=z^lp(u)$, where $\displaystyle u=z^2$. We note that $z^lY_{l,\nu}$ is harmonic. A direct computation shows that the equation $\displaystyle \mathcal{L}_bF=\lambda F$ reduces to
\begin{align*}
    z^lY_{l,\nu}(-4\kappa u(1-u)p''+(-2\kappa(2l+d)(1-u)+4\kappa(b+1)z)p'+2\kappa(b+1)lp)&=\lambda z^lp(u)Y_{l,\nu}.
\end{align*}
We set $\beta_l=l+\frac{d}{2}-1$, and the equation for $p$ becomes
\begin{align*}
    u(1-u)p''+(\beta_l+1-(b+\beta_l+2)u)p'+\frac{\lambda-2\kappa(b+1)l}{4\kappa}p&=0.
\end{align*}
Let $\displaystyle \mu=\frac{\lambda-2\kappa(b+1)l}{4\kappa}$ and $\displaystyle \mathcal{F}_p(u)=u^{\beta_l+1}(1-u)^{b+1}p'(u)$. We multiply our differential equation for $p$ by $\displaystyle u^{\beta_l}(1-u)^b$, and we get that $\displaystyle -\mathcal{F}_p'(u)=\mu u^{\beta_l}(1-u)^bp(u)$. Let now $F$ be an eigenfunction. Since $\displaystyle p(u)\in L^2((0,1),u^{\beta_l}(1-u)^b)$ and $\displaystyle u^{\beta_l}(1-u)^b\in L^1((0,1))$, it follows that $\displaystyle \lim_{\substack{u\rightarrow1_-}}\mathcal{F}_p(u)$ is well-defined. Let now $q$ be a smooth radial function which vanishes around the origin and equals $1$ near $1$.

Let $\displaystyle \Psi_q(u)=z^lq(z^2)Y_{l,\nu}(\omega)$. We recall that $\displaystyle \mathcal{L}_bF=-\kappa R^{-b}\text{div}(R^{b+1}\nabla)$, where $R$ is given by $R(y)=1-|y|^2$. Since $b+1>0$, the equation $\mathcal{L}_bF=\lambda F$, along with integration by parts and the fact that $\displaystyle R^{b+1}$ tends to zero towards $\partial\mathbb{B}$, implies that 
\begin{align*}
\int_{\mathbb{B}}R^{b}F\cdot\overline{\mathcal{L}_b\Psi_q}\,dy&=\kappa\int_{\mathbb{B}}R^{b+1}\nabla F\cdot\overline{\nabla\Psi_q}\,dy=\int_{\mathbb{B}}R^{b}\mathcal{L}_bF\cdot\overline{\Psi_q}\,dy+\kappa\int_{\partial\overline{B}}R^{b+1}(\partial_\nu F\overline{\Psi_q}-F\partial_\nu\overline{\Psi_q})\,dS\\
&=\int_{\mathbb{B}}R^{b}\mathcal{L}_bF\cdot\overline{\Psi_q}\,dy=\lambda\int_{\mathbb{B}}R^{b}F\cdot\overline{\Psi_q}\,dy.
\end{align*}
For $u=z^2$, we obtain as before that
\begin{align*}
    \mathcal{L}_b(z^lq(u)Y_{l,\nu})&=z^lY_{l,\nu}(-4\kappa u(1-u)q''+(-2\kappa(2l+d)(1-u)+4\kappa(b+1)u)q'+2\kappa(b+1)lq)\\
    &=z^lY_{l,\nu}\left(2\kappa(b+1)lq-4\kappa\frac{\frac{d}{du}(u^{\beta_l+1}(1-u)^{b+1}q'(u))}{u^{\beta_l}(1-u)^b}\right)
\end{align*}
Upon making the change of variable $u=z^2$, and using the fact that $\|Y_{l,\nu}\|_{L^2(\mathbb{S}^{d-1})}=1$, we get that
\begin{align*}
 \int_{\mathbb{B}}R^{b}F\cdot\overline{\mathcal{L}_b\Psi_q}\,dy&=\kappa(b+1)l\int_0^1u^{\beta_l}(1-u)^{b}p(u)\overline{q(u)}\,du-2\kappa\int_0^1 p(u)\frac{d}{du}(u^{\beta_l+1}(1-u)^{b+1}\overline{q'(u)})\,du \\
 \int_{\mathbb{B}}R^{b}F\cdot\overline{\Psi_q}\,dy&=\frac{1}{2}\int_0^1u^{\beta_l}(1-u)^{b}p(u)\overline{q(u)}\,du
\end{align*}
Now, the equality $\displaystyle \int_{\mathbb{B}}R^{b}F\cdot\overline{\mathcal{L}_b\Psi_q}\,dy=\lambda\int_{\mathbb{B}}R^{b}F\cdot\overline{\Psi_q}\,dy$ implies that
\begin{align*}
    \frac{\lambda}{2}\int_0^1u^{\beta_l}(1-u)^{b}p(u)\overline{q(u)}\,du&=\kappa(b+1)l\int_0^1u^{\beta_l}(1-u)^{b}p(u)\overline{q(u)}\,du\\
    &-2\kappa\int_0^1 p(u)\frac{d}{du}(u^{\beta_l+1}(1-u)^{b+1}\overline{q'(u)})\,du\\
    &=\kappa(b+1)l\int_0^1u^{\beta_l}(1-u)^{b}p(u)\overline{q(u)}\,du\\
    &+2\kappa\int_0^1 p'(u)u^{\beta_l+1}(1-u)^{b+1}\overline{q'(u)}\,du
\end{align*}
This implies that
\begin{align*}
    \int_0^1u^{\beta_l+1}(1-u)^{b+1} p'(u)\overline{q'(u)}\,du&=\frac{\lambda-2\kappa(b+1)l}{4\kappa}\int_0^1u^{\beta_l}(1-u)^{b}p(u)\overline{q(u)}\,du,
\end{align*}
hence
\begin{align*}
    \int_0^1\mathcal{F}_p(u)\overline{q'(u)}\,du&=\mu\int_0^1 u^{\beta_l}(1-u)^bp(u)\overline{q(u)}\,du.
\end{align*}

On the other hand,
\begin{align*}
    \int_0^1\mathcal{F}_p(u)\overline{q'(u)}\,du&=\lim_{\substack{u\rightarrow1_-}}\mathcal{F}_p(u)-\int_0^1\mathcal{F}_p'(u)\overline{q(u)}\,du\\
    &=\lim_{\substack{u\rightarrow1_-}}\mathcal{F}_p(u)+\mu\int_0^1 u^{\beta_l}(1-u)^bp(u)\overline{q(u)}\,du
\end{align*}

We therefore need $\displaystyle \lim_{\substack{u\rightarrow1_-}}\mathcal{F}_p(u)=0$. Around the origin, our operator is uniformly elliptic, so the solution of
\begin{align*}
    u(1-u)p''+(\beta_l+1-(b+\beta_l+2)u)p'+\frac{\lambda-2\kappa(b+1)l}{4\kappa}p&=0
\end{align*}
must be also be regular around the origin. It follows from \cite{DLMF}, $15.10.1$ and $15.10.2$, that $\displaystyle p(u)=_2F(a_+,a_-;\beta_l+1;u)$, where $a_+$ and $a_-$ satisfy $\displaystyle a_++a_-=b+\beta_l+1$, and $\displaystyle a_+a_-=-\mu$. From formula $15.5.1$ in \cite{DLMF}, it follows that
\begin{align*}
    \mathcal{F}_p(u)&=u^{\beta_l+1}(1-u)^{b+1}\frac{a_+a_-}{c_l}_2F(a_++1,a_-+1;\beta_l+2;u)\\
    &=u^{\beta_l+1}\frac{a_+a_-}{c_l}_2F(\beta_l+1-a_+,\beta_l+1-a_-+1;\beta_l+2;u)
\end{align*}
Since the excess of the lower parameter over the sum of the upper parameter is $\displaystyle \beta_l+2-(\beta_l+1-a_+)-(\beta_l+1-a_-)=b+1>0$, by the Gauss formula (see $15.4.20$ in \cite{DLMF}) we have
\begin{align*}
    \lim_{\substack{u\rightarrow 1_-}}\mathcal{F}_p(u)&=\frac{\Gamma(\beta_l+1)\Gamma(b+1)}{\Gamma(a_+)\Gamma(a_-)}. 
\end{align*}
Since $\beta_l+1$ and $b+1$ are both nonzero and finite, it follows that $a_+$ or $a_-$ must be a zero of $\displaystyle\frac{1}{\Gamma}$, hence $a_+=-n$ or $a_-=-n$ for some nonnegative integer $n$. If for instance $a_+=-n$, we have $a_-=n+b+\beta_l+1$ and $\mu=n(n+b+\beta_l+1)$, which in turn give
\begin{align*}
 \lambda^{(b)}_{n,l}=2\kappa(b+1)l+4\kappa n\left(n+l+b+\frac{d}{2}\right)   
\end{align*}
and the eigenfunctions 
\begin{align*}F^{(b)}_{n,l,\nu}(y)=z^lP^{(b,\beta_l)}_n(2z^2-1)Y_{l,\nu}(\omega),
\end{align*}
where $P^{(b,\beta_l)}_n$ is a Jacobi polynomial. We now normalize our resulting eigenfunctions. We can of course check that these solutions have indeed flux zero at $1$.

For the Jacobi polynomials, we have
\begin{align*}
  \int_0^1 z^{2l+d-1}(1-z^2)^b|P^{(b,\beta_l)}(2z^2-1)|^2\,dz&= \frac{\Gamma(n+b+1)\Gamma\left(n+l+\frac{d}{2}\right)}{2n!\left(2n+b+l+\frac{d}{2}\right)\Gamma\left(n+b+l+\frac{d}{2}\right)} 
\end{align*}
We denote this common value by $h^{(b)}_{n,l}$, and we can then define $\displaystyle \phi^{(b)}_{n,l}(z)=\frac{z^lP_n^{(b,\beta_l)}(2z^2-1)}{\sqrt{h^{(b)}_{n,l}}}$, so that our normalized eigenfunctions are given by
\begin{align*}
    \Phi^{(b)}_{n,l,\nu}(y)=\phi^{(b)}_{n,l}(z)Y_{l,\nu}(\omega).
\end{align*}
They form an orthonormal basis in $\displaystyle L^2(R^b\,dy)$ due to their completeness, which in turn is a consequence of the completeness of the spherical harmonics, and the the one of the Jacobi polynomials in $L^2((0,1),u^{\beta_l}(1-u)^{b}\,du)$ in each angular sector.

Upon rescaling, we can also check that the eigenvalues and eigenvalues that we have obtained coincide with the ones from \cite{Seis}, see equation $(14)$ and Theorem 1.
\section{The proof of the Strichartz estimates}
In order to prove our third main result, which provides Strichartz estimates for our specific isotropic Euler background, we are going to decompose our solution into dyadic blocks depending on the dyadic normal and angular frequencies, and then carry out a $TT^{\ast}$ argument for each resulting dyadic block.

We use the completeness of our set of eigenfunctions, and we write $\displaystyle \Phi^{\left(\frac{1}{\kappa}-1\right)}_{n,l,\nu}=\Phi_j$, where $\displaystyle j=(n,l,\nu)$. We then decompose 
\begin{align*}
    V(0)&=\sum_{\substack{j}}f_j\Phi_j\\
    \partial_\tau V(0)&=\sum_{\substack{j}}g_j\Phi_j.
\end{align*}
We note that the only time zero occurs as an eigenvalue is as $\lambda_{(0,0)}$. Since this corresponds to a mode that is constant in space which commutes with $Q(\tau)$ and $\mathcal{L}_{\frac{1}{\kappa}-1}$ and vanishes under the Hessian, we can project it away and only retain the nonzero modes. In what follows, we shall also denote this contribution by $V$. In particular, we can also see that the smallest nonzero eigenvalue is $2$.

We immediately get that
\begin{align*}
    V(\tau)&=\sum_{\substack{j=(n,\nu)}}\left(\cos(\tau\sqrt{\lambda_j})f_j+\frac{\sin(\tau\sqrt{\lambda_j})}{\sqrt{\lambda_j}}g_j\right)\Phi_j.
\end{align*}
We also similarly decompose $\displaystyle V_{N,L}=V^+_{N,L}+V^-_{N,L}$.

For the base acoustic metric we then obtain
\begin{align*}
    \|(V_\tau,\nabla_yV)\|^2_{\mathcal{H}}&=\sum_{\substack{j}}(|g_j|^2+\lambda_j|f_j|^2)
\end{align*}
Let $\displaystyle s_{q,r}=\frac{d+b+1}{2}-\frac{d}{r}-\frac{1}{2q}$ and $\displaystyle \delta_{q,r}=\frac{d-1}{2}\left(1-\frac{2}{r}\right)-\frac{2}{q}>0$. We would like to prove that
\begin{align*}
    \|\nabla_y^2V\|_{L_\tau^q L_y^r(I\times\mathbb{B})}&\lesssim \|(V_\tau,\nabla_y V)(0)\|_{\mathcal{H}_R^{2s}}\,\text{, for } 2s>2s_{q,r}.
\end{align*}
We shall work on dyadic blocks on which $n+1\sim N$ and $l+1\sim L$, where $N$ and $L$ are both dyadic numbers. 

We set $\displaystyle V_{N,L}=\sum_{\substack{j\\
n+1\sim N,l+1\sim L}}V_j$, where $\displaystyle j=(n,l,\nu)$. For the sake of simplicity, we also introduce the notations $\displaystyle W_{N,L}=N(N+L)$ and $\displaystyle \Lambda_{N,L}=W_{N,L}^{\frac{1}{2}}$. When $n+1\sim N$ and $l+1\sim L$, since $\displaystyle \lambda_{n,l}=2l+4\kappa n\left(n+l-1+\frac{1}{\kappa}+\frac{d}{2}\right)$, it follows that $\displaystyle \lambda_{n,l}\simeq W_{N,L}$. 

We also introduce the notation
\begin{align*}
    F^{\pm}_j&=\frac{1}{2}(f_j\mp i\lambda_{n,l}^{-\frac{1}{2}}g_j).
\end{align*}
It is then clear that $\displaystyle V=V^++V^-$, where $\displaystyle V^\pm(\tau)=\sum_{\substack{j=(n,\nu)}} e^{\pm i\tau\sqrt{\lambda_{n,l}}}F^\pm_j\Phi_j$.
We also denote by $\displaystyle\mathcal{E}_{N,L}$ the energy of this given dyadic block. Namely, it is given by the relation
\begin{align*}
    \mathcal{E}_{N,L}&=\|\partial_\tau V_{N,L}(0)\|_{L_y^2(R^{\frac{1}{\kappa}-1})}+\kappa \|\nabla_y V_{N,L}(0)\|_{L_y^2(R^{\frac{1}{\kappa}})}. 
\end{align*}
We shall aim to prove that
\begin{align*}
    \|\nabla_y^2V_{N,L}\|_{L_\tau^q L_y^r(I\times\mathbb{B})}&\lesssim W_{N,L}^{s_{q,r}}N^{-\delta_{q,r}}\mathcal{E}_{N,L},
\end{align*}
where $\displaystyle s_{q,r}=\frac{d+b+1}{2}-\frac{d}{r}-\frac{1}{2q}$ and $\displaystyle \delta_{q,r}=\frac{d-1}{2}\left(1-\frac{2}{r}\right)-\frac{2}{q}$.

It suffices to prove this estimate for an arbitrary pair of partial derivatives. In what follows, we fix such a pair and we set $\displaystyle T_{ij}=\partial_i\partial_j$. We are going to prove that
\begin{align*}
    \|T_{ij}V^+_{N,L}\|_{L_\tau^q L_y^r(I\times\mathbb{B})}&\lesssim W_{N,L}^{s_{q,r}}N^{-\delta_{q,r}}\mathcal{E}_{N,L}.
\end{align*}
The estimate involving $\displaystyle V^-_{N,L}$ is similar. For the sake of simplicity, we shall use $V_{N,L}$ instead of $V^+_{N,L}$, suppressing the plus superscript.

Our eigenfunctions have some useful properties which are the object of the following
\begin{lemma}\label{Derivative Expansion}
    Let $\displaystyle d\geq 2$, $\displaystyle \kappa>0$, and $\displaystyle \gamma>-1$. We set $\displaystyle R(y)=1-|y|^2$ and $\displaystyle H_\gamma=L^2(\mathbb{B},R^\gamma\,dy)$. We denote $\displaystyle \mathcal{L}_\gamma=-\kappa R\Delta_y+2\kappa(\gamma+1)y\cdot\nabla_y$. Let $\displaystyle y=z\omega$, where $\displaystyle z=|y|$, and $\displaystyle\omega\in\mathbb{S}^{d-1}$. We recall the form of the normalized eigenfunctions
    \begin{align*}
        \Phi^{(b)}_{n,l,\nu}(y)&=\frac{z^lP^{(b,\beta_l)}(2z^2-1)}{\sqrt{h^{(b)}_{n,l}}}Y_{l,\nu}(\omega),
    \end{align*}
    where $\displaystyle \beta_l=l+\frac{d}{2}-1$, and associated eigenvalues $\displaystyle \lambda_{n,l}=4\kappa n\left(n+l+\gamma+\frac{d}{2}\right)+2\kappa(\gamma+1)l$. Let $\displaystyle \mathcal{Y}^{(\gamma)}_{n,l}$ be the space generated by the eigenfunctions $\displaystyle \Phi^{(b)}_{n,l,\nu}$. We then have the following
    \begin{itemize}
        \item $\displaystyle\partial_i\mathcal{Y}^{(\gamma)}_{n,l}\subset \mathcal{Y}^{(\gamma+1)}_{n,l-1}\oplus \mathcal{Y}^{(\gamma+1)}_{n-1,l+1}$.
        \item $\displaystyle\partial_i\partial_j\mathcal{Y}^{(\gamma)}_{n,l}\subset \mathcal{Y}^{(\gamma+2)}_{n,l-2}\oplus \mathcal{Y}^{(\gamma+2)}_{n-1,l}\oplus \mathcal{Y}^{(\gamma+2)}_{n-2,l+2}$
        \item For an expansion of the form $\displaystyle F_{N,L}=\sum_{\substack{n+1\sim N\\l+1\sim L}}\sum_{\substack{\nu}}c_{n,l,\nu}\Phi^{(b)}_{n,l,\nu}$, we have the estimates
        \begin{align*}
            \|\nabla F_{N,L}\|_{H_{\gamma+1}}&\lesssim \Lambda_{N,L}\|F_{N,L}\|_{H_{\gamma}}\\
            \|\nabla^2 F_{N,L}\|_{H_{\gamma+2}}&\lesssim \Lambda^2_{N,L}\|F_{N,L}\|_{H_{\gamma}}
        \end{align*}

        Moreover, if we write $\displaystyle \partial_i\partial_j\Phi^{(\gamma)}_{n,l,\nu}=\sum_{\substack{a=0}}^2\sum_{\substack{\mu}}d^{(a,\mu)}_{ij;n,l,\nu}\Phi^{(\gamma+2)}_{n-a,l-2+2a,\mu}$ and then define \newline $\displaystyle G^{(a)}_{n-a,l-2+2a,\mu}=\sum_{\substack{\nu}}d^{(a,\mu)}_{ij;n,l,\nu}c_{n,l,\nu}$, we have
        \begin{align*}
            \sum_{\substack{a=0}}^2\|G^{(a)}\|^2_{l^2}&\lesssim \Lambda^4_{N,L}\|c\|^2_{l^2}.
        \end{align*}
    \end{itemize}
\end{lemma}
\begin{proof}
   Let $\mathcal{P}^{(\gamma)}_k$ be the space of polynomials of degree at most $k$ which are orthogonal to all polynomials of degree at most $k-1$ in $H_\gamma$. It is clear that $\displaystyle \Phi^{(\gamma)}_{n,l,\nu}\in\mathcal{P}^{(\gamma)}_{l+2n}$.

   Let now $F\in\mathcal{P}^{(\gamma)}_k$ be arbitrary and $G$ a polynomial of degree at most $k-2$. We then have
   \begin{align*}
       \int_{\mathbb{B}} \partial_iF\overline{G}R^{\gamma+1}\,dy&=-\int_{\mathbb{B}}F\partial_i(\overline{G}R^{\gamma+1})\,dy=-\int_{\mathbb{B}}F\partial_i\overline{G}R^{\gamma+1}\,dy+2(\gamma+1)\int_{\mathbb{B}}F\overline{G}R^{\gamma}y_i\,dy\\
       &=0.
   \end{align*}
   Here we have used the fact that after integrating by parts, $F$ ended up being integrated against a polynomial of degree at most $k-1$.

   Therefore, for every index $i$, $\displaystyle \partial_i\mathcal{P}^{(\gamma)}_k\subset\mathcal{P}^{(\gamma+1)}_{k-1}$, which encapsulates the weight shift. We now analyze the change in angular degree. For an arbitrary eigenfunction $F$ we write
   \begin{align*}
   F(y)=H_l(y)p(z^2),
   \end{align*}
   where $z=|y|$ and $\displaystyle H_l(y)=z^lY_{l,\nu}(\omega)$ is a homogeneous harmonic polynomial of degree $l$. We have $\displaystyle \partial_iF=\partial_i(z^lY_{l,\nu}(\omega))p(z^2)+2H_l(y)y_ip'(z^2)$. We can immediately see that the first term has angular degree $l-1$. As for the second one, we note that $\displaystyle \Delta(y_iH_l)=2\partial_iH_l+y_i\Delta H_l$ and that $\displaystyle \Delta(|y|^2\partial_iH_l)=2(2l+d-2)\partial_iH_l$, hence
   \begin{align*}
       y_iH_l=H_{l+1,i}+\frac{|y|^2}{2l+d-2}\partial_iH_l,
   \end{align*}
   where $H_{l+1,i}$ is a homogeneous harmonic polynomial of degree $l+1$. Thus, only the angular degrees $l-1$ and $l+1$ occur. When $l=0$, we can immediately see that $y_iH_0$ is already harmonic. In all instances, the total polynomial degree is $l+2n-1$, hence the conclusion follows. By iterating, we obtain the conclusion for second order derivatives.

   We now move on to proving the derivative estimates. Whenever $F$ is a polynomial, we have $\displaystyle \langle F,\mathcal{L}_\gamma F\rangle_{H_\gamma}=\kappa \|\nabla F\|^2_{H_{\gamma+1}}$. For every polynomial which appears in the decomposition of $F$, its associated eigenvalue $\lambda_{n,l}$ satisfies $\displaystyle 1+\lambda_{n,l}\lesssim W_{N,L}$. This immediately implies that $\displaystyle \|\nabla F\|^2_{H_{\gamma+1}}\lesssim W_{N,L}\| F\|^2_{H_{\gamma}}$. Repeating the argument also yields the second order estimate $\displaystyle \|\nabla^2 F\|^2_{H_{\gamma+2}}\lesssim W_{N,L}^2\| F\|^2_{H_{\gamma}}$.

   Finally, in a given $(n,l)$ sector, the orthogonality in the three branches described by$G^{(a)}_{n-a,l-2+2a,\mu}$ gives the estimate 
   \begin{align*}
       \sum_{\substack{a=0}}^2\|G^{(a)}\|^2_{l^2}\lesssim W_{N,L}^2\|c\|^2_{l^2},
   \end{align*}
   as desired.
\end{proof}
By taking $\displaystyle \gamma=\alpha$ and $\displaystyle b=\alpha+2$,  the estimate $\displaystyle \|F^{\pm}_{N,L}\|_{l^2}\lesssim \Lambda_{N,L}^{-1}\mathcal{E}_{N,L}$ now implies
\begin{align*}
    \sum_{\substack{a=0}}^2\|G^{(a)}\|_{l^2}&\lesssim \Lambda_{N,L}^2\|F^{\pm}_{N,L}\|_{l^2}\lesssim \Lambda_{N,L}\mathcal{E}_{N,L}.
\end{align*}

We can now write 
\begin{align*}
    T_{ij}V^+_{N,L}&=\sum_{\substack{a=0}}^2\sum_{\substack{n,l,\nu}}e^{i\tau\varpi_n^{(a)}(l)}G^{(a),+}_{n,l,\nu}\Phi^{(b)}_{n,l,\nu},
\end{align*}
where $\displaystyle \varpi_n^{(a)}(l)=\sqrt{\lambda_{n+a,l+2-2a}}$, $\displaystyle a=\overline{0,2}$.

We shall first prove our desired Strichartz estimate around the boundary of the ball, where the elliptic operators $\mathcal{L}_\gamma$ are degenerate. We choose a smooth radial cut-off function $\displaystyle \chi_{\rm bd}$ which is identically equal to one on $\displaystyle \left\{z||z|\geq \frac{2}{3}\right\}$ and is supported in $\displaystyle \left\{z||z|\geq \frac{1}{2}\right\}$.

For $\displaystyle G\in l^2$ indexed by $(n,l,\nu)$, we define the operator $\mathcal{V}_{N,L}$ by
\begin{align*}
    \mathcal{V}_{N,L}(\tau)G(y)=\chi_{\rm bd}(z)\sum_{\substack{n+1\sim N\\l,\nu}}e^{i\tau\varpi_n^{(a)}(l)}\chi(l/L)\phi^{(b)}_{n,l}(z)Y_{l,\nu}(\omega)G_{n,l,\nu},
\end{align*}
where $\displaystyle z=|y|$ and $\chi$ is a smooth compactly supported angular cut-off function which is identically equal to $1$ on the support of the coefficients.

We can see that the kernel of the operator $\displaystyle \mathcal{V}_{N,L}(\tau)\mathcal{V}_{N,L}(\tau')^\ast$ is given by
\begin{align*}
    K_{N,L}(u;y,\overline{y})&=\chi_{\rm bd}(z)\chi_{\rm bd}(\overline{z})\sum_{\substack{n+1\sim N\\l,\nu}}e^{iu\varpi_n^{(a)}(l)}\chi^2(l/L)\phi^{(b)}_{n,l}(z)\phi^{(b)}_{n,l}(\overline{z})Y_{l,\nu}(\omega)\overline{Y_{l,\nu}(\overline{\omega})}.
\end{align*}
Here, $\displaystyle \overline{z}=|\overline{y}|$, $\overline{\omega}\in\mathbb{S}^{d-1}$, and $u=\tau-\tau'$. 

We shall first consider the regime $\displaystyle N\leq L^{\frac{1}{3}}$. It turns out that the radial components $\phi^{(b)}_{n,l}(z)$ of our eigenfunctions have some useful properties which are the object of the following
\begin{lemma}\label{Uniform Lemma}
    Whenever $b\geq 0$ we have the bounds
    \begin{align*}
        |\phi^{(b)}_{n,l}(z)|^2&\lesssim (n+1)^b(n+l+1)^{b+1}\,,z\geq\frac{1}{2}.\\
        \sup_{\substack{z\geq\frac{1}{2}}}\sum_{\substack{n\geq 0\\
        1+\lambda^{(b)}_{n,l}\leq C\Lambda^2}}|\phi^{(b)}_{n,l}(z)|^2&\lesssim\Lambda^{2b+2}.
    \end{align*}
\end{lemma}
\begin{proof}
    From \cite[18.14.3-18.14.5]{DLMF}, we have the Jacobi inequality
    \begin{align*}
        \left(\frac{1+x}{2}\right)^{\frac{\beta}{2}}|P_n^{(b,\beta)}(x)|\leq P^{(b,\beta)}(1)=\frac{(b+1)_n}{n!}\,, b,\beta\geq 0\,, x\in[-1,1].
    \end{align*}
    Thus, for $\displaystyle z\geq\frac{1}{2}$, we have
    \begin{align*}
        z^{\beta_l}|P_n^{(b,\beta)}(2z^2-1)|\lesssim_b(n+1)^b.
    \end{align*}
    Since $\displaystyle \beta_l=l-1+\frac{d}{2}$, we get that
    \begin{align*}
        z^{l}|P_n^{(b,\beta)}(2z^2-1)|\leq 2^{\frac{d}{2}-1}\frac{(b+1)_n}{n!}\simeq_{b,d}(n+1)^b.
    \end{align*}
    The usual asymptotics for the $\Gamma$ function further imply that $\displaystyle h^{(b)}_{n,l}\simeq (n+1)^b(n+l+1)^{-b-1}$. Thus,
    \begin{align*}
        |\phi^{(b)}_{n,l}(z)|&\lesssim (n+1)^{\frac{b}{2}}(n+l+1)^{\frac{b+1}{2}}\,,z\geq\frac{1}{2}.
    \end{align*}
    Thus, for every $\displaystyle z\geq\frac{1}{2}$,
    \begin{align*}
        \sum_{\substack{n\geq 0\\
        1+\lambda^{(b)}_{n,l}\leq C\Lambda^2}}|\phi^{(b)}_{n,l}(z)|^2&\lesssim \sum_{\substack{n\geq 0\\
        1+\lambda^{(b)}_{n,l}\leq C\Lambda^2}}(n+1)^{b}(n+l+1)^{b+1}\lesssim \Lambda^{2b+2}
    \end{align*}
\end{proof}
Moreover, when we apply derivatives with respect to the angular parameter $l$, we have the following
\begin{lemma}\label{Uniform Derivative bounds lemma}
    When $\displaystyle N\leq L^{\frac{1}{2}}$, we have
    \begin{align*}
      \sup_{\substack{z\geq\frac{1}{2}}}|(L\partial_l)^j\phi^{(b)}_{n,l}(z)|&\lesssim_j L^{\frac{b+1}{2}}N^{\frac{b}{2}+j}. 
    \end{align*}
    This in particular implies the estimate
    \begin{align*}
      \sup_{\substack{z,\overline{z}\geq\frac{1}{2}}}|(L\partial_l)^j(\phi^{(b)}_{n,l}(z)\phi^{(b)}_{n,l}(\overline{z}))|&\lesssim_j L^{b+1}N^{b+j}. 
    \end{align*}
\end{lemma}
\begin{proof}
    Let $\displaystyle S=1-z^2$. It is clear that $\displaystyle 0\leq S\leq \frac{3}{4}$. For the sake of simplicity, we also set $\displaystyle B=b+n+l+\frac{d}{2}$.
    From \cite[18.5.7,18.5.12]{DLMF}, we deduce that the expansions for the Jacobi and Laguerre polynomials imply that
    \begin{align*}
        P_n^{(b,\beta_l)}&=\frac{(b+1)_n}{n!}\sum_{\substack{k=0}}^n\frac{(-n)_k(B)_k}{(b+1)_kk!}S^k\\
        L_n^{(b)}(Su)&=\frac{(b+1)_n}{n!}\sum_{\substack{k=0}}^n\frac{(-n)_k}{(b+1)_kk!}(Su)^k
    \end{align*}
    We define the weighted averaged operator $\mathcal{T}_B$ by
    \begin{align*}
        \mathcal{T}_Bf=\frac{1}{\Gamma(B)}\int_0^\infty u^{B-1}e^{-u}f(u)\,du.
    \end{align*}
    We can also immediately check that for every $k\geq 0$, $\displaystyle \mathcal{T}_Bu^k=\frac{\Gamma(B+k)}{\Gamma(B)}=(B)_k$. This implies that
    \begin{align*}
      P_n^{(b,\beta_l)}(1-2S)&=\mathcal{T}_B[L_n^{(b)}(S\cdot)]  .
    \end{align*}

    We also introduce the weighted Laguerre profiles $\displaystyle \mathcal{B}_n(X)=e^{-X}L_n^{(b)}(2X)$. In particular, $\displaystyle L_n^{(b)}(Su)=e^{\frac{Su}{2}}\mathcal{B}_n\left(\frac{Su}{2}\right)$. Thus, we have the following
    \begin{align*}
      P_n^{(b,\beta_l)}(1-2S)&=\frac{1}{\Gamma(B)}\int_0^\infty u^{B-1}e^{-u\left(1-\frac{S}{2}\right)}\mathcal{B}_n\left(\frac{Su}{2}\right)\,du\\
      &=\frac{\left(1-\frac{S}{2}\right)^{-B}}{\Gamma(B)}\int_0^\infty v^{B-1}e^{-v}\mathcal{B}_n\left(\frac{Sv}{2-S}\right)\,dv\\
      &=\left(1-\frac{S}{2}\right)^{-B}\mathcal{T}_B\left[\mathcal{B}_n\left(\frac{S}{2-S}\cdot\right)\right]. 
    \end{align*}
    Our eigenfunctions then satisfy
    \begin{align*}
        \phi^{(b)}_{n,l}(z)&=(h^{(b)}_{n,l})^{-\frac{1}{2}}C_{n,l}(S)\mathcal{T}_Bf_S,
    \end{align*}
    where
    \begin{align*}
     C_{n,l}(S)&=z^l\left(1-\frac{S}{2}\right)^{-B}\\ 
     f_S(u)&=\mathcal{B}_n\left(\frac{S}{2-S}u\right).
    \end{align*}
    We first analyze the case when derivatives fall onto the Laguerre factors. For this purpose, we shall derive a few bounds for their derivatives. From \cite[18.14.8]{DLMF}, we have
    \begin{align*}
       e^{-X}|L_n^{(b)}(2X)|\leq L_n^{(b)}(0)=\frac{(b+1)_n}{n!}\lesssim_b(n+1)^b.  
    \end{align*}
    This implies that $\displaystyle\|\mathcal{B}_n\|_{L^\infty}\lesssim_b(n+1)^b$. For the sake of simplicity, we also introduce the operator $\displaystyle \mathcal{D}_X=X\partial_X$. From \cite[18.9.13, 18.9.14, 18.9.23]{DLMF}, we have the recurrence relation
    \begin{align*}
        \mathcal{D}_X\mathcal{B}_n&=\frac{n+1}{2}\mathcal{B}_{n+1}-\frac{b+1}{2}\mathcal{B}_{n}-\frac{n+b}{2}\mathcal{B}_{n-1},
    \end{align*}
    where $n\geq 1$. This in particular implies that for every $\displaystyle j\geq 1$,
    \begin{align*}
        \mathcal{D}_X^j\mathcal{B}_n&=\sum_{\substack{k\geq 0\\|k-n|\leq j}}a^{(j)}_{n,k}\mathcal{B}_k,
    \end{align*}
    where $\displaystyle |a^{(j)}_{n,k}|\lesssim_{b,j}(n+1)^j$. This immediately implies that $\displaystyle \|\mathcal{D}_X^j\mathcal{B}_n\|_{L^\infty}\lesssim_{j,b}(n+1)^{j+b}\lesssim_{j,b}N^{j+b}$.
    
    On the other hand, we have the identity $\displaystyle X^j\partial^j_X=\Pi_{\substack{k=0}}^{j-1}(\mathcal{D}_X-k)$, hence 
    \begin{align*}
        \|X^j\partial_X^j\mathcal{B}_n\|_{L^\infty}&\lesssim_{j,b}N^{j+b}
    \end{align*}
    We immediately get that
    \begin{align*}
        |f_S^{(j)}(u)|&\lesssim N^{j+b}u^{-j}.
    \end{align*}
    We now introduce the operator $\displaystyle \mathcal{A}$ defined by
    \begin{align*}
        \mathcal{A}f(u)&=\int_0^\infty e^{-v}\frac{f(u+v)-f(u)}{v}\,dv
    \end{align*}
    We claim that whenever $f$ is a product between a decaying or constant exponential and a polynomial, $\displaystyle \partial_B^j\mathcal{T}_Bf=\mathcal{T}_B[\mathcal{A}^jf]$.  For $f(u)=e^{-cu}$, where $c\geq 0$, we have
    \begin{align*}
       \mathcal{T}_Bf&=\frac{1}{\Gamma(B)}\int_0^\infty u^{B-1}e^{-(c+1)u}\,du=(c+1)^{-B}\\
       \mathcal{A}f(u)&=e^{-cu}\int_0^\infty \frac{e^{-(c+1)v}-e^{-v}}{v}\,dv=-e^{-cu}\log(1+c).
    \end{align*}
    This now immediately implies that
    \begin{align*}
       \mathcal{T}_B[\mathcal{A}^jf]&=(-\log(1+c))^j(1+c)^{-B}=\partial_B^j \mathcal{T}_B[f]
    \end{align*}
    We also note that when $k\geq 1$, $\displaystyle u^ke^{-cu}=(-\partial_c)^ke^{-cu}$. This immediately implies the desired property for the case when $f$ is a product between a nonincreasing exponential and a monomial $u^k$, which immediately allows us to infer our conclusion.

    We can now prove that $\displaystyle|\partial_B^k\mathcal{T}_Bf_S|\lesssim B^{-k}N^{b+k}$. We shall first prove it for $k=1$.

    We have $\displaystyle f_S(u+v)-f_S(u)=\int_0^vf_S'(u+t)\,dt$. In this case, for every $u\geq 0$,
    \begin{align*}
        |\mathcal{A}f_S(u)|&\leq\int_0^\infty e^{-v}\frac{|f(u+v)-f(u)|}{v}\,dv\leq \int_0^\infty e^{-v}\frac{vN^{1+b}u^{-1}}{v}\,dv\lesssim N^{1+b}u^{-1}.
    \end{align*}
    In this case,
    \begin{align*}
        |\partial_B\mathcal{T}_Bf_S|&\lesssim |\mathcal{T}_B[\mathcal{A}f_S]|\lesssim \frac{1}{\Gamma(B)}\int_0^\infty e^{-u}u^{B-2}N^{1+b}\,du\\
        &=N^{1+b}\frac{\Gamma(B-1)}{\Gamma(B)}=\frac{N^{1+b}}{B-1}
    \end{align*}
    It immediately follows that $\displaystyle |\partial_B\mathcal{T}_Bf_S|\lesssim N^{1+b}B^{-1}\lesssim N^{1+b}L^{-1}$ for $\displaystyle B\geq 2$. The bound in the case $B=1$ can also be obtained immediately.

    By iterating this estimate, we get that \begin{align*} |\mathcal{A}^jf_s(u)|\lesssim_{j,b}N^{j+b}u^{-j}.
    \end{align*}
    In the case $B\geq 2j+1$, a direct integration then allows us to deduce the bound 
\begin{align*}
|\partial_B^j\mathcal{T}_Bf_S|\lesssim_{j,b}N^{j+b}\frac{\Gamma(B-j)}{\Gamma(B)}\lesssim_{j,b}N^{j+b}B^{-j}.
    \end{align*}
    The remaining cases, which consist of bounded frequencies, can be estimated directly.

    We now analyze the case in which the derivative falls onto 
    \begin{align*}
    (h^{(b)}_{n,l})^{-\frac{1}{2}}&=\left(2\left(2n+b+l+\frac{d}{2}\right)\frac{\Gamma(n+1)}{\Gamma(n+b+1)}\frac{\Gamma\left(n+b+l+\frac{d}{2}\right)}{\Gamma\left(n+l+\frac{d}{2}\right)}\right)^{\frac{1}{2}}.
    \end{align*}
    In our regime, $\displaystyle n+l+\frac{d}{2}\simeq L$, hence $\displaystyle (h^{(b)}_{n,l})^{-\frac{1}{2}}\simeq L^{\frac{b+1}{2}}N^{-\frac{b}{2}}$. Let $\displaystyle \psi(x)=\frac{\Gamma'(x)}{\Gamma(x)}$ be the digamma function. We can see that
    \begin{align*}
        \frac{d^l}{dl^r}\log((h^{(b)}_{n,l})^{-\frac{1}{2}})&=\frac{1}{2}\left[\frac{(-1)^{r-1}(r-1)!}{\left(2n+b+l+\frac{d}{2}\right)^r}+\psi^{(r-1)}\left(n+b+l+\frac{d}{2}\right)-\psi^{(r-1)}\left(n+l+\frac{d}{2}\right)\right]
    \end{align*}
    For the difference involving the derivatives of the digamma function we have
    \begin{align*}
        \psi^{(r-1)}\left(n+b+l+\frac{d}{2}\right)-\psi^{(r-1)}\left(n+l+\frac{d}{2}\right)&=\int_0^b\psi^{(r)}\left(n+t+l+\frac{d}{2}\right)\,dt.
    \end{align*}
    The asymptotics for the digamma function from \cite[5.15.9]{DLMF} give $\displaystyle \psi^{(r)}(x)\lesssim x^{-r}$ for every $\displaystyle x\geq 1$, hence
    $\displaystyle  \left|\frac{d^l}{dl^r}\log((h^{(b)}_{n,l})^{-\frac{1}{2}})\right|\lesssim_{b,d,r} L^{-r}$. This in turn implies $\displaystyle \left|\frac{d^l}{dl^j}(h^{(b)}_{n,l})^{-\frac{1}{2}}\right|\lesssim_{b,d,j} L^{\frac{b+1}{2}-j}N^{-\frac{b}{2}}$.

    We finally analyze the case in which the derivatives fall onto $\displaystyle C_{n,l}(S)$. We have
    \begin{align*}
        \log C_{n,l}(S)&=\frac{l}{2}\log(1-S)-\left(n+b+l+\frac{d}{2}\right)\log\left(1-\frac{S}{2}\right)\\
        &=l\left(\frac{1}{2}\log(1-S)-\log\left(1-\frac{S}{2}\right)\right)-\left(n+b+\frac{d}{2}\right)\log\left(1-\frac{S}{2}\right)
    \end{align*}
    Since $\displaystyle l\simeq L$ and $\displaystyle n+b+\frac{d}{2}\lesssim_{b,d}N$, we immediately get that 
    \begin{align*}
       \log C_{n,l}(S)\lesssim -cLS^2+CNS\lesssim -\frac{c}{2}LS^2+C'\frac{N^2}{L}. 
    \end{align*} 
    The relation $\displaystyle N\leq L^{\frac{1}{2}}$ now readily implies that $\displaystyle \log C_{n,l}(S)\lesssim-\frac{c}{2}LS^2+C'$. We also have the relation $\displaystyle \frac{d^r}{dl^j}C_{n,l}(S)=\left(\frac{1}{2}\log(1-S)-\log\left(1-\frac{S}{2}\right)\right)^jC_{n,l}(S)$. This immediately implies that
    \begin{align*}
        |(L\partial_l)^jC_{n,l}(S)|&\lesssim (LS^2)^je^{-\frac{c}{2}LS^2}\lesssim C_j.
    \end{align*}
    Now, the desired estimate for $\displaystyle (L\partial_l)^j\phi^{(b)}_{n,l}(z)$ immediately follows by the product rule. Indeed, we have
    \begin{align*}
        (L\partial_l)^j\phi^{(b)}_{n,l}(z)=\sum_{\substack{j_1,j_2,j_3\geq 0\\j_1+j_2+j_3=j}}\frac{j!}{(j_1)!(j_2)!(j_3)!}(L\partial_l)^{j_1}(h^{(b)}_{n,l})^{-\frac{1}{2}}(L\partial_l)^{j_2}C_{n,l}(S)(L\partial_l)^{j_3}\mathcal{T}_Bf_S,
    \end{align*}
    and we can then apply the estimates we have proved for each factor.

    The second estimate also immediately follows from the product rule.
\end{proof}
When studying the kernel $K_{N,L}$ of $\displaystyle \mathcal{V}_{N,L}(\tau)\mathcal{V}^\ast_{N,L}(\tau)$, we shall also need some properties of the degree operator. They are the object of the following
\begin{lemma}\label{Degree Operator bounds}
    Let $\displaystyle \Pi_l$ be the orthogonal projector on degree $l$ spherical harmonics, as well as the degree operator $\mathcal{D}$ for spherical harmonics on $\mathbb{S}^{d-1}$ given by $\displaystyle \mathcal{D}Y_{l,\nu}=lY_{l,\nu}$, where $\displaystyle Y_{l,\nu}$ is a spherical harmonic. Let $\displaystyle N,L\geq 1$, and we assume that $a$ is a smooth compactly supported function with support in $\displaystyle (c_-,c_+)$. We assume that $\|a^{(j)}\|_{L^\infty(\mathbb{R})}\leq C_jK_0N^j$ for some constant $K_0>0$ every $j$ with $\displaystyle 0\leq j\leq\max\{d-1,2\}$. Then,
    \begin{align*}
        \|a(\mathcal{D}/L)\|_{L^\infty(\mathbb{S}^{d-1})\rightarrow L^\infty(\mathbb{S}^{d-1})}+\|a(\mathcal{D}/L)\|_{L^1(\mathbb{S}^{d-1})\rightarrow L^1(\mathbb{S}^{d-1})}\leq CK_0N^{d-1},
    \end{align*}
    where the operator $\displaystyle a(\mathcal{D}/L)$ is defined as 
    \begin{align*}
        a(\mathcal{D}/L)f&=\sum_{\substack{l\geq 0}}a(l/L)\Pi_lf.
    \end{align*}
\end{lemma}
\begin{proof}
    We first prove the result when $d\geq 3$, $j,r\geq 0$. Let $\displaystyle A_j^{(r)}=\binom{j+r}{r}$ and the Ces\`aro means $\displaystyle C_j^{(r)}=\frac{1}{A_j^{(r)}}\sum_{\substack{l=0}}^jA_{j-l}^{(r)}\Pi_l$. From Theorem 2.2 in \cite{Dai-Xu}, we have the bound
    \begin{align*}
        \sup_{\substack{j\geq 0}}\|C_j^{(r)}\|_{L^\infty\rightarrow L^\infty}+\|C_j^{(r)}\|_{L^1\rightarrow L^1}<\infty.
    \end{align*}
    We choose $r=d-2$, so that $\displaystyle d-2>\frac{d-2}{2}$, and 
    \begin{align*}
        \sup_{\substack{j\geq 0}}\|C_j^{(d-2)}\|_{L^\infty\rightarrow L^\infty}+\|C_j^{(d-2)}\|_{L^1\rightarrow L^1}<\infty.
    \end{align*}
    One can prove the identity
    \begin{align*}
        a(\mathcal{D}/L)&=\sum_{\substack{j\geq 0}}\Delta_{d-1}a_jA_j^{(d-2)}C_j^{(d-2)}
    \end{align*}
    We also have the identity
    \begin{align*}
        \Delta_{d-1}a_j&=\frac{(-1)^{d-1}}{L^{d-1}}\int_{[0,1]^{d-1}}a^{(d-1)}\left(\frac{j+t_1+\dots+t_{d-1}}{L}\right)\,dt_1\,dt_2\dots \,dt_{d-1}
    \end{align*}
    This immediately implies that $\displaystyle |\Delta_{d-1}a_j|\leq \frac{\|a^{(d-1)}\|_{L^\infty}}{L^{d-1}}\leq \frac{CK_0N^{d-1}}{L^{d-1}}$. It is also clear that $\displaystyle A_j^{(d-2)}=\binom{j+d-2}{d-2}\leq C_{d}(1+j)^{d-2}$. Thus,
    \begin{align*}
\sum_{\substack{j\geq 0}}|\Delta_{d-1}a_j|A_j^{(d-2)}\leq \sum_{0\leq j\leq c_+L}\frac{CK_0N^{d-1}}{L^{d-1}}(1+j)^{d-2}\leq CK_0N^{d-1}
    \end{align*}
    Let now $\displaystyle p_\ast\in\{1,\infty\}$. We have
    \begin{align*}
        \|a(\mathcal{D}/L)\|_{L^{p_\ast}}&\leq\sum_{\substack{j\geq 0}}|\Delta_{d-1}a_j|A_j^{(d-2)}\|C_j^{(d-2)}f\|_{L^{p_\ast}}\\
        &\leq C_{d-1}\sum_{\substack{j\geq 0}}|\Delta_{d-1}a_j|A_j^{(d-2)}\|f\|_{L^{p_\ast}}\\
        &\leq C_{d-1}K_0L^{d-2}\|f\|_{L^{p_\ast}},
    \end{align*}
    as desired.
    
    We now analyze the case $d=2$. In this case, we can no longer apply the estimate for Ces\`aro means from Dai-Xu \cite{Dai-Xu}. We parametrize the circle $\mathbb{S}^1$ by $\iota\in[-\pi,\pi]$. In this case, the degree operator $\mathcal{D}$ is given by $\displaystyle \mathcal{D}=|D_\iota|$, and we have
    \begin{align*}
        a(\mathcal{D}/L)f(\iota)&=\sum_{\substack{k\in\mathbb{Z}}}a(|k|/L)\hat{f}(k)e^{ik\iota},
    \end{align*}
    where $\displaystyle \hat{f}(k)=\frac{1}{2\pi}\int_{-\pi}^\pi f(\iota)e^{-ik\iota}\,d\iota$. Let $\tilde{a}$ be the even extension of $a$ given by $\displaystyle \tilde{a}(\xi)=a(|\xi|)$ for $\displaystyle \xi\in\mathbb{R}$. It is clear that $\displaystyle \supp(\tilde{a})\subset[-c_+,-c_-]\cup[c_-,c_+]$, hence it is in particular supported away from zero and smooth. From the hypotheses, we immediately have $\displaystyle \|\tilde{a}\|_{L^1}\leq CK_0$ and $\displaystyle \|\tilde{a}''\|_{L^1}\leq CK_0N^2$. Its inverse Fourier transform is given by
    \begin{align*}
        \check{\tilde{a}}(u)=\frac{1}{2\pi}\int_{\mathbb{R}}e^{iu\xi}\tilde{a}(\xi)\,d\xi.
    \end{align*}
    We can immediately see that $\displaystyle \|\check{\tilde{a}}\|_{L^\infty}\lesssim \|\tilde{a}\|_{L^1}\leq CK_0$. When $u\neq 0$, two integrations by parts show that $\displaystyle \|\check{\tilde{a}}\|_{L^\infty}\lesssim \|\tilde{a}''\|_{L^1}|u|^{-2}\leq CK_0N^2u^{-2}$. These in turn imply that
    \begin{align*}
        \|\check{\tilde{a}}\|_{L^1}\leq CK_0\left(\int_{\substack{|u|\leq N}}1\,du+\int_{\substack{|u|> N}}N^2u^{-2}\,du\right)\lesssim CK_0N.
    \end{align*}
    Of course, the convolution kernel associated to the operator $\displaystyle a(\mathcal{D}/L)$ is given by
    \begin{align*}
        K_L(\iota)&=\frac{1}{2\pi}\sum_{\substack{k\in\mathbb{Z}}}\tilde{a}(k/L)e^{ik\iota}.
    \end{align*}
    Equivalently, this is given by the periodization
    \begin{align*}
 K_L(\iota)&=L\sum_{\substack{j\in\mathbb{Z}}}\check{\tilde{a}}(L(\iota+2\pi j))
    \end{align*}
    To see this, we note that its $k$-th Fourier coefficient is given by
    \begin{align*}
        \frac{1}{2\pi}&\int_{-\pi}^\pi L\sum_{\substack{j\in\mathbb{Z}}}\check{\tilde{a}}(L(\iota+2\pi j))e^{-ik\iota}\,d\iota=\frac{L}{2\pi}\sum_{\substack{j\in\mathbb{Z}}}\int_{\pi(2j-1)}^{\pi(2j+1)} \check{\tilde{a}}(L\iota)e^{-ik\iota}\,d\iota\\
        &=\frac{L}{2\pi}\int_{\mathbb{R}} \check{\tilde{a}}(L\iota)e^{-ik\iota}\,d\iota=\frac{1}{2\pi}\int_{\mathbb{R}} \check{\tilde{a}}(\iota)e^{-ik/L\iota}\,d\iota\\
        &=\frac{1}{2\pi}\tilde{a}(k/L).
    \end{align*}
    This implies that
    \begin{align*}
        \|K_L\|_{L^1((-\pi,\pi))}&\leq L\sum_{\substack{j\in\mathbb{Z}}}\|\check{\tilde{a}}(L(\iota+2\pi j))\|_{L^1((-\pi,\pi))}=L\sum_{\substack{j\in\mathbb{Z}}}\|\check{\tilde{a}}(L\iota)\|_{L^1((\pi(2j-1),\pi(2j+1)))}\\
        &=L\|\check{\tilde{a}}(L\cdot)\|_{L^1}=\|\check{\tilde{a}}\|_{L^1}\leq CK_0N.
    \end{align*}
    On the other hand,
    \begin{align*}
        a(\mathcal{D}/L)f(\iota)=\int_{-\pi}^\pi K_L(\iota-\omega)f(\omega)\,d\omega,
    \end{align*}
    hence Young's inequality immediately implies that
    \begin{align*}
        \|a(\mathcal{D}/L)\|_{L^\infty(\mathbb{S}^{1})\rightarrow L^\infty(\mathbb{S}^{1})}+\|a(\mathcal{D}/L)\|_{L^1(\mathbb{S}^{1})\rightarrow L^1(\mathbb{S}^{1})}\leq CK_0N.
    \end{align*}
    This finishes the proof.
\end{proof}
The main dispersive estimate that we are going to need is the object of the following
\begin{lemma}\label{Kernel Dispersive Bounds}
    Let us define $\displaystyle \alpha=\frac{1}{\kappa}-1$, and the degree operator $\mathcal{D}$ for spherical harmonics on $\mathbb{S}^{d-1}$ given by $\displaystyle \mathcal{D}Y_{l,\nu}=lY_{l,\nu}$, where $\displaystyle Y_{l,\nu}$ is a spherical harmonic. We also consider $\displaystyle \chi\in C_c^\infty((c_-,c_+))$, with $\displaystyle 0<c_-<c_+<\infty$, as well as a transformed time length $\displaystyle T_\tau$. For a physical eigenvalue $\displaystyle \lambda_{n,l}=2l+4\kappa n\left(n+l+\frac{d}{2}+\alpha\right)$, we also define $\displaystyle \varpi_{n}(l)=\sqrt{\lambda_{n,l}}$. Let $\displaystyle n\geq 0,N,L\geq 1$ with $\displaystyle n+1\sim N$ and $\displaystyle N\leq L^{\frac{1}{3}}$. In this case, we have the following dispersive estimate
    \begin{align*}
        \|e^{iu\varpi_n(\mathcal{D})}\chi\left(\mathcal{D}/L\right)\|_{L^1(\mathbb{S}^{d-1})\rightarrow L^\infty(\mathbb{S}^{d-1})}&\leq C_{T_\tau}L^{d-1}(1+|u|\sqrt{NL})^{-\frac{d-1}{2}}\,\text{, for }|u|\leq T_\tau.
    \end{align*}
    Moreover, the same estimate is true if we replace $\displaystyle \varpi_n$ by $\displaystyle \varpi_n^{(j)}=\sqrt{\lambda_{n+j,l+2-2j}}$, with $\displaystyle j=\overline{0,2}$.
\end{lemma}
\begin{proof}
    For the sake of simplicity, we denote $\displaystyle P=-\Delta_{\mathbb{S}^{d-1}}$, $\displaystyle \mathcal{A}=P^{\frac{1}{2}}$, such that for a spherical harmonic $Y_{l,\nu}$ we have $\displaystyle PY_{l,\nu}=l(l+d-2)Y_{l,\nu}$, and $\displaystyle \mathcal{A}Y_{l,\nu}=\sqrt{l(l+d-2)}Y_{l,\nu}$. We are going to apply Dinh's dispersive estimate for smooth compact manifolds without boundary, which is equation (3.8) in \cite{Dinh}. In our case, it states that for a fixed cut-off $\displaystyle \varphi\in C_c^\infty((0,\infty))$, there exist constants $t_0>0$ and $C>0$ such that
    \begin{align*}
        \|e^{i\zeta h^{-1}(h\mathcal{A})^{\frac{1}{2}}}\varphi(h^2P)\|_{L^1(\mathbb{S}^{d-1})\rightarrow L^\infty(\mathbb{S}^{d-1})}&\lesssim Ch^{-(d-1)}\left(1+\frac{|\zeta|}{h}\right)^{-\frac{d-1}{2}},
    \end{align*}
    whenever $\displaystyle h\in(0,1]$ and $\displaystyle |\zeta|\leq t_0$.

    We take our semiclassical parameter to be $\displaystyle h=L^{-1}$, $\displaystyle \zeta=\frac{s}{\sqrt{L}}$, and we therefore deduce that for every $\tilde{\chi}\in C_c^\infty((0,\infty))$, we have
    \begin{align*}
         \|e^{is(\mathcal{A})^{\frac{1}{2}}}\tilde{\chi}(\mathcal{A}/L)\|_{L^1(\mathbb{S}^{d-1})\rightarrow L^\infty(\mathbb{S}^{d-1})}&\lesssim CL^{d-1}\left(1+|s|\sqrt{L}\right)^{-\frac{d-1}{2}}
    \end{align*}
    for $\displaystyle |s|\leq t_0\sqrt{L}$.

    Let $\displaystyle A_n=2+4\kappa n$ and $\displaystyle B_n=4\kappa n\left(n+\frac{d}{2}+\alpha\right)$, so that $\displaystyle \lambda_{n,l}=lA_n+B_n$. Since $\displaystyle n+1\sim N$, we have $\displaystyle A_n\sim N$ and $\displaystyle |B_n|\leq CN^2$. For $\displaystyle j=\overline{0,2}$, we define $\displaystyle A_{n,j}=A_{n+j}$ and $\displaystyle B_{n,j}=B_{n+j}+(2-2j)A_{n+j}$. It is clear that we still have $\displaystyle A_{n,j}\sim N$ and $\displaystyle |B_{n,j}|\leq CN^2$. For the purpose of our estimates, we can therefore consider $\displaystyle \varpi_n(l)=\sqrt{Al+B}$, where $\displaystyle A\sim N$ and $\displaystyle |B|\leq CN^2$.

    We shall now compare the two phases $\displaystyle \varphi_n(l)$ with $\displaystyle \sqrt{A}(l(l+d-2))^{\frac{1}{4}}$. We define $\displaystyle \theta=\frac{l}{L}$, and the difference 
    \begin{align*}
    \mathcal{R}_{n,L}(\theta)&=\sqrt{AL\theta+B}-\sqrt{A}(L\theta(L\theta+d-2))^{\frac{1}{4}}=\sqrt{AL}\left(\sqrt{\theta+\frac{B}{AL}}-\left(\theta\left(\theta+\frac{d-2}{L}\right)\right)^{\frac{1}{4}}\right).
    \end{align*}
    Since $\displaystyle \supp\tilde{\chi}\Subset(0,\infty)$ and $\displaystyle \frac{B}{LA}\leq\frac{N}{L}\leq L^{-\frac{2}{3}}$, there exists $c>0$ such that $\theta\geq c$. Upon multiplying the expression by its conjugate, we can immediately see that for every $\displaystyle j\geq 0$, $\displaystyle \left|\partial_\theta^j\left(\sqrt{\theta+\frac{B}{AL}}-\left(\theta\left(\theta+\frac{d-2}{L}\right)\right)^{\frac{1}{4}}\right)\right|\lesssim_j \frac{N}{L}$, which immediately implies that $\displaystyle \left|\partial_\theta^j\mathcal{R}_{n,L}(\theta)\right|\lesssim_j \frac{N^{\frac{3}{2}}}{L^{\frac{1}{2}}}$. For every $\displaystyle u\in[-T_\tau,T_\tau]$, we define $\displaystyle a_{u,n,l}(\theta)=\chi(\theta)e^{iu\mathcal{R}_{n,l}(\theta)}$. Since $\displaystyle N\leq L^{\frac{1}{3}}$, it thus follows that for every $\displaystyle j\geq 0$, $\displaystyle \sup_{\substack{|u|\leq T_\tau}}\left\|\partial_\theta^ja_{u,n,l}(\theta)\right\|_{L^\infty}\leq C_{j,T_\tau,\chi}$. From Lemma \ref{Degree Operator bounds}, it follows that $\displaystyle \sup_{\substack{|u|\leq T_\tau}}\|a_{u,n,L}(\mathcal{D}/L)\|_{L^\infty\rightarrow L^\infty}\leq C_{T_\tau}$. Let $\displaystyle s=\sqrt{A}u$. We now claim that
    \begin{align*}
        e^{iu\varpi_n(\mathcal{D})}\chi\left(\mathcal{D}/L\right)&=a_{u,n,L}(\mathcal{D}/L)e^{iu\sqrt{A}(\mathcal{A})^{\frac{1}{2}}}\tilde{\chi}(\mathcal{A}/L).
    \end{align*}
    Indeed, on degree $l$, we note that $\displaystyle\frac{\mathcal{A}}{L}$ has eigenvalue $\displaystyle \frac{\sqrt{l(l+d-2)}}{L}$. When $\displaystyle \chi(l/L)=0$, $a_{u,n,L}(l/L)=0$, hence the equality is true. When $\displaystyle \chi(l/L)\neq 0$, it is clear that $\displaystyle c_-\leq l/L\leq c_+$, and in particular $\displaystyle c-_\leq \frac{\sqrt{l(l+d-2)}}{L}\leq \sqrt{c_+\left(c_++\frac{d-2}{L}\right)}$. As long as we choose $\tilde{\chi}$ so that it equals $1$ on $\displaystyle \left[c_-,\sqrt{c_+\left(c_++\frac{d-2}{L}\right)}\right]$, the desired identity immediately follows from the definition of $\mathcal{R}_{n,L}$ and the fact that $\displaystyle a_{u,n,l}(\theta)=\chi(\theta)e^{iu\mathcal{R}_{n,l}(\theta)}$. As long as $\displaystyle |s|\leq t_0\sqrt{L}$, the conclusion follows. This is equivalent to $\displaystyle T_\tau\frac{\sqrt{N}}{\sqrt{L}}\leq t_0$, which can be guaranteed as long as $L$ is sufficiently large, given that $\displaystyle T_\tau\frac{\sqrt{N}}{\sqrt{L}}\leq \frac{T_\tau}{L^{\frac{1}{3}}}$. Thus, there exists a threshold $L_\ast$ depending on $T_\tau$, $d$, $\kappa$, and $\chi$, for which the argument above is valid for $L\geq L_\ast$.

    We finally analyze the case $L\leq L_\ast$. Our operator has kernel
    \begin{align*}
        K_{n,L}(u,\omega,\overline{\omega})&=\sum_{\substack{l,\nu}}e^{iu\varpi_n(l)}\chi(l/L)Y_{l,\nu}(\omega)\overline{Y_{l,\nu}(\overline{\omega})}.
    \end{align*}
    By the spherical addition formula and the Cauchy-Buniakowski-Schwartz inequality,
     \begin{align*}
        |K_{n,L}(u,\omega,\overline{\omega})|&\lesssim C_{\chi}L^{d-1}.
    \end{align*}
    We have that $\displaystyle (1+|u|\sqrt{NL})^{-\frac{d-1}{2}}\geq (1+T_\tau L_\ast^{\frac{2}{3}})^{-\frac{d-1}{2}}$, which immediately implies that
     \begin{align*}
        |K_{n,L}(u,\omega,\overline{\omega})|&\lesssim C_{T_\tau}L^{d-1}(1+|u|\sqrt{NL})^{-\frac{d-1}{2}},
    \end{align*}
    as desired. This finishes the proof.
\end{proof}
We recall that our kernel is given by the formula
\begin{align*}
    K_{N,L}(u;y,\overline{y})&=\chi_{\rm bd}(z)\chi_{\rm bd}(\overline{z})\sum_{\substack{n+1\sim N\\l,\nu}}e^{iu\varpi_n^{(a)}(l)}\chi^2(l/L)\phi^{(b)}_{n,l}(z)\phi^{(b)}_{n,l}(\overline{z})Y_{l,\nu}(\omega)\overline{Y_{l,\nu}(\overline{\omega})}.
\end{align*}
By Lemma \ref{Uniform Lemma}, we have the bounds
\begin{align*}
    |\phi^{(b)}_{n,l}(z)\phi^{(b)}_{n,l}(\overline{z})|&\lesssim L^{b+1}N^b.
\end{align*}
The spherical addition formula and the Cauchy-Buniakowski-Schwartz inequality now imply that
\begin{align*}
    \sum_{\substack{l+1\sim L}}\sum_{\substack{\nu}}|Y_{l,\nu}(\omega)||Y_{l,\nu}(\overline{\omega})|&\lesssim L^{d-1}.
\end{align*}
Since there are $O(N)$ radial indices $n$ with the property $n+1\sim N$, we get that
\begin{align*}
    | K_{N,L}(u;y,\overline{y})|&\lesssim L^{d+b}N^{b+1}.
\end{align*}
In what follows, we show that when $|u|\leq T_\tau$ and $\displaystyle N\leq L^{\frac{1}{3}}$,
\begin{align*}
    |K_{N,L}(u;y,\overline{y})|&\lesssim_{T_\tau}L^{d+b}N^{b+d}(1+|u|\sqrt{NL})^{-\frac{d-1}{2}}
\end{align*}
To see that, we fix $n$ with $n+1\sim N$, and $z,\overline{z}\in\left[\frac{1}{2},1\right]$, and define
\begin{align*}
    k_{n,L}(u;\omega,\overline{\omega};z,\overline{z})&=\sum_{\substack{l,\nu}}e^{iu\varpi_n^{(a)}(l)}\chi^2(l/L)\phi^{(b)}_{n,l}(z)\phi^{(b)}_{n,l}(\overline{z})Y_{l,\nu}(\omega)\overline{Y_{l,\nu}(\overline{\omega})}
\end{align*}
In this case,
\begin{align*}
    K_{N,L}(u;y,\overline{y})&=\chi_{\rm bd}(z)\chi_{\rm bd}(\overline{z})\sum_{\substack{n+1\sim N}}k_{n,L}(u;\omega,\overline{\omega};z,\overline{z}).
\end{align*}
We now define $\displaystyle a_{n,z,\overline{z}}(\theta)=\chi(\theta)^2\phi^{(b)}_{n,L\theta}(z)\phi^{(b)}_{n,L\theta}(\overline{z})$, as well as the operator $\displaystyle \mathcal{M}_{n;z,\overline{z}}=a_{n,z,\overline{z}}(\mathcal{D}/L)$. The latter is given by
\begin{align*}
 \mathcal{M}_{n;z,\overline{z}}f&=\sum_{\substack{l\geq 0}}\chi^2(l/L) \phi^{(b)}_{n,l}(z)\phi^{(b)}_{n,l}(\overline{z})\Pi_lf  
\end{align*}
By Lemma \ref{Uniform Derivative bounds lemma}, we have
\begin{align*}
    |(L\partial_l)^j\phi^{(b)}_{n,l}(z)|+|(L\partial_l)^j\phi^{(b)}_{n,l}(\overline{z})|&\lesssim_j L^{\frac{b+1}{2}}N^{\frac{b}{2}+j}.
\end{align*}
By the product rule, it immediately follows that
\begin{align*}
    \|(L\partial_l)^ja_{n,z,\overline{z}}\|_{L^\infty}&\lesssim L^{b+1}N^{b+j}.
\end{align*}
By Lemma \ref{Degree Operator bounds}, we have
\begin{align*}
   \|\mathcal{M}_{n;z,\overline{z}}\|_{L^\infty(\mathbb{S}^{d-1})\rightarrow L^\infty(\mathbb{S}^{d-1})}&\lesssim L^{b+1}N^{b+d-1}.
\end{align*}
Let $M_{n;z,\overline{z}}$ be the kernel of $\mathcal{M}_{n;z,\overline{z}}$. Since $\displaystyle \|\mathcal{M}_{n;z,\overline{z}}\|_{L^\infty(\mathbb{S}^{d-1})\rightarrow L^\infty(\mathbb{S}^{d-1})}=\sup_{\substack{\omega}}\int_{\mathbb{S}^{d-1}}|M_{n;z,\overline{z}}(\omega,\eta)|\,d\eta$, it in particular follows that
\begin{align*}
    \sup_{\substack{\omega}}\int_{\mathbb{S}^{d-1}}|M_{n;z,\overline{z}}(\omega,\eta)|\,d\eta&\lesssim L^{b+1}N^{b+d-1}.
\end{align*}
We now consider a smooth compactly supported cut-off function $\displaystyle\chi_1$ which is identically equal to $1$ on $\displaystyle \supp\chi$. Let us define the operator
\begin{align*}
    \mathcal{S}_{n,L}(u)&=e^{iu\overline{\varpi}^{(a)}_n(\mathcal{D})}\chi_1(\mathcal{D}/L).
\end{align*}
By Lemma \ref{Kernel Dispersive Bounds}, we have 
\begin{align*}
    \|\mathcal{S}_{n,L}(u)\|_{L^1(\mathbb{S}^{d-1})\rightarrow L^\infty(\mathbb{S}^{d-1})}&\lesssim_{T_\tau}L^{d-1}(1+|u|\sqrt{NL})^{-\frac{d-1}{2}}.
\end{align*}
Let $\displaystyle S_{n,L}(u;\eta,\overline{\omega})$ be the kernel of $\mathcal{S}_{n,L}$.
Since $\displaystyle \|\mathcal{S}_{n,L}(u)\|_{L^1(\mathbb{S}^{d-1})\rightarrow L^\infty(\mathbb{S}^{d-1})}=\sup_{\eta,\overline{\omega}}|S_{n,L}(u;\eta,\overline{\omega})|$, it follows that
\begin{align*}
    \sup_{\eta,\overline{\omega}}|S_{n,L}(u;\eta,\overline{\omega})|&\lesssim _{T_\tau}L^{d-1}(1+|u|\sqrt{NL})^{-\frac{d-1}{2}}. 
\end{align*}
We can also check that $\displaystyle k_{n,L}(u;\omega,\overline{\omega};z,\overline{z})$ is the kernel of $\displaystyle\mathcal{M}_{n;z,\overline{z}}S_{n,L}(u)$, due to the support conditions on $\chi$ and $\chi_1$. Thus,
\begin{align*}
    k_{n,L}(u;\omega,\overline{\omega};z,\overline{z})&=\int_{\mathbb{S}^{d-1}}M_{n;z,\overline{z}}(\omega,\eta)S_{n,L}(u;\eta,\overline{\omega})\,d\eta,
\end{align*}
from which we immediately get the bound
\begin{align*}
    |k_{n,L}(u;\omega,\overline{\omega};z,\overline{z})|&\leq \sup_{\substack{\omega}}\int_{\mathbb{S}^{d-1}}|M_{n;z,\overline{z}}(\omega,\eta)|\,d\eta\sup_{\eta,\overline{\omega}}|S_{n,L}(u;\eta,\overline{\omega})|\\
    &\lesssim L^{b+d}N^{b+d-1}(1+|u|\sqrt{NL})^{-\frac{d-1}{2}}.
\end{align*}
Since there are $O(N)$ radial terms with $n+1\sim N$, it follows that
\begin{align*}
    |K_{N,L}(u;y,\overline{y})|&\lesssim L^{b+d}N^{b+d}(1+|u|\sqrt{NL})^{-\frac{d-1}{2}}.
\end{align*}
By interpolating the two bounds for $K_{N,L}$, we get that
\begin{align*}
    |K_{N,L}(u;y,\overline{y})|\lesssim L^{b+d}N^{b+1}\left(1+\frac{|u|\sqrt{NL}}{N^2}\right)^{-\frac{d-1}{2}},
\end{align*}
uniformly for $L$, $N$, and $u$, with $N\leq L^{\frac{1}{3}}$ and $|u|\leq T_\tau$.
It then follows that
\begin{align*}
    \|K_{N,L}\|_{L_u^{\frac{q}{2}}}&\lesssim L^{d+b-\frac{1}{q}}N^{b+1+\frac{3}{q}}
\end{align*}
For $\displaystyle G\in l^2$, we also define the reconstruction operator $\displaystyle \mathfrak{V}_{N,L}G(\tau)=\mathcal{V}_{N,L}(\tau)G$. By using a smooth cutoff around $l=0$, we can also define this operator for the lowest angular block. We can also immediately see that its norm is controlled by the usual static estimate that we shall derive momentarily. We also note that if we enlarge the dyadic supports by fixed constant factors, we can also allow for fixed shifts in the radial and angular indices.

For an arbitrary test function $F$ we then have
\begin{align*}
 \|\mathfrak{V}_{N,L}\mathfrak{V}_{N,L}^\ast F(\tau)\|_{L_{y}^\infty}&\lesssim \int |K_{N,L}(\tau-\tau')\|F(\tau')\|_{L_y^1}\,d\tau',   
\end{align*}
By Young's inequality and the previous kernel bound we also have

\begin{align*}
    \|\mathfrak{V}_{N,L}\mathfrak{V}_{N,L}^\ast F\|_{L_\tau^qL_y^\infty}&\lesssim L^{d+b-\frac{1}{q}}N^{b+1+\frac{3}{q}}\|F\|_{L_\tau^{q'}L_y^1}.
\end{align*}

We also claim that $\displaystyle \|\mathfrak{V}_{N,L}(\tau)\mathfrak{V}^\ast_{N,L}(\tau')\|_{L^2(\,dy)\rightarrow L^2(\,dy)}\lesssim W^b_{N,L}$. For $G\in l^2$, we set for brevity $\displaystyle F(y)=\mathcal{V}_{N,L}(\tau)G(y)$. We know that the function $\Phi^{(b)}_{n,l,\nu}$ form an orthogonal basis of $\displaystyle L^2(R^b\,dy)$, hence
\begin{align*}
    \|F\|^2_{L^2(R^b\,dy)}=\int_{\mathbb{B}}R^b|F|^2\,dy\lesssim \|G\|^2_{l^2}.
\end{align*}
When $R\geq W_{N,L}^{-1}$, this implies the bound
\begin{align*}
    \int_{R\geq W_{N,L}^{-1}}|F(y)|^2\,dy\leq W^b_{N,L}\int_{\mathbb{B}}R^b|F|^2\,dy\lesssim W^b_{N,L}\|G\|^2_{l^2}.
\end{align*}
We now analyze the case $\displaystyle R<W^{-1}_{N,L}$. By Lemma \ref{Uniform Lemma}, when $l+1\sim L$ and $z\geq \frac{1}{2}$, we have
\begin{align*}
    \sum_{\substack{n+1\sim N}}|\phi^{(b)}_{n,l}(z)|^2\lesssim \sum_{\substack{n+1\sim N}}(n+1)^b(n+l+1)^{b+1}\lesssim N^{b+1}(N+L)^{b+1}=W^{b+1}_{N,L}.
\end{align*}
By the orthogonality of the spherical harmonics, for fixed $z$ we have
\begin{align*}
    \int_{\mathbb{S}^{d-1}}|F(z\omega)|^2\,d\omega&=|\chi_{\rm bd}(z)|^2\sum_{\substack{l,\nu}}|\chi(l/L)|^2\left|\sum_{\substack{n+1\sim N}}e^{it\varpi^{(a)}_n(l)}\phi^{(b)}_{n,l}(z)G_{n,l,\nu}\right|^2\\
    &\lesssim \sum_{\substack{l,\nu}}\left(\sum_{\substack{n+1\sim N}}|\phi^{(b)}_{n,l}(z)|^2\right)\left(\sum_{\substack{n+1\sim N}}|G_{n,l,\nu}|^2\right)\\
    &\lesssim W^{b+1}_{N,L}\|G\|^2_{l^2}.
\end{align*}
Thus,
\begin{align*}
    \int_{R<W^{-1}_{N,L}}|F(y)|^2\,dy&=\int_{\sqrt{1-W^{-1}_{N,L}}}^1z^{d-1}\int_{\mathbb{S}^{d-1}}|F(z\omega)|^2\,d\omega\,dz\\
    &\lesssim W^{b+1}_{N,L}\|G\|^2_{l^2}\int_{\sqrt{1-W^{-1}_{N,L}}}^1z^{d-1}\,dz\\
    &\lesssim W^{b}_{N,L}\|G\|^2_{l^2}
\end{align*}
By combining the two bounds, we deduce that $\displaystyle \|F\|_{L^2(\,dy)}\lesssim W^{\frac{b}{2}}_{N,L}\|G\|_{l^2}$, hence 
\begin{align*}
\|\mathcal{V}_{N,L}(\tau)\|_{l^2\rightarrow L^2(\,dy)}\lesssim W^{\frac{b}{2}}_{N,L}.
\end{align*}
By duality, we get that $\displaystyle \|\mathcal{V}^\ast_{N,L}(\tau')\|_{L^2(\,dy)\rightarrow l^2}\lesssim W^{\frac{b}{2}}_{N,L}$, hence 
\begin{align*} 
\|\mathcal{V}_{N,L}(\tau)\mathcal{V}^\ast_{N,L}(\tau')\|_{L^2(\,dy)\rightarrow L^2(\,dy)}\lesssim W^{b}_{N,L}.
\end{align*}

Now, since for an arbitrary branch we have $\displaystyle\|G^{(a)}\|_{l^2}\lesssim  \Lambda_{N,L}\mathcal{E}_{N,L}$, and in particular, when $N\leq L^{\frac{1}{3}}$, $\displaystyle \Lambda_{N,L}\simeq\sqrt{NL}$, interpolating with the previous kernel bound gives
\begin{align*}
    \|\mathcal{V}_{N,L}(\tau)\mathcal{V}_{N,L}^\ast (\tau')\|_{L^{r'}\rightarrow L^r}&\lesssim L^{b+d\left(1-\frac{2}{r}\right)}N^{b+1-\frac{2}{r}}\left(1+\frac{|\tau-\tau'|\sqrt{NL}}{N^2}\right)^{-\frac{d-1}{2}\left(1-\frac{2}{r}\right)}.
\end{align*}
If $(q,r)$ is not sharp admissible, for a test function $F$, this implies that
\begin{align*}
    \|\mathfrak{V}\mathfrak{V}^\ast F\|_{L_\tau^qL_y^r}&\lesssim L^{b+d\left(1-\frac{2}{r}\right)-\frac{1}{q}}N^{b+1-\frac{2}{r}+\frac{3}{q}}\|F\|_{L_\tau^{q'}L_y^{r'}}
\end{align*}
Since $\displaystyle \langle \mathfrak{V}_{N,L}\mathfrak{V}_{N,L}^\ast F,F\rangle=\|\mathfrak{V}_{N,L}^\ast F\|^2_{l^2}$, it follows that
\begin{align*}
    \|\mathcal{V}_{N,L}(\cdot)G\|_{L_\tau^qL_y^r}&\lesssim L^{\frac{b}{2}+\frac{d}{2}\left(1-\frac{2}{r}\right)-\frac{1}{2q}}N^{\frac{b+1}{2}+\frac{3}{2q}-\frac{1}{r}}\|G\|_{l^2}.
\end{align*}
This in turn implies that
\begin{align*}
    \|\chi_{\rm bd}T_{ij}V_{N,L}\|_{L_\tau^qL_y^\infty}&\lesssim L^{\frac{b+d+1}{2}-\frac{d}{r}-\frac{1}{2q}}N^{\frac{b+2}{2}+\frac{3}{2q}-\frac{1}{r}}\mathcal{E}_{N,L}.
\end{align*}
We now recall that $\displaystyle s_{q,r}=\frac{d+b+1}{2}-\frac{d}{r}-\frac{1}{2q}$ and $\displaystyle \frac{b+2}{2}+\frac{3}{2q}-\frac{1}{r}=s_{q,r}-\delta_{q,r}$, and since in this regime $W_{N,L}\simeq NL$, we get that
\begin{align*}
    \|\chi_{\rm bd}\nabla^2V_{N,L}\|_{L_\tau^qL_y^\infty}&\lesssim_{T_\tau}W_{N,L}^{s_{q,r}}N^{-\delta_{q,r}}\mathcal{E}_{N,L}
\end{align*}
for every $N\leq L^{\frac{1}{3}}$. When $(q,r)$ is sharp admissible, the previous bounds include an additional $\displaystyle \left[\log\left(2+\frac{T_\tau}{N^{\frac{3}{2}}L^{-\frac{1}{2}}}\right)\right]$, which we can absorb into the $W_{N,L}^s$ factor, since $s>s_{q,r}$.

We now analyze the case $N\geq L^{\frac{1}{3}}$. By Lemma \ref{Uniform Lemma}, for every $z\geq\frac{1}{2}$, we have
\begin{align*}
    \sum_{\substack{n\geq 0\\1+\lambda^{(b)}_{n,l}\lesssim\Lambda^2_{N,L}}}|\phi^{(b)}_{n,l}(z)|^2&\lesssim\Lambda^{2b+2}_{N,L}.
\end{align*}
This implies that
\begin{align*}
    \sum_{\substack{l+1\sim L,\nu\\1+\lambda^{(b)}_{n,l}\lesssim\Lambda^2_{N,L}}}|\Phi^{(b)}_{n,l}(z)|^2&\lesssim\Lambda^{2b+2}_{N,L}L^{d-1}.
\end{align*}
By Cauchy-Buniakowski-Schwartz's inequality, for $G\in l^2$ we have
\begin{align*}
    \|\mathcal{V}_{N,L}(\tau)G\|_{L_y^\infty}&\lesssim L^{\frac{d-1}{2}}\Lambda^{b+1}_{N,L}\|G\|_{l^2}.
\end{align*}
We also recall that $\displaystyle \|\mathcal{V}_{N,L}(\tau)G\|_{L_y^2}\lesssim \Lambda^b_{N,L}\|G\|_{l^2}$. By interpolating these two bounds, we get that
\begin{align*}
    \|\mathcal{V}_{N,L}(\tau)G\|_{L_y^r}&\lesssim (\Lambda^b_{N,L})^{\frac{2}{r}}(L^{\frac{d-1}{2}}\Lambda^{b+1}_{N,L})^{1-\frac{2}{r}}\|G\|_{l^2}.
\end{align*}

Since $\displaystyle\|G^{(a)}\|_{l^2}\lesssim  \Lambda_{N,L}\mathcal{E}_{N,L}$, upon integrating on a time interval we have

   \begin{align*}
\|\chi_{\rm bd}\nabla^2V_{N,L}\|_{L_\tau^qL_y^r}&\lesssim T_\tau^{\frac{1}{q}}L^{\frac{d-1}{2}\left(1-\frac{2}{r}\right)}\Lambda^{b+2-\frac{2}{r}}_{N,L}\mathcal{E}_{N,L}.
\end{align*} 
However, $W_{N,L}=\Lambda^2_{N,L}$, and $\displaystyle L^{\frac{d-1}{2}\left(1-\frac{2}{r}\right)}\Lambda^{b+2-\frac{2}{r}}_{N,L}=W^{s_{q,r}}_{N,L}L^{\frac{d-1}{2}\left(1-\frac{2}{r}\right)}\Lambda_{N,L}^{-(d-1)\left(1-\frac{2}{r}\right)+\frac{1}{q}}$. Whenever $N\geq L^{\frac{1}{3}}$, we can immediately show that $L^{\frac{d-1}{2}\left(1-\frac{2}{r}\right)}\Lambda_{N,L}^{-(d-1)\left(1-\frac{2}{r}\right)+\frac{1}{q}}\leq N^{-\delta_{q,r}}$. We can then infer that

\begin{align*}
    \|\chi_{\rm bd}\nabla^2V_{N,L}\|_{L_\tau^qL_y^r}&\lesssim_{T_\tau}W_{N,L}^{s_{q,r}}N^{-\delta_{q,r}}\mathcal{E}_{N,L}.
\end{align*}

We now analyze the interior of the ball. Here, the operators $\mathcal{L}_\gamma$ are uniformly elliptic. Let $\chi_{\rm int}$ be a smooth radial cut-off function supported in $\left\{z||z|\leq \frac{3}{4}\right\}$. Now, elliptic estimates for $\mathcal{L}_\gamma$, along with the spectral properties of $V_{N,L}$ imply that
\begin{align*}
    \|\chi_{\rm int}V_{N,L}(\tau)\|_{H^\sigma}&\lesssim \Lambda^\sigma_{N,L}\|V_{N,L}\|_{L^2(R^\alpha)}.
\end{align*}
By combining these with Sobolev embeddings for $\sigma=\frac{d}{2}\left(1-\frac{2}{r}\right)+2+\varepsilon$, we get that
\begin{align*}
    \|\chi_{\rm int}\nabla^2V_{N,L}\|_{L_\tau^qL_y^r}&\lesssim \Lambda^\sigma_{N,L}\|V_{N,L}\|_{L^2(R^\alpha)}.
\end{align*}
Since $\displaystyle \|V_{N,L}\|_{L^2(R^\alpha)}\lesssim \Lambda_{N,L}^{-1}\mathcal{E}_{N,L}$, it follows that
\begin{align*}
    \|\chi_{\rm int}\nabla^2V_{N,L}\|_{L_\tau^qL_y^r}&\lesssim \Lambda^{\frac{d}{2}\left(1-\frac{2}{r}\right)+2+\varepsilon}_{N,L}\mathcal{E}_{N,L}.
\end{align*}
Since $\displaystyle N\leq \Lambda_{N,L}$, if we choose $\varepsilon<b+\frac{1}{2}+\frac{1}{q}-\frac{1}{r}=2s_q-\delta_q$, we deduce the inequality $\Lambda^{\frac{d}{2}\left(1-\frac{2}{r}\right)+2+\varepsilon}_{N,L}\leq W^{2s_{q,r}}N^{-\delta_{q,r}}$, which implies the estimate
\begin{align*}
    \|\chi_{\rm int}\nabla^2V_{N,L}\|_{L_\tau^qL_y^r}&\lesssim_{T_\tau}W_{N,L}^{s_{q,r}}N^{-\delta_{q,r}}\mathcal{E}_{N,L}.
\end{align*}
Combining this with the boundary estimate yields
\begin{align*}
    \|\nabla^2V_{N,L}\|_{L_\tau^qL_y^r}&\lesssim_{T_\tau}W_{N,L}^{s_{q,r}}N^{-\delta_{q,r}}\mathcal{E}_{N,L}.
\end{align*}

Our last ingredient consists of a comparison between the dyadic block energies and the global $\mathcal{H}^{2s}$ norm of the initial data. This is the object of the following
\begin{lemma}
    When $N,L\geq 1$ and $s\geq 0$, we have
    \begin{align*}
        \sum_{\substack{N,L}}W^{2s}_{N,L}\mathcal{E}^2_{N,L}&\lesssim \|(V_\tau,\nabla_y V)(0)\|^2_{\mathcal{H}^{2s}_R}.
    \end{align*}
\end{lemma}

\begin{proof}
    When $n+1\sim N$ and $l+1\sim L$, we have $1+\lambda_{n,l}\sim W_{N,L}$. When $s=m\geq 0$ is an integer, we have
    \begin{align*}
    \sum_{\substack{N,L}}W^{2m}_{N,L}\mathcal{E}^2_{N,L}&\lesssim \sum_{\substack{N,L}}\sum_{\substack{n+1\sim N,l+1\sim L}}(1+\lambda_{n,l})^{2m}(|g_{n,l}|^2+\lambda_{n,l}|f_{n,l}|^2)\\
    &\simeq_m\sum_{\substack{k=0}}^m\|((\mathcal{L}_{\frac{1}{\kappa}-1})^kg,\nabla_y(\mathcal{L}_{\frac{1}{\kappa}-1})^kf)\|^2_{\mathcal{H}_R}\\
    &\lesssim_{d,m,\kappa} \|(g,\nabla_yf)\|^2_{\mathcal{H}^{2m}_R}.
    \end{align*}
    The last estimate follows from the definition of the weighted spaces $\mathcal{H}^{2k}$ and the elliptic estimates from Ifrim-Tataru \cite{Ifrim-Tataru} (see section 2.2 and Theorem 7 in section 7 in \cite{Ifrim-Tataru}).

    We now treat the noninteger case. We consider the fixed coefficient map
    \begin{align*}
      (g,w)\rightarrow\left(\langle g,\Phi_j\rangle_{L^2(R^{\frac{1}{\kappa}-1}\,dy},\frac{\kappa}{\sqrt{\lambda_j}} \langle g,\nabla\Phi_j\rangle_{L^2(R^{\frac{1}{\kappa}}\,dy)}\right)_j  
    \end{align*}
    We note that $\mathcal{L}_{\frac{1}{\kappa}-1}=(\nabla)^\ast\nabla$, where the adjoint is taken with respect to the inner product in the Hilbert space $L^2(\kappa R^{\frac{1}{\kappa}}\,dy)$. As we already have the estimates in the corresponding Hilbert spaces when $s=m$ is an integer, the desired estimate in the noninteger case directly follows from interpolation.
\end{proof}
We can now finish the proof of Theorem \ref{Strichartz Estimates} in the case $Q(\tau)=0$. Let $s>s_{q,r}$, and set $\epsilon=s-s_{q,r}>0$. By Cauchy-Buniakowski-Schwartz's inequality we have
\begin{align*}
    \|\nabla_y^2 V\|_{L_\tau^qL_y^r}&\lesssim_{T_\tau}\sum_{\substack{N,L}} \|\nabla_y^2 V_{N,L}\|_{L_\tau^qL_y^r}\lesssim_{T_\tau} \sum_{\substack{N,L}} W_{N,L}^{s_{q,r}}N^{-\delta_{q,r}}\mathcal{E}_{N,L}\\
    &\lesssim_{T_\tau} \left(\sum_{\substack{N,L}}W_{N,L}^{-2\epsilon}N^{-2\delta_{q,r}}\right)^{\frac{1}{2}}\left(\sum_{\substack{N,L}}W_{N,L}^{2s}\mathcal{E}^2_{N,L}\right)^{\frac{1}{2}}\\
    &\lesssim_{T_\tau} \left(\sum_{\substack{N,L}}W_{N,L}^{-2\epsilon}N^{-2\delta_{q,r}}\right)^{\frac{1}{2}}\|(V_\tau,\nabla_y V)(0)\|_{\mathcal{H}^{2s}_R}
\end{align*}
Since $\displaystyle W_{N,L}=N(N+L)\geq NL$, we deduce that
\begin{align*}
    \|\nabla_y^2 V\|_{L_\tau^qL_y^r}&\lesssim_{T_\tau} \left(\sum_{\substack{N,L}}L^{-2\epsilon}N^{-2\epsilon-2\delta_{q,r}}\right)^{\frac{1}{2}}\|(V_\tau,\nabla_y V)(0)\|_{\mathcal{H}^{2s}_R}\\
    &\lesssim_{T_\tau,q,s}\|(V_\tau,\nabla_y V)(0)\|_{\mathcal{H}^{2s}_R},
\end{align*}
which is the desired estimate.

The special case also follows immediately, since we may choose $q_1>4$ sufficiently close to $4$ such that $s>s_{q_1,\infty}$, and we obtain
\begin{align*}
  \|\nabla_y^2 V\|_{L_\tau^4L_y^\infty}\lesssim_{T_\tau,s}  T^{\frac{1}{4}-\frac{1}{q_1}}\|\nabla_y^2 V\|_{L_\tau^{q_1}L_y^\infty}\lesssim_{T_\tau,s}\|(V_\tau,\nabla_y V)(0)\|_{\mathcal{H}^{2s}_R}. 
\end{align*}
This finishes the proof in the case $\displaystyle Q(\tau)=0$. We now extend the result to general $Q\in L_\tau^1$. Our equation is
\begin{align*}
    \partial_\tau^2V+\mathcal{L}_\alpha V+Q(\tau)V&=0
\end{align*}

By Duhamel's formula, we can write
\begin{align*}
    V(\tau)&=\cos(\tau\sqrt{\mathcal{L}_\alpha})V(0)+\sqrt{\mathcal{L}_\alpha}^{-1}\sin(\tau\sqrt{\mathcal{L}_\alpha})V_\tau(0)\\
    &=-\int_0^\tau \sqrt{\mathcal{L}_\alpha}^{-1}\sin((\tau-\sigma)\sqrt{\mathcal{L}_\alpha})(Q(\sigma)V(\sigma))\,d\sigma:=V^1-V^2
\end{align*}
From the previous case, we get that
\begin{align*}
    \|\nabla_y^2 V^1\|_{L_\tau^qL_y^r}&\lesssim_{T_\tau,s}\|(V_\tau,\nabla_y V)(0)\|_{\mathcal{H}^{2s}_R}
\end{align*}
For $V^2$, we get that
\begin{align*}
  \|\nabla_y^2 V^2\|_{L_\tau^qL_y^r}&\lesssim \int_0^\tau \|\nabla_y^2 (\sqrt{\mathcal{L}_\alpha}^{-1}\sin((\tau-\sigma)\sqrt{\mathcal{L}_\alpha})(Q(\sigma)V(\sigma)))\|_{L_\tau^qL_y^\infty}\,d\sigma\\
  &\lesssim_{T_\tau,s}\int_0^{T_\tau}\|Q(\sigma)V(\sigma)\|_{\mathcal{H}^{2s}_R}\,d\sigma \\
  &\lesssim_{T_\tau,s}\int_0^{T_\tau}|Q(\sigma)|\|(\mathcal{L}_\alpha)^sV(\sigma)\|_{\mathcal{H}_R}\,d\sigma\\
  &\lesssim_{T_\tau,s}\|Q\|_{L_\tau^1}\|(V_\tau,\nabla_y V)(0)\|_{\mathcal{H}^{2s}_R}
\end{align*}
Here, the last step follows by the energy estimate in Proposition 3.1 in \cite{Ifrim-Tataru}, and as well as the elliptic estimates and the definition of weighted $\mathcal{H}^{2k}$ spaces (see section 2.2 and Theorem 7 in section 7 in \cite{Ifrim-Tataru}), since in this case, $\mathcal{L}_\alpha$ commutes with the equation.

This finishes the proof.
\bibliography{compressible_euler_strichartz.bib}
\bibliographystyle{plain}
\end{document}